\documentclass{article}
\usepackage{amsfonts}
\usepackage{mathrsfs}

\usepackage{amsmath}
\usepackage{latexsym,amsfonts,amssymb}

\newcommand{\bfv}{{\bf{v}}}
\newcommand{\bfn}{{\bf{n}}}

\newcommand{\bfe}{{\bf{e}}}

\begin{document}

\title{Linearized stability theorem for invariant and quasi-invariant parabolic differential equations in Banach manifolds with applications to
  free boundary problems}
\author{Shangbin Cui\\[0.2cm]
  {\small School of Mathematics, Sun Yat-Sen University, Guangzhou, Guangdong 510275,}\\
  {\small People's Republic of China. E-mail:\,cuishb@mail.sysu.edu.cn}}
\date{}
 \maketitle

\begin{abstract}
  If a differential equation in a Banach manifold is invariant or quasi-invariant under the action of one or more Lie groups, then its stationary
  points cannot be isolated, so that classical linearized stability theorem does not apply to it. The first main purpose of this paper is to establish
  a linearized stability theorem for parabolic differential equations in Banach manifolds which are either invariant or quasi-invariant under actions
  of a number of Lie groups. The second purpose of this paper is to apply this theorem to analyze stability of stationary solutions of some free
  boundary problems. In order to apply the abstract result to concrete free boundary problems, Banach manifold made up of certain kind of domains
  such as simple domains in ${\mathbf{R}}^n$ is a fundamental tool which seems to have not been well-studied in the literature yet. Hence in this paper
  we also make some basic investigation to a such manifold. In Section 5 we use Nash-Moser implicit function theorem to prove an interesting result
  for an obstacle problem which says that if the domain $\Omega$ of this obstacle problem is a small perturbation of a sphere then its interface
  $\Gamma$ is smooth and depends on $\Omega$ smoothly. By using these results, in the last section we prove asymptotic stability of radial stationary
  solution of a free boundary problem modeling the growth of necrotic tumors, which has been kept open for over ten years.
\medskip

\medskip

\end{abstract}

\section{Introduction}
\setcounter{equation}{0}

\hskip 2em
  Classical linearized stability theorem is an important fundamental theorem in the theory of ordinary differential equations. It states that for a
  differential equation $x'=F(x)$ in ${\mathbf{R}}^n$ with an isolated stationary point $x^*$, i.e., $F(x^*)=0$, where $F\in C^1(O,{\mathbf{R}}^n)$ for
  some open subset $O$ of ${\mathbf{R}}^n$ and $x^*\in O$, if $s:=\displaystyle\max_{1\leqslant j\leqslant n}{\rm Re}\lambda_j<0$, where $\lambda_j$'s
  are all eigenvalues of $F'(x^*)$, then $x^*$ is asymptotically stable, whereas if $s>0$ then $x^*$ is unstable.

  The above theorem has been successfully extended to quasilinear and fully nonlinear parabolic differential equations in Banach spaces by
  Poitier-Ferry \cite{Poi}, Lunardi \cite{Lun1}, Drangeid \cite{Dra}, and Da Prato and Lunardi \cite{Pra}, extending an earlier well-known result
  on semilnear parabolic differential equations in Banach spaces; see Chapter 9 of \cite{Lun2} for an exposition on this topic. The extended
  theorem states as follows: Let $X$ be a Banach space and $X_0$ an embedded Banach subspace of $X$ (see Section 3 for this concept). Let $O$ be an
  open subset of $X_0$ and $F\in C^{2-0}(O,X)$. Consider the differential equation $x'=F(x)$ in $X$. We say this equation is of parabolic type
  if for any $x\in O$, the operator $F'(x)\in L(X_0,X)$ is a sectorial operator when regarded as an unbounded linear operator in $X$ with domain
  $X_0$, and the graph norm of $X_0$ is equivalent to its own norm $\|\cdot\|_{X_0}$. Assume that the differential equation $x'=F(x)$ is of parabolic
  type and it has an isolated stationary point $x^*\in O$, i.e., $F(x^*)=0$. Let $s$ be the spectrum bound of $F'(x^*)$, i.e., $s=\sup\{{\rm Re}
  \lambda:\lambda\in\sigma(F'(x^*))\}$, where as usual $\sigma(\cdot)$ denotes the spectrum of a linear operator. Then we have the following assertion:
  If $s<0$ then $x^*$ is asymptotically stable, whereas if $s>0$ then $x^*$ is unstable. See \cite{Lun2} for the proof of this theorem\footnotemark[1].
\footnotetext[1]{Extension of linearized stability theorem to non-parabolic differential equations in Banach spaces is a very interesting topic worthy
  of great efforts. We refer the reader to see \cite{EngNag} for related topic.}

  In applications, however, we often encounter the situation $s=0$ and, more unpleasantly, the stationary point $x^*$ is not isolated but is contained
  in a manifold made up of stationary points. This situation is typical in the study of two kinds of evolutionary free boundary problems: The motion
  of liquid drops and the growth of tumors. In such problems, the related differential equations (i.e., the differential equations in Banach spaces or
  manifolds reduced from original free boundary problems) usually have some significant symmetric properties, i.e., they are invariant or
  quasi-invariant under some Lie group actions. It is clear that if a differential equation $x'=F(x)$ in a Banach space or manifold $X$ is invariant
  or quasi-invariant under a Lie group action $(G,p)$, where $G$ is a non-discrete Lie group and $p:G\times X\to X$ is an action of the group $G$ to
  $X$, then for stationary point $x^*$ of this equation, all points in $X$ of the form $p(a,x^*)$ with $a\in G$ are also its stationary points, so
  that generally speaking no stationary points are isolated and, as a result, if $\sigma(F'(x^*))$ does not intersect the open right-half plane
  then $s=0$, and the linearized stability theorem reviewed above does not apply, at least directly, to this kind of equations.

  Sometimes such deficiency of non-isolation of the stationary point can be remedied through imposition of additional conditions related to
  conservative quantities. A typical example in this line is the following free boundary problem modeling the motion of liquid drops (cf.
  \cite{Sol1,Sol2}):
\begin{equation}
\left\{
\begin{array}{cl}
   \partial_t\bfv-\nu\Delta\bfv+(\bfv\cdot\nabla)\bfv+\nabla\pi=0, &\quad\;\; x\in\Omega(t), \;\; t>0,\\
   \nabla\cdot\bfv=0,   &\quad \;\; x\in\Omega(t), \;\; t>0,\\
   T(\bfv,\pi)\bfn=-\gamma\kappa\bfn,   &\quad \;\; x\in\partial\Omega(t), \;\; t>0,\\
   V_n=\bfv\cdot\bfn, &\quad \;\; x\in\partial\Omega(t), \;\; t>0.
\end{array}
\right.
\end{equation}
  Here $\Omega(t)$ is the domain in ${\mathbf{R}}^3$ occupied by the liquid drop at time $t$, $\bfv$ and $\pi$ are velocity and pressure fields in the
  liquid drop, respectively, $\nu>0$ is the viscosity coefficient constant, $\gamma>0$ is the surface tension coefficient constant, $\kappa$ and
  $\bfn$ are the mean curvature field (whose sign is designated to be positive for spheres) and the unit outward normal field on $\partial\Omega(t)$,
  respectively, $V_n$ is the normal velocity of the free boundary $\partial\Omega(t)$, and $T(\bfv,\pi)=\nu S(\bfv)-\pi I$ ($I=$ the third-order
  identity matrix) and $S(\bfv)=\nabla\otimes\bfv+(\nabla\otimes\bfv)^T$
  are the stress tensor and the doubled strain tensor, respectively. It is clear that (1.1) is either invariant or quasi-invariant under the following
  three types of group actions in ${\mathbf{R}}^3$: scaling (quasi-invariant), translation (invariant), and rotation (invariant). Thus if
  $(\bfv_s,\pi_s,\Omega_s)$ is a stationary solution of the above problem and it is not a fixed point under some of these group actions or subgroup
  actions then all points in its trajectory of such group or subgroup actions are also stationary solutions, so that stationary solutions form a
  nontrivial manifold. For instance, for the trivial stationary solution
$$
   \bfv_s(x)=(0,0,0), \qquad  \pi_s(x)=\gamma, \qquad \Omega_s=B(0,1)=\{x\in{\mathbf{R}}^n:|x|<1\},
$$
  actions of the scaling and translation groups to it yield a four-dimensional manifold made up of stationary solutions (note that the action of
  rotation group does not yield new stationary solutions). It is known that the problem (1.1) has also the following non-trivial stationary solution
  (cf. \cite{Sol1}):
$$
   \bfv_s(x)=\mu(x_2,-x_1,0), \qquad  \pi_s(x)=\frac{1}{2}\mu^2(x_1^2+x_2^2)+1, \qquad  \Omega_s=\{x\in{\mathbf{R}}^n:r<\rho_s(\omega)\}.
$$
  Here and throughout this paper $(r,\omega)$ denotes the polar coordinate of a point $x\in{\mathbf{R}}^3\backslash\{0\}$ (in Section 4
  $x\in{\mathbf{R}}^n\backslash\{0\}$ with $n\geqslant 2$), i.e., $r=|x|$ and $\omega=x/r$. Besides, $\mu$ is a nonzero constant, and $\rho_s$ is
  the unique solution of the following
  equation:
$$
  2\gamma\kappa(\rho_s)-\frac{\mu^2}{2}\rho_s^2(\omega)(\omega_1^2+\omega_2^2)=1, \quad  \omega\in{\mathbf{S}}^2,
$$
  where $\kappa(\rho)$ is the mean curvature of the surface $r=\rho(\omega)$:
$$
  \kappa(\rho)=-\frac{1}{2}\frac{1}{\rho}\Big[\frac{\Delta_{\omega}\rho}{\sqrt{g}}+\nabla_{\omega}\rho\cdot
  \nabla_{\omega}\Big(\frac{1}{\sqrt{g}}\Big)-\frac{2\rho}{\sqrt{g}}\Big], \quad  g=\rho^2+|\nabla_{\omega}\rho|^2.
$$
  Here and throughout this paper $\Delta_{\omega}$ and $\nabla_{\omega}$ denote the Laplace-Beltrami operator and the gradient field on the sphere
  ${\mathbf{S}}^2$ (in Section 4 on ${\mathbf{S}}^{n-1}$ with $n\geqslant 2$), respectively.
  For this stationary solution, not only actions of scaling and translation groups to it yield new stationary solutions, but also the action of the
  two-dimensional  subgroup $O_{12}(3)$ of the rotation group $O(3)$ generated by rotations with rotation axes $Ox_1$ and $Ox_2$ also yields new
  stationary solutions (note that the action of a rotation with rotation axes $Ox_3$ does not produce new stationary solutions), so that the above
  stationary solution is contained in a six-dimensional manifold made up of stationary solutions. But such lack of non-isolation of the stationary
  point can be remedied by imposing the following additional conditions:
\begin{equation}
\left\{
\begin{array}{cl}
   |\Omega(t)|=c, &\quad\;\; t>0,\\
   \displaystyle\int_{\Omega(t)}\bfv dx=0, \quad \int_{\Omega(t)}(\bfv\times x)dx=\beta\bfe_3, &\quad\;\; t>0,\\
   \displaystyle\int_{\Omega(t)}xdx=0,  &\quad\;\; t>0,
\end{array}
\right.
\end{equation}
  where $c>0$ and $\beta\in{\mathbf{R}}$ are given constants. The condition $(1.2)_1$ is due to conservation of volume of $\Omega(t)$:
  $\displaystyle\frac{d}{dt}\Big(\int_{\Omega(t)}dx\Big)=0$, which eliminates the possibility to get new stationary solutions through scaling (by
  fixing $c$); the conditions in $(1.2)_2$ are due to conservation of momentum and conservation of angular momentum, which eliminates the possibility
  to get new stationary solutions through rotation (by fixing $\beta$); the condition in $(1.2)_3$ is obtained from the momentum-free and
  incompressibility conditions, which eliminates the possibility to get new stationary solutions through translation. Hence, by imposing
  these additional conditions, the stationary solution $(\bfv_s,p_s,\Omega_s)$ becomes isolated and $0$ is removed from spectrum of the linearized
  operator, so that the linearized stability theorem reviewed above becomes applicable (note that the problem (1.1) can be reduced into a parabolic
  differential equation in some Banach manifold; cf. \cite{WuCui} for discussion to a free boundary problem modeling tumor growth which has some
  similar features as the problem (1.1)).

  We note that some other free boundary problems modeling fluid drops by using Stokes equations or stationary Stokes equations rather than
  Navier-Stokes equations have similar features as the above model. For work on such models we refer the reader to see e.g. \cite{EscPro},
  \cite{GunPro} and references therein.

  However, there are some other free boundary problems for which similar conservative quantities do not exist, so that the above technique is not
  always effective. A typical example is the following free boundary problem modeling the growth of a tumor cultivated in laboratory:
\begin{equation}
\left\{
\begin{array}{rll}
   \Delta\sigma=&f(\sigma) &\quad\;\; \mbox{in} \;\; \Omega(t),\;\; t>0,\\
   -\Delta\pi=&g(\sigma)   &\quad\;\; \mbox{in} \;\; \Omega(t),\;\; t>0,\\
   \sigma=&\bar{\sigma}   &\quad\;\; \mbox{on} \;\; \partial\Omega(t),\;\; t>0,\\
   \pi=&\gamma\kappa  &\quad\;\; \mbox{on} \;\; \partial\Omega(t),\;\; t>0,\\
   V_n=&-\partial_{n}\pi  &\quad\;\; \mbox{on} \;\; \partial\Omega(t),\;\; t>0.
\end{array}
\right.
\end{equation}
  Here $\Omega(t)$ is the domain in ${\mathbf{R}}^3$ occupied by the tumor at time $t$, $\sigma$ is the nutrient concentration in the tumor, $\pi$ is
  the pressure in the tumor, $\bar{\sigma}$ is a positive constant reflecting constant supply of nutrient from tumor surface, $\gamma$, $\kappa$ and
  $V_n$ are as before, $\partial_{n}$ denotes the derivative in normal direction of $\partial\Omega(t)$, and $f(\sigma)$, $g(\sigma)$ are given
  functions representing nutrient consumption rate and tumor cell proliferation rate, respectively, that nutrient concentration in level $\sigma$
  can sustain. Under suitable conditions on $f$, $g$ and $\bar{\sigma}$, the above equations have a unique radial stationary solution (cf. \cite{Cui1},
  \cite{Cui2}, \cite{Cui5}). It is clear that the above equations are also invariant under actions of the translation and rotation groups in
  ${\mathbf{R}}^3$ (but is not invariant or quasi-invariant under the action of scaling group). Hence, by acting to the radial stationary solution with
  the translation group, we get a three-dimensional manifold made up of stationary solutions (note that the action of rotation group does not yield
  new stationary solutions). For this free boundary problem, there does not exist a conservative quantity for us to isolate the radial stationary
  solution from other stationary solutions obtained from translating the radial stationary solution.

  The purpose of this paper is twofold. The first goal is to prove an abstract linearized stability theorem for parabolic differential equations
  in Banach manifolds which are invariant or quasi-invariant under some Lie group actions. This abstract theorem can be roughly stated as follows:
  Let $\mathfrak{M}$ be a Banach manifold and $G_i$, $i=1,2,\cdots,N$, be a number of Lie groups acting on $\mathfrak{M}$ through group actions
  $p_i:G_i\times\mathfrak{M}\to\mathfrak{M}$, $i=1,2,\cdots,N$, respectively. We assume that these Lie group actions are quasi-commutative, i.e., for
  any $1\leqslant i,j\leqslant N$ with $i\neq j$ there exists corresponding mapping $f_{ij}:G_i\times G_j\to G_j$ such that the following relation
  holds:
$$
  p_i(a,p_j(b,x))=p_j(f_{ij}(a,b),p_i(a,x)), \quad \forall x\in\mathfrak{M}, \;\;\forall a\in G_i,\;\; \forall b\in G_j.
$$
  We denote $g(a,x)=p_1(a_1,p_2(a_2,\cdots, p_N(a_N,x)))$ for $x\in\mathfrak{M}$ and $a=(a_1,a_2,\cdots,a_N)\in G:=G_1\times G_2\times\cdots \times
  G_N$, and assume that ${\rm rank}\,\partial_ag(a,x)=n:=\dim G_1+\dim G_2+\cdots+\dim G_N$. Consider an autonomous differential equation in
  $\mathfrak{M}$:
\begin{equation}
  x'=\mathscr{F}(x),
\end{equation}
  where $\mathscr{F}$ is a vector field on $\mathfrak{M}$ with domain $\mathfrak{M}_0$, which is an embedded Banach submanifold of $\mathfrak{M}$
  (see Section 3 for definitions of these concepts). We assume this equation is of parabolic type and is quasi-invariant under all Lie group actions
  $(G_i,p_i)$, $i=1,2,\cdots,N$, i.e., for each $1\leqslant i\leqslant N$ there exists corresponding group homomorphism $\theta_i:G_i\to {\mathbf{R}}_+
  =(0,\infty)$ such that the following relation holds:
$$
  \mathscr{F}(p_i(a,x))=\theta_i(a)\partial_xp_i(a,x)\mathscr{F}(x), \quad \forall x\in\mathfrak{M}, \;\;\forall a\in G_i.
$$
  If $\theta_i(a)\equiv 1$ for some $i$ then the equation (1.4) is said to be invariant under the Lie group action $(G_i,p_i)$.
  Let $x^*\in\mathfrak{M}_0$ be a stationary point of the equation (1.4), i.e., $\mathscr{F}(x^*)=0$, and set $M_c=\{g(a,x^*):a\in G\}$,
  i.e., $M_c$ is the combined trajectory of $x^*$ under these Lie group actions. Assume, in addition to some other conditions (see Section 3 for
  details), that $\dim{\rm Ker}\mathscr{F}'(x^*)=n$.
  Let $s^*=\sup\{{\rm Re}\lambda:\lambda\in\sigma(\mathscr{F}'(x^*))\backslash\{0\}\}$. Then we have the following assertions: If $s^*<0$ then there
  exists a Banach submanifold $M_s$ of $\mathfrak{M}_0$ of codimension $n$ such that the solution $x=x(t)$ of (1.4) satisfies $\displaystyle
  \lim_{t\to\infty}x(t)=x^*$ if and only if $x(0)\in M_s$ (in this case $x(t)\in M_s$ for all $t\geqslant 0$), and for any $x_0$ in a small
  neighborhood of $x^*$, there exist unique $a\in G$ and $\hat{x}_0\in M_s$ such that $x_0=g(a,\hat{x}_0)$, and the solution $x=x(t)$ of (1.4) with
  initial data $x(0)=x_0$ satisfies $\displaystyle\lim_{t\to\infty}x(t)=g(a,x^*)$. See Theorem 3.4 in Section 3 for precise statement of this result.
  We say $x^*$ is asymptotically stable module the groups $G_1$, $G_2$, $\cdots$, $G_N$. If $s^*>0$ then $x^*$ is clearly unstable.

  In the case $N=1$, a local version of the above result was already established in the reference \cite{Cui3}, where the theorem was proved for
  differential equations in Banach spaces rather than Banach manifolds and, correspondingly, the concept of local Lie group action rather than Lie
  group action was employed. Application of the concept of local Lie group action makes the theorem proved in \cite{Cui3} look somewhat awkward. Surely,
  if the concept of Lie group action rather than local Lie group action were used in \cite{Cui3}, the theorem would look more pleasant but it would
  no longer be applicable to free boundary problems. To solve this problem in the present paper we re-establish the theorem in Banach manifolds, and
  furthermore, we consider actions of more than one Lie groups so that the abstract theorem can be applied to more general free boundary problems.

  The second purpose of this paper is to apply the abstract result talked above to study asymptotic stability of the radial stationary solution of
  the problem (1.3) in the case that $f$, $g$ are the following discontinuous functions:
\begin{equation}
   f(\sigma)=\lambda\sigma H(\sigma-\hat{\sigma}),\qquad
   g(\sigma)=a(\sigma-\tilde{\sigma})H(\sigma-\hat{\sigma})-b,
\end{equation}
  where $H$ is the Heaviside function: $H(s)=1$ for $s>0$ and $H(s)=0$ for $s\leqslant 0$, and $\lambda$, $a$, $b$, $\hat{\sigma}$ and $\tilde{\sigma}$
  are positive constants: $\lambda$ is the consumption rate coefficient of nutrient by tumor cells, $a$ is the proliferation rate coefficient of tumor
  cells (=the birth rate of tumor cells that a unit amount of nutrient can sustain), $b$ is the dissolution rate of dead cells, $\hat{\sigma}$ is a
  threshold value of nutrient concentration to sustain tumor cells alive and proliferating, i.e., only in the region where $\sigma>\hat{\sigma}$ tumor
  cells are alive and proliferating, and $\tilde{\sigma}=\hat{\sigma}-(b/a)$. We assume that $0<\hat{\sigma}<\bar{\sigma}$, $b<a\hat{\sigma}$ (so that
  $0<\tilde{\sigma}<\hat{\sigma}$) and, for simplicity of notations $\lambda=\bar{\sigma}=1$, which can always be achieved through rescaling. For
  $f$, $g$ given by (1.5), the problem (1.3) models growth of necrotic tumors, cf. \cite{ByrC2} for details.

  If instead of (1.5) the functions $f$, $g$ are smooth monotone increasing functions in $[0,\infty)$ and satisfy the properties
  $f(0)=0$, $g(0)<0$ and $g(\infty)>0$, the problem (1.3) models the growth of nonnecrotic tumors, cf. \cite{ByrC1}, which has been intensively
  studied, cf. \cite{BazFri1}, \cite{BazFri2}, \cite{Cui1}, \cite{Cui4}, \cite{FH1}, \cite{FH2} and references therein. It was proved that there
  exists a threshold value $\gamma^*>0$ for the surface tension coefficient $\gamma$, such that if $\gamma>\gamma^*$ then the unique radial stationary
  solution is asymptotically stable module translations, whereas if $\gamma<\gamma^*$ then it is unstable. In the necrotic case, analysis of the
  problem (1.3) is much harder: In addition to the outer free boundary $\partial\Omega(t)$ whose evolution is governed by the equation $V_n=
  -\partial_{\nu}\pi$, discontinuity of the functions $f$, $g$ at $\sigma=\hat{\sigma}$ produces an inner free boundary or interface
  $\Gamma(t)$ dividing the domain $\Omega(t)$ into two disjoint regions, with the outer region $\Omega_{\rm liv}(t)=\{x\in\Omega(t):\sigma(x,t)>
  \hat{\sigma}\}$ being the living shell and the inner region $\Omega_{\rm nec}(t)={\rm int}\{x\in\Omega(t):\sigma(x,t)=\hat{\sigma}\}$ the necrotic
  core of the tumor. Main difficulty of analysis is caused by existence of the inner free boundary $\Gamma(t)$, which is of obstacle type so that we
  have not an obvious equation governing its evolution to use. In \cite{Cui2}, \cite{Cui5} the spherically symmetric version of this problem was
  studied. It was proved that this problem has a unique radial stationary solution which is asymptotically stable under spherically symmetric
  perturbations. However, whether this stationary solution is asymptotically stable under spherically non-symmetric perturbations has been kept
  unknown for over ten years. We note that recently some numerical results on this problem were obtained by Hao et al \cite{Hao1}. In this paper
  we use the abstract result mentioned above to prove that similar result as for the nonnecrotic case is also valid for the necrotic case, i.e.,
  there exists a threshold value $\gamma^*$ for the surface tension coefficient $\gamma$ such that if $\gamma>\gamma^*$ then the radial stationary
  solution is asymptotically stable module translations under spherically non-symmetric perturbations, whereas if $\gamma<\gamma^*$ then it is
  unstable under spherically non-symmetric perturbations. See Theorem 6.8 in Section 6 for precise statement of this result.

  In getting the above result, a key step is to prove that the inner free boundary $\Gamma(t)$ is smooth and depends on the outer free boundary
  $\partial\Omega(t)$ smoothly. To be more precise, let us consider the following obstacle problem:
\begin{equation}
\left\{
\begin{array}{cl}
   -\Delta\sigma+\sigma\geqslant 0, \quad \sigma\geqslant\hat{\sigma}, \quad (-\Delta\sigma+\sigma)(\sigma-\hat{\sigma})=0 &\quad\;\,
   \mbox{in}\;\;\Omega,\\
   \sigma=1   &\quad \mbox{on} \;\; \partial\Omega,
\end{array}
\right.
\end{equation}
  where $\Omega$ is a given bounded domain in ${\mathbf{R}}^3$ with a $C^2$ boundary and $0<\hat{\sigma}<1$ is a given constant. It is easy to prove
  that this problem has a unique solution $\sigma\in\displaystyle\cap_{1\leqslant p<\infty}W^{2,p}(\Omega)$ which satisfies $\hat{\sigma}\leqslant
  \sigma\leqslant 1$. Let $\Omega_{\rm liv}=\{x\in\Omega:\sigma(x)>\hat{\sigma}\}$, $\Omega_{\rm nec}={\rm int}\{x\in\Omega:\sigma(x)=\hat{\sigma}\}$
  and $\Gamma=\partial\Omega_{\rm liv}\cap\partial\Omega_{\rm nec}$. $\Gamma$ is called the free boundary or interface of the above obstacle problem.
  It is well-known that regularity of free boundaries of obstacle problems is a very hard topic. We refer the reader to see \cite{Caff} and
  \cite{Fried} and references therein for some classical results on it. For our purpose, in this paper we shall use Nash-Moser implicit function
  theorem to prove the following interesting result: For the problem (1.6), the free boundary $\Gamma$ is smooth and depends on $\Omega$ smoothly if
  $\Omega$ is a small perturbation of a sphere; see Theorem 5.2 in Section 5.

  In order to use the abstract theorem talked before to study the problem (1.3), Banach manifolds of certain domains in ${\mathbf{R}}^n$ are an
  essential tool which, however, seems to have not yet been well-studied in the literature. Note that the concept of Frech\'{e}t manifold of smooth
  domains in ${\mathbf{R}}^n$ was already introduced thirty years ago by Hamilton in his famous work \cite{Ham1}, but Frech\'{e}t manifold is
  not a good tool for the study of asymptotic behavior of solutions of free boundary problems: A useful linearized stability theorem for differential
  equations in Frech\'{e}t manifolds is very hard to establish. Due to this reason, in this paper we make some basic investigation
  to the Banach manifold of simple domains in ${\mathbf{R}}^n$. A basic difficulty on this topic is that such a manifold is only a topological Banach
  manifold but not a differentiable Banach manifold, i.e., it does not possess a differentiable structure. It follows that the concept of differential
  equations in a such manifold cannot be understood in usual sense. Our observation is that though the manifold of simple domains itself does
  not possess a differentiable structure, certain of its embedded submanifolds possess differentiable structures in the topology inherited from the
  whole manifold, which is sufficient for the concept of differential equations in a such manifold. See Section 2 for details. Here we mention that
  Banach manifold of closed hypersurfaces in ${\mathbf{R}}^n$ was recently considered by Pr\"{u}ss and Simonett \cite{PruSim}. But the idea and
  results on this topic obtained in the present work are completely new.

  The organization of the rest part is as follows. In Section 2 we make some basic investigation to the Banach manifold made by simple domains
  in ${\mathbf{R}}^n$. The purpose of this section is to provide a basic tool for the study of evolutionary free boundary problems: In Sections 4
  and 6 we shall reduce evolutionary free boundary problems into differential equations in this manifold. A significant property of this manifold
  is that it does not possess a differentiable structure, but some of its embedded submanifolds possess differentiable structures in its topology.
  This motivates the study in Section 3 to differential equations in Banach manifolds without differentiable structure, where the linearized stability
  theorem for parabolic differential equations in non-differentiable Banach manifolds which are invariant or quasi-invariant under a number of Lie
  group actions is established. In Section 4 we show some simple applications of this abstract result. Section 5 is a preparational section for
  Section 6. We prove that the free boundary $\Gamma$ of the obstacle problem (1.6) is smooth and depends on $\Omega$ smoothly, provided
  $\Omega$ is a small perturbation of a sphere. We use this result to prove the surface tension free part $\pi_0$ of the solution of the second
  equation in (1.3) has the property that the map $\Omega\mapsto \partial_{\bfn}\pi_0|_{\partial\Omega}$ is smooth, even if $f$, $g$ are
  discontinuous functions given by (1.5). In the final section we use the abstract result obtained in Section 3 to study asymptotic stability of the
  radial stationary solution of the tumor model (1.3) in the necrotic case.

\section{Banach manifolds of simple domains in ${\mathbf{R}}^n$}
\setcounter{equation}{0}

\hskip 2em
  In this section we make a basic study to Banach manifolds of simple domains in ${\mathbf{R}}^n$. Such a manifold is only a topological manifold and
  does not possess a differentiable structure. However, some of its embedded Banach submanifolds possess differentiable structures in its topology,
  so that the concept of differential equations in a such manifold still makes sense. The purpose of this section is to provide a basic tool for study
  of evolutionary free boundary problems, and the results obtained in this section will be applied in Sections 4 and 6. Also, from the discussion
  of this section it will be clear why in the next section differential equations will be considered in non-differentiable Banach manifolds, not in
  differentiable Banach manifolds.

\subsection{Basic concepts}

\hskip 2em
  As usual, given a nonnegative integer $m$, a real number $\mu\in [0,1]$, a bounded open set $\Omega\subseteq{\mathbf{R}}^n$ and a sufficiently
  smooth bounded closed hypersurface $S\subseteq{\mathbf{R}}^n$, the notations $C^{m+\mu}(\overline{\Omega})$ and $C^{m+\mu}(S)$ denote the usual
  $m\!+\!\mu$-th order H\"{o}lder spaces on $\overline{\Omega}$ and $S$, respectively, and the notation $C^{m+\mu}(\overline{\Omega},{\mathbf{R}}^n)$
  denotes the usual $m\!+\!\mu$-th order $n$-vector H\"{o}lder space on $\overline{\Omega}$, with the cases $\mu=0,1$ understood in conventional sense.
  We use the notation $\dot{C}^{m+\mu}(\overline{\Omega})$ to denote the closure of $C^{\infty}(\overline{\Omega})$ in $C^{m+\mu}(\overline{\Omega})$,
  and similarly for the notations $\dot{C}^{m+\mu}(S)$ and $\dot{C}^{m+\mu}(\overline{\Omega},{\mathbf{R}}^n)$. The last three spaces are called
  $m\!+\!\mu$-th order little H\"{o}lder spaces. A significant difference between little H\"{o}lder spaces and H\"{o}lder spaces is that for
  nonnegative integers $k,m$ and real numbers $\mu,\nu\in [0,1]$, if $k+\nu>m+\mu$ then $\dot{C}^{k+\nu}(\overline{\Omega})$ (resp.
  $\dot{C}^{k+\nu}(S)$, $\dot{C}^{k+\nu}(\overline{\Omega},{\mathbf{R}}^n)$) is dense in $\dot{C}^{m+\mu}(\overline{\Omega})$ (resp.
  $\dot{C}^{m+\mu}(S)$, $\dot{C}^{m+\mu}(\overline{\Omega},{\mathbf{R}}^n)$), but $C^{k+\nu}(\overline{\Omega})$ (resp. $C^{k+\nu}(S)$,
  $C^{k+\nu}(\overline{\Omega},{\mathbf{R}}^n)$) is not dense in $C^{m+\mu}(\overline{\Omega})$ (resp. $C^{m+\mu}(S)$,
  $C^{m+\mu}(\overline{\Omega},{\mathbf{R}}^n)$).
\medskip

  {\bf Definition 2.1}\ \ {\em Let $m$ be a positive integer and $0\leqslant\mu\leqslant 1$. An open set $\Omega\subseteq{\mathbf{R}}^n$ is said to be
  a simple $C^{m+\mu}$-domain if $\Omega$ is $C^{m+\mu}$-diffeomorphic to the open unit sphere $B(0,1)$ in ${\mathbf{R}}^n$, i.e., there exists a bijective
  mapping $\Phi:\overline{B(0,1)}\to\overline{\Omega}$ satisfying the following properties:
$$
  \Phi\in C^{m+\mu}(\overline{B(0,1)},{\mathbf{R}}^n) \quad \mbox{and} \quad \Phi^{-1}\in C^{m+\mu}(\overline{\Omega},{\mathbf{R}}^n).
$$
  We use the notation $\mathfrak{D}^{m+\mu}({\mathbf{R}}^n)$ to denote the set of all simple $C^{m+\mu}$-domains in ${\mathbf{R}}^n$. If instead of
  $C^{m+\mu}$ the notation $\dot{C}^{m+\mu}$ is used in the above relations, then the notation $\dot{\mathfrak{D}}^{m+\mu}({\mathbf{R}}^n)$ is used
  correspondingly.}
\medskip

  From \cite{Ham1} we know that all smooth simple domains in ${\mathbf{R}}^n$ form a Frech\'{e}t manifold $\mathfrak{D}^{\infty}({\mathbf{R}}^n)$ built
  on the Frech\'{e}t space $C^{\infty}({\mathbf{S}}^{n-1})$, with tangent space at the point $\Omega\in\mathfrak{D}^{\infty}({\mathbf{R}}^n)$ being
  $T_{\Omega}(\mathfrak{D}^{\infty}({\mathbf{R}}^n))=C^{\infty}(\partial\Omega)$. In applications, however, just as Frech\'{e}t spaces are not as
  convenient to use as Banach spaces, Frech\'{e}t manifolds are not as convenient to use as Banach manifolds. Hence, in what follows we introduce a
  local chart for $\mathfrak{D}^{m+\mu}({\mathbf{R}}^n)$ ($m\in \mathbb{N}$, $0\leqslant\mu\leqslant 1$) at each of its point to make it into a
  Banach manifold.

  In the following of this section we always assume that $m,n$ are positive integers $\geqslant 2$, and $0\leqslant\mu\leqslant 1$.
  Let $\Omega\in\mathfrak{D}^{m+\mu}({\mathbf{R}}^n)$. Choose a closed $C^{\infty}$-hypersurface $S\subseteq{\mathbf{R}}^n$ sufficiently closed to
  $\partial\Omega$ such that the following three conditions are satisfied: $(a)$\ $S$ encloses a simple $C^{\infty}$-domain $Q$; $(b)$\ there exists
  $\delta>0$ such that by letting $\mathcal{R}=\{x\in{\mathbf{R}}^n:d(x,S)<\delta\}$, the mapping $\Psi:S\times(-\delta,\delta)\to \mathcal{R}$,
$$
   \Psi(x,t)=x+t\bfn(x), \quad \forall x\in S,\;\; \forall t\in (-\delta,\delta)
$$
  is a $C^{\infty}$-diffeomorphism of $S\times(-\delta,\delta)$ onto $\mathcal{R}$, where $\bfn$ denotes the unit normal field of $S$, outward pointing
  with respect to $Q$; $(c)$\ $\partial\Omega\subseteq\mathcal{R}$. Existence of a such hypersurface $S$ is well-known. Let $\Pi$ and $\Lambda$ be
  compositions of $\Psi^{-1}:\mathcal{R}\to S\times(-\delta,\delta)$ with the projections of $S\times (-\delta,\delta)$ onto $S$ and $(-\delta,\delta)$,
  respectively, i.e., $\Pi:\mathcal{R}\to S$, $\Lambda:\mathcal{R}\to(-\delta,\delta)$,
$$
  \Pi(\Psi(x,t))=x, \quad \Lambda(\Psi(x,t))=t, \quad \forall x\in S, \;\; \forall t\in(-\delta,\delta),
$$
  and
$$
  y=\Pi(y)+\Lambda(y)\bfn(\Pi(y)),\;\; \forall y\in \mathcal{R}.
$$
  It is easy to see that if $\delta$ is small enough then for any $y\in \mathcal{R}$, $\Pi(y)$ is the point in $S$ nearest to $y$ and $\Lambda(y)$ is
  the algebraic distance of $y$ to $S$, i.e, $\Lambda(y)=-d(y,S)$ if $y\in\mathcal{R}\cap Q$ and $\Lambda(y)=d(y,S)$ if $y\in\mathcal{R}\backslash
  Q$. Let $\mathcal{O}$ be the $\delta$-neighborhood of the null function in $C^{m+\mu}(S)$, i.e., $\mathcal{O}=\{\rho\in C^{m+\mu}(S):\,
  \|\rho\|_{C^{m+\mu}(S)}<\delta\}$. Since $\rho\in\mathcal{O}$ implies $\displaystyle\max_{x\in S}|\rho(x)|<\delta$, it follows that for any $\rho
  \in\mathcal{O}$ the mapping $\theta_{\rho}:S\to{\mathbf{R}}^n$,
$$
    \theta_{\rho}(x)=x+\rho(x)\bfn(x), \quad \forall x\in S
$$
  is a $C^{m+\mu}$-diffeomorphism of $S$ onto its image $S_{\rho}=\theta_{\rho}(S)$, which is a closed $C^{m+\mu}$-hypersurface contained in
  $\mathcal{R}$ enclosing a simple $C^{m+\mu}$-domain $\Omega_{\rho}$. Since $\partial\Omega\subseteq \mathcal{R}$, there exists a unique
  $\rho_0\in\mathcal{O}$ such that $\partial\Omega=S_{\rho_0}$ or $\Omega=\Omega_{\rho_0}$. Let
\begin{equation}
  \mathcal{U}=\{\Omega_{\rho}:\rho\in\mathcal{O}\},
\end{equation}
  and
\begin{equation}
  \varphi:\mathcal{U}\to C^{m+\mu}(S), \quad \varphi(\Omega_{\rho})=\rho, \;\;\forall\rho\in\mathcal{O}.
\end{equation}
  We call $(\mathcal{U},\varphi)$ a {\em regular local chart} of $\mathfrak{D}^{m+\mu}({\mathbf{R}}^n)$ at the point $\Omega$, and call the closed
  hypersurface $S$ the {\em base hypersurface} of this local chart (the phrase ``regular'' refers to the fact that the base hypersurface $S$ is smooth).
  Clearly, $C^{m+\mu}(S)\approx C^{m+\mu}({\mathbf{S}}^{n-1})$, i.e., $C^{m+\mu}(S)$ and $C^{m+\mu}({\mathbf{S}}^{n-1})$ are isomorphic to each other as
  Banach spaces.
\medskip

  {\bf Theorem 2.2}\ \ {\em With local chart $(\mathcal{U},\varphi)$ defined as above, $\mathfrak{D}^{m+\mu}({\mathbf{R}}^n)$ is a $($topological or
  $C^0)$ Banach manifold built on the Banach space $C^{m+\mu}({\mathbf{S}}^{n-1})$.}
\medskip

  {\em Proof}.\ \ Let $(\mathcal{U}_1,\varphi_1)$, $(\mathcal{U}_2,\varphi_2)$ be two local charts such that $\mathcal{U}_1\cap\mathcal{U}_2
  \neq\varnothing$, with base hypersurfaces $S_1$, $S_2$, respectively. Let $\bfn_i$ be the unit outward pointing normal field of $S_i$,
  $\mathcal{R}_i$ the neighborhood of $S_i$ and $\Pi_i:\mathcal{R}_i\to S_i$ the projection as defined above, $i=1,2$. Since $C^{m+\mu}(S_i)\approx
  C^{m+\mu}({\mathbf{S}}^{n-1})$, $i=1,2$, we only need to prove $\varphi_2\circ\varphi_1^{-1}\in C(\mathcal{O}_1,\mathcal{O}_2)$, where $\mathcal{O}_i
  =\varphi_i(\mathcal{U}_1\cap\mathcal{U}_2)\subseteq C^{m+\mu}(S_i)$,  $i=1,2$. Let $\Omega\in\mathcal{U}_1\cap\mathcal{U}_2$. Then there exist
  $\rho_i\in\mathcal{O}_i$, $i=1,2$, such that
$$
  \partial\Omega=\{x+\rho_1(x)\bfn_1(x):x\in S_1\}=\{y+\rho_2(y)\bfn_2(y):y\in S_2\}.
$$
  By definition, $(\varphi_2\circ\varphi_1^{-1})(\rho_1)=\rho_2$, $\forall\rho_1\in\mathcal{O}_1$. The above equality implies
$$
  x+\rho_1(x)\bfn_1(x)=y+\rho_2(y)\bfn_2(y), \quad x\in S_1, \;\; y\in S_2,
$$
  i.e, for any $x\in S_1$ there exists a unique corresponding $y\in S_2$ such that the above equality holds, and vice versa. Since for given $\rho_1
  \in\mathcal{O}_1$ and $y\in S_2$, the point $x\in S_1$ is uniquely determined, we write $x=\chi(\rho_1,y)$. Substituting this expression into the
  above equality and computing inner products of both sides of it with $\bfn_2(y)$, we get
\begin{equation}
  \rho_2(y)=\langle\chi(\rho_1,y)-y,\bfn_2(y)\rangle+\rho_1(\chi(\rho_1,y))\langle\bfn_1(\chi(\rho_1,y)),\bfn_2(y)\rangle, \quad \forall y\in S_2.
\end{equation}
  Clearly, the function $x=\chi(\rho_1,y)$ is the implicit function defined by the equation
\begin{equation}
  y=\Pi_2(x+\rho_1(x)\bfn_1(x)).
\end{equation}
  This implicit function equation is regular, i.e., the derivative $\partial y/\partial x: T_x(S_1)\to T_y(S_2)$ has bounded inverse because
  $(\partial y/\partial x)^{-1}=\partial x/\partial y$ is the derivative of the function $y\mapsto x=\Pi_1(y+\rho_2(y)\bfn_2(y))$ which belongs to
  $C^{m+\mu}(S_2,S_1)$. Since the mapping $(\rho_1,x)\mapsto\Pi_2(x+\rho_1(x)\bfn_1(x))$ belongs to $C^{m+\mu}(\mathcal{O}_1\times S_1,S_2)$, by
  the implicit function theorem it follows that $\chi\in C^{m+\mu}(\mathcal{O}_1\times S_2,S_1)$. This implies that
\begin{equation}
  [\rho_1\mapsto\chi(\rho_1,\cdot)]\in C^k(\mathcal{O}_1,C^{m-k+\mu}(S_2,S_1)), \quad k=0,1,\cdots,m.
\end{equation}
  From (2.3) and (2.5) we easily obtain
\begin{equation}
  [\rho_1\mapsto\rho_2]\in C^k(\mathcal{O}_1,C^{m-k+\mu}(S_2)), \quad k=0,1,\cdots,m.
\end{equation}
  Hence, in particular, $\varphi_2\circ\varphi_1^{-1}\in C(\mathcal{O}_1,\mathcal{O}_2)$, as desired. This proves the theorem. $\quad\Box$
\medskip

  {\em Remark}.\ \ The proof of the above theorem shows that for a regular local chart $(\mathcal{U},\varphi)$ of $\mathfrak{D}^{m+\mu}({\mathbf{R}}^n)$
  and any $k\in{\mathbf{Z}}_+$, if $\Omega\in\mathcal{U}\cap\mathfrak{D}^{m+k+\mu}({\mathbf{R}}^n)$ then $\varphi(\Omega)\in C^{m+k+\mu}(S)$. Beside,
  from the above proof we see that $\mathfrak{D}^{m+\mu}({\mathbf{R}}^n)$ does not possess a differentiable structure. However, it is still possible to
  define differentiable points in $\mathfrak{D}^{m+\mu}({\mathbf{R}}^n)$ and tangent spaces at such points. These concepts will enable us to define
  differentiable curves in $\mathfrak{D}^{m+\mu}({\mathbf{R}}^n)$ and their tangent fields.
\medskip

  {\bf Lemma 2.3}\ \ {\em Let $(\mathcal{U}_1,\varphi_1)$, $(\mathcal{U}_2,\varphi_2)$ be two local charts of $\mathfrak{D}^{m+\mu}({\mathbf{R}}^n)$
  such that $\mathcal{U}_1\cap\mathcal{U}_2\neq\varnothing$.  Let $\mathcal{O}_i=\varphi_i(\mathcal{U}_1\cap\mathcal{U}_2)\subseteq C^{m+\mu}(S_i)$,
  $i=1,2$. Then for any $\Omega\in\mathcal{U}_1\cap\mathcal{U}_2\cap\mathfrak{D}^{m+\mu+1}({\mathbf{R}}^n)$, $\varphi_2\circ\varphi_1^{-1}$ is
  differentiable at $\rho=\varphi_1(\Omega)$: $(\varphi_2\circ\varphi_1^{-1})'(\rho)\in L(C^{m+\mu}(S_1),C^{m+\mu}(S_2))$, and $[\rho\mapsto
  (\varphi_2\circ\varphi_1^{-1})'(\rho)]\in C(\mathcal{O}_1\cap C^{m+\mu+1}(S_1),L(C^{m+\mu}(S_1),C^{m+\mu}(S_2)))$.}
\medskip

  {\em Proof}.\ \ Let $\sigma:\mathcal{O}_1\times S_2\to{\mathbf{R}}$ be the mapping given by the right-hand side of (2.3), i.e.,
$$
  \sigma(\rho,y)=\langle\chi(\rho,y)-y,\bfn_2(y)\rangle+\rho(\chi(\rho,y))\langle\bfn_1(\chi(\rho,y)),\bfn_2(y)\rangle, \quad
   \rho\in\mathcal{O}_1, \;\; y\in S_2,
$$
  where $\chi:\mathcal{O}_1\times S_2\to S_1$ is as before. Then $(\varphi_2\circ\varphi_1^{-1})(\rho)=\sigma(\rho,\cdot)$.
  From the above expression and the fact $\chi\in C^{m+\mu}(\mathcal{O}_1\times S_2,S_1)$ we see $\sigma\in C^{m+\mu}(\mathcal{O}_1\times S_2)$.
  Moreover, for any $\zeta\in C^{m+\mu}(S_1)$,
\begin{equation}
\begin{array}{rl}
   \partial_{\rho}\sigma(\rho,y)\zeta=&\langle\partial_{\rho}\chi(\rho,y)\zeta,\bfn_2(y)\rangle+[\rho'(\chi(\rho,y))
   \partial_{\rho}\chi(\rho,y)\zeta+\zeta(\chi(\rho,y))]\langle\bfn_1(\chi(\rho,y)),\bfn_2(y)\rangle
\\ [0.1cm]
   &+\rho(\chi(\rho,y))\langle\bfn_1'(\chi(\rho,y))\partial_{\rho}\chi(\rho,y)\zeta,\bfn_2(y)\rangle, \quad y\in S_2, \;\;
   \rho\in\mathcal{O}_1.
\end{array}
\end{equation}
  In what follows we prove:
\begin{equation}
  [\rho\mapsto[y\mapsto\partial_{\rho}\sigma(\rho,y)]]\in C(\mathcal{O}_1\cap C^{m+\mu+1}(S_1),C^{m+\mu}(S_2,L(C^{m+\mu}(S_1),{\mathbf{R}}))).
\end{equation}
  We first prove:
\begin{equation}
  [\rho\mapsto[y\mapsto\partial_{\rho}\chi(\rho,y)]]\in C(\mathcal{O}_1\cap C^{m+\mu+1}(S_1),C^{m+\mu}(S_2,L(C^{m+\mu}(S_1),{\mathbf{R}}^n))).
\end{equation}
  Indeed, letting $G:\mathcal{O}_1\times S_1\to S_2$ be the mapping given by the right-hand side of (2.4), we have
$$
  y=G(\rho,\chi(\rho,y)), \quad y\in S_2,\;\; \rho\in\mathcal{O}_1.
$$
  Differentiating both sides of the above equation in $\rho$, we get
\begin{equation}
  \partial_{\rho}\chi(\rho,y)=-[\partial_xG(\rho,\chi(\rho,y))]^{-1}\partial_{\rho}G(\rho,\chi(\rho,y)).
\end{equation}
  We have
\begin{equation}
  \partial_xG(\rho,x)\xi=\Pi_2'(x+\rho(x)\bfn_1(x))\{\xi+[\rho'(x)\xi]\bfn_1(x)+\rho(x)\bfn_1'(x)\xi\}, \quad \forall\xi\in T_x(S_1),
\end{equation}
\begin{equation}
  \partial_{\rho}G(\rho,x)\zeta=\Pi_2'(x+\rho(x)\bfn_1(x))[\zeta(x)\bfn_1(x)], \quad \forall\zeta\in C^{m+\mu}(S_1).
\end{equation}
  From (2.11) it is easy to see
$$
  [\rho\mapsto[(x,\xi)\mapsto(G(\rho,x),\partial_xG(\rho,x)\xi)]]\in C(\mathcal{O}_1\cap C^{m+\mu+1}(S_1),C^{m+\mu}(T(S_1),T(S_2))),
$$
  which combined with the fact $\chi\in C^{m+\mu}(\mathcal{O}_1\times S_2,S_1)$ implies
\begin{equation}
   [\rho\mapsto[(y,\eta)\mapsto(\chi(\rho,y),[\partial_xG(\rho,\chi(\rho,y))]^{-1}\eta)]]\in
   C(\mathcal{O}_1\cap C^{m+\mu+1}(S_1),C^{m+\mu}(T(S_2),T(S_1))).
\end{equation}
  From (2.12) it is also easy to see
\begin{equation}
   [\rho\mapsto[(x,\zeta)\mapsto(G(\rho,x),\partial_{\rho}G(\rho,x)\zeta)]]\in C(\mathcal{O}_1,C^{m+\mu}(S_1\times C^{m+\mu}(S_1),T(S_2))).
\end{equation}
  Combining (2.10), (2.13), (2.14) and using the fact $\chi\in C^{m+\mu}(\mathcal{O}_1\times S_2,S_1)\subseteq C(\mathcal{O}_1,C^{m+\mu}(S_2,S_1))$
  we see that (2.9) follows. From (2.7), (2.9) and the fact $\chi\in C^{m+\mu}(\mathcal{O}_1\times S_2,S_1)\subseteq C(\mathcal{O}_1,
  C^{m+\mu}(S_2,S_1))$ we obtain (2.8). Hence
$$
  [\rho\mapsto(\psi\circ\varphi^{-1})'(\rho)]=[\rho\mapsto[\zeta\mapsto\partial_{\rho}\sigma(\rho,\cdot)\zeta]]\in
  C(\mathcal{O}_1\cap C^{m+\mu+1}(S_1),L(C^{m+\mu}(S_1),C^{m+\mu}(S_2))).
$$
  This proves the desired assertion. $\quad\Box$
\medskip

  The above lemma ensures that the following definition makes sense:
\medskip

  {\bf Definition 2.4}\ \ {\em Let $\Omega\in\mathfrak{D}^{m+\mu+1}({\mathbf{R}}^n)$, and regard it as a point in $\mathfrak{D}^{m+\mu}({\mathbf{R}}^n)$.
  We have the following concepts:

  $(1)$\ \ For a function $F:\mathcal{O}\to{\mathbf{R}}$ defined in a neighborhood $\mathcal{O}\subseteq\mathfrak{D}^{m+\mu}({\mathbf{R}}^n)$ of $\Omega$,
  we say that $F$ is differentiable at $\Omega$ if for any regular local chart $(\mathcal{U},\varphi)$ of $\mathfrak{D}^{m+\mu}({\mathbf{R}}^n)$ at
  $\Omega$, the function $F\circ\varphi^{-1}:\varphi(\mathcal{O}\cap\mathcal{U})\to{\mathbf{R}}$ is differentiable at $\Omega$. We denote by
  $\mathscr{D}^1_{\Omega}$ the set of all functions $F:\mathfrak{D}^{m+\mu}({\mathbf{R}}^n)\to{\mathbf{R}}$ which are differentiable at $\Omega$.

  $(2)$\ \ Let $f:(-\varepsilon,\varepsilon)\to\mathfrak{D}^{m+\mu}({\mathbf{R}}^n)$ $(\varepsilon>0)$ be a curve in $\mathfrak{D}^{m+\mu}({\mathbf{R}}^n)$
  passing $\Omega$, i.e., $f(0)=\Omega$. We say that $f(t)$ is differentiable at $t=0$ if for any regular local chart $(\mathcal{U},\varphi)$ of
  $\mathfrak{D}^{m+\mu}({\mathbf{R}}^n)$ at $\Omega$, the function $t\mapsto\varphi(f(t))$ is differentiable at $t=0$. Moreover, we define the tangent
  vector of this curve at $\Omega$, or the derivative $f'(0)$ of $f(t)$ at $t=0$, to be the mapping $f'(0):\mathscr{D}^1_{\Omega}\to{\mathbf{R}}$
  defined by
$$
  f'(0)F=(F\circ f)'(0), \quad \forall F\in\mathscr{D}^1_{\Omega}.
$$

  $(3)$\ \ We denote
$$
  T_{\Omega}(\mathfrak{D}^{m+\mu}({\mathbf{R}}^n))=\{f'(0):f:(-\varepsilon,\varepsilon)\to\mathfrak{D}^{m+\mu}({\mathbf{R}}^n),f(0)=\Omega,
  \mbox{$f(t)$ is differentiable at $t=0$}\},
$$
  and call it the tangent space of $\mathfrak{D}^{m+\mu}({\mathbf{R}}^n)$ at $\Omega$.}
\medskip

  For $\Omega\in\mathfrak{D}^{m+\mu+1}({\mathbf{R}}^n)$, the tangent space $T_{\Omega}(\mathfrak{D}^{m+\mu}({\mathbf{R}}^n))$ can be expressed in a
  different form which is very useful from viewpoint of applications. Indeed, since $\Omega\in\mathfrak{D}^{m+\mu+1}({\mathbf{R}}^n)$ implies that
  $\Gamma:=\partial\Omega$ is a $C^{m+\mu+1}$-hypersurface and its normal field $\bfn$ is of $C^{m+\mu}$-class: $\bfn\in C^{m+\mu}(\Gamma,
  {\mathbf{R}}^n)$, by using $\Gamma$ as a base hypersurface and repeating the argument before, we get a (irregular) local chart of
  $\mathfrak{D}^{m+\mu}({\mathbf{R}}^n)$ at $\Omega$, which we denote as $(\mathcal{U}_{\Omega},\varphi_{\Omega})$ and call as the {\em standard
  local chart} of $\mathfrak{D}^{m+\mu}({\mathbf{R}}^n)$ at $\Omega$. Clearly the projection $\Pi_{\Omega}:\mathcal{R}_{\Omega}\to\Gamma$ is of
  $C^{m+\mu}$-class: $\Pi_{\Omega}\in C^{m+\mu}(\mathcal{R}_{\Omega},\Gamma)$, where $\mathcal{R}_{\Omega}$ is a neighborhood of $\Gamma$ similar to
  $\mathcal{R}$ for the regular local chart. Using these facts and a similar argument as in the proof of Lemma 2.3, we see that for any regular local
  chart $(\mathcal{U},\varphi)$ of $\mathfrak{D}^{m+\mu}({\mathbf{R}}^n)$ at $\Omega$ with a smooth base hypersurface $S$, the coordinate transformation
  mappings $\varphi_{\Omega}\circ\varphi^{-1}$ and $\varphi\circ\varphi_{\Omega}^{-1}$ are differentiable at the points $\varphi(\Omega)$ and
  $\varphi_{\Omega}(\Omega)=0$, respectively (note that $\varphi(\Omega)\in C^{m+\mu+1}(S)$ and, since $\varphi_{\Omega}(\Omega)=0$, the term
  $\bfn'$ (cf. (2.11)) does not appear in computation of $(\varphi\circ\varphi_{\Omega}^{-1})'(\varphi_{\Omega}(\Omega))$).
  It follows that for a curve $f:(-\varepsilon,\varepsilon)\to\mathfrak{D}^{m+\mu}({\mathbf{R}}^n)$ $(\varepsilon>0)$ in
  $\mathfrak{D}^{m+\mu}({\mathbf{R}}^n)$ passing $\Omega$, $f(t)$ is differentiable at $t=0$ if and only if the function $t\mapsto\varphi_{\Omega}(f(t))$
  is differentiable at $t=0$. Since in this case we have $(\varphi_{\Omega}\circ f)'(0)\in C^{m+\mu}(\Gamma)$ and, furthermore, for any $\xi\in
  C^{m+\mu}(\Gamma)$ there exists a curve $f:(-\varepsilon,\varepsilon)\to\mathfrak{D}^{m+\mu}({\mathbf{R}}^n)$ $(\varepsilon>0)$ in
  $\mathfrak{D}^{m+\mu}({\mathbf{R}}^n)$ passing $\Omega$ which is differentiable at $t=0$ such that $(\varphi_{\Omega}\circ f)'(0)=\xi$ (choose, for
  instance, $f(t)=\varphi_{\Omega}^{-1}(t\xi)$ for $t\in(-\varepsilon,\varepsilon)$, or $f(t,x)=x+t\xi(x)\bfn(x)$ for $x\in\Gamma$ and $t\in
  (-\varepsilon,\varepsilon)$), we have
\begin{equation}
  T_{\Omega}(\mathfrak{D}^{m+\mu}({\mathbf{R}}^n))\approx C^{m+\mu}(\Gamma).
\end{equation}
  Hence, the tangent space $T_{\Omega}(\mathfrak{D}^{m+\mu}({\mathbf{R}}^n))$ of $\mathfrak{D}^{m+\mu}({\mathbf{R}}^n)$ at a point $\Omega\in
  \mathfrak{D}^{m+\mu+1}({\mathbf{R}}^n)$ can be alternatively defined to be the Banach space $C^{m+\mu}(\Gamma)$:
  $T_{\Omega}(\mathfrak{D}^{m+\mu}({\mathbf{R}}^n))=C^{m+\mu}(\Gamma)$, but for this purpose we need to identify any $Q\in\mathcal{U}_{\Omega}$ with
  $\varphi_{\Omega}(Q)$ ($f'(0)$ is then identified with $(\varphi_{\Omega}\circ f)'(0)$). Note that in this definition the
  tangent vector $f'(0)$ of a curve $f:(-\varepsilon,\varepsilon)\to\mathfrak{D}^{m+\mu}({\mathbf{R}}^n)$ $(\varepsilon>0)$ has a simple physical
  explanation: Regarding $f$ as a flow of simple $C^{m+\mu}$-domains,  $f'(0)$ is the normal velocity of the boundary $\partial\Omega=\Gamma$. To
  see this let $f:(-\varepsilon,\varepsilon)\to\mathfrak{D}^{m+\mu}({\mathbf{R}}^n)$ be a curve in $\mathfrak{D}^{m+\mu}({\mathbf{R}}^n)$ passing
  $\Omega$, differentiable at $t=0$. Let $f(t)=\Omega_t$, $\forall t\in (-\varepsilon,\varepsilon)$ (so that $\Omega_0=\Omega$). There exists a
  function $\rho:(-\varepsilon,\varepsilon)\times\Gamma\to{\mathbf{R}}$ such that
$$
  \partial\Omega_t=\{x+\rho(t,x)\bfn(x):x\in\Gamma\}, \quad \forall t\in (-\varepsilon,\varepsilon); \quad \rho(0,\cdot)=0.
$$
  Clearly, $(\varphi_{\Omega}\circ f)(t)=\rho(t,\cdot)$, $\forall t\in (-\varepsilon,\varepsilon)$. Hence $(\varphi_{\Omega}\circ f)'(0)=
  \partial_t\rho(0,\cdot)$, which is clearly the normal velocity of $\partial\Omega$.
\medskip

  From the proof of Lemma 2.3 we see that if $(\mathcal{U}_1,\varphi_1)$, $(\mathcal{U}_2,\varphi_2)$ are two regular local charts of
  $\mathfrak{D}^{m+\mu}({\mathbf{R}}^n)$ such that $\mathcal{U}_1\cap\mathcal{U}_2\neq\varnothing$, then for any positive integer $k$ and $\Omega\in
  \mathcal{U}_1\cap\mathcal{U}_2\cap\mathfrak{D}^{m+\mu+k}({\mathbf{R}}^n)$, $\varphi_2\circ\varphi_1^{-1}$ is $k$-th order differentiable at $\rho=
  \varphi_1(\Omega)$: $(\varphi_2\circ\varphi_1^{-1})^{(k)}(\rho)\in L^k(C^{m+\mu}(S_1),C^{m+\mu}(S_2))$, and
$$
\begin{array}{rl}
  [\rho\mapsto (\varphi_2\circ\varphi_1^{-1})^{(j)}(\rho)]\in & C^{k-j}(\mathcal{O}_1\cap C^{m+\mu+k}(S_1),L^j(C^{m+\mu}(S_1),C^{m+\mu}(S_2))),\\
  & \qquad\qquad\qquad j=1,2,\cdots,k,
\end{array}
$$
  where $\mathcal{O}_i=\varphi_i(\mathcal{U}_1\cap\mathcal{U}_2)\subseteq C^{m+\mu}(S_i)$, $i=1,2$, and, for two Banach spaces $X$ and $Y$, the
  notation $L^k(X,Y)$ denotes the Banach space of all bounded symmetric $k$-linear operators defined in $X$ and valued in $Y$. It follows that for
  a mapping $f:(-\varepsilon,\varepsilon)\to\mathfrak{D}^{m+\mu}({\mathbf{R}}^n)$ $(\varepsilon>0)$, if $f(t)\in\mathfrak{D}^{m+\mu+k}({\mathbf{R}}^n)$,
  $\forall t\in (-\varepsilon,\varepsilon)$, then we can define $k$-th order differentiability of $f$ as follows: Assume that $\varepsilon$ is
  sufficiently small such that $f(-\varepsilon,\varepsilon)$ is contained in $\varphi(\mathcal{U})$ for some regular local chart
  $(\mathcal{U},\varphi)$ with a smooth base hypersurface $S$. Then $f\in C^k((-\varepsilon,\varepsilon),\mathfrak{D}^{m+\mu}({\mathbf{R}}^n))$ if and
  only if
$$
  \varphi\circ f\in C((-\varepsilon,\varepsilon),C^{m+k+\mu}(S))\cap C^1((-\varepsilon,\varepsilon),C^{m+k-1+\mu}(S))\cap\cdots
  \cap C^k((-\varepsilon,\varepsilon),C^{m+\mu}(S)).
$$
  Later on we shall call a domain in $\mathfrak{D}^{m+\mu+k}({\mathbf{R}}^n)$ a {\em $C^k$-point} in $\mathfrak{D}^{m+\mu}({\mathbf{R}}^n)$.
\medskip

  {\bf Definition 2.5}\ \ {\em Let $X$ be a Banach space and $F:\mathcal{O}\subseteq\mathfrak{D}^{m+\mu}({\mathbf{R}}^n)\to X$, where $\mathcal{O}$ is
  an open subset of $\mathfrak{D}^{m+\mu}({\mathbf{R}}^n)$. Let $\Omega\in\mathcal{O}\cap\mathfrak{D}^{m+\mu+1}({\mathbf{R}}^n)$. We say that $F$ is
  differentiable at $\Omega$ if for a regular local chart $(\mathcal{U},\varphi)$ of $\mathfrak{D}^{m+\mu}({\mathbf{R}}^n)$ at $\Omega$ with a smooth
  base hypersurface $S$, the mapping $F\circ\varphi^{-1}:\varphi(\mathcal{U})\subseteq C^{m+\mu}(S)\to X$ is differentiable at $\varphi(\Omega)$, and
  define $F'(\Omega)\in L(T_{\Omega}(\mathfrak{D}^{m+\mu}({\mathbf{R}}^n)),X)$ as follows:
$$
  F'(\Omega)\xi=(F\circ f)'(0)\in X, \quad \forall\xi\in T_{\Omega}(\mathfrak{D}^{m+\mu}({\mathbf{R}}^n)),
$$
  where $f\in C((-\varepsilon,\varepsilon),\mathfrak{D}^{m+\mu}({\mathbf{R}}^n))\,(\varepsilon>0)$, $f(0)=\Omega$, $f(t)$ is differentiable at $t=0$
  and $f'(0)=\xi$. We denote $F\in\mathfrak{C}^1(\mathcal{O}\cap\mathfrak{D}^{m+\mu+1}({\mathbf{R}}^n),X)$ if $F$ is differentiable at every point
  $\Omega\in\mathcal{O}\cap\mathfrak{D}^{m+\mu+1}({\mathbf{R}}^n)$ and $(F\circ\varphi^{-1})'\in C(\varphi(\mathcal{O}\cap\mathcal{U}\cap
  \mathfrak{D}^{m+\mu+1}({\mathbf{R}}^n)),L(C^{m+\mu}(S),X))$ for any regular local chart $(\mathcal{U},\varphi)$ of $\mathfrak{D}^{m+\mu}({\mathbf{R}}^n)$
  with a smooth base hypersurface $S$.}
\medskip

  In particular, for any regular local chart $(\mathcal{U},\varphi)$ of $\mathfrak{D}^{m+\mu}({\mathbf{R}}^n)$ with a smooth base hypersurface $S$, the
  notation $\varphi'(\Omega)$ makes sense for any $C^1$-point $\Omega\in\mathcal{U}$, and $\varphi'(\Omega)\in
  L(T_{\Omega}(\mathfrak{D}^{m+\mu}({\mathbf{R}}^n)),C^{m+\mu}(S))$.

  For a mapping $F:\mathfrak{D}^{m+\mu}({\mathbf{R}}^n)\to\mathfrak{D}^{m+\mu}({\mathbf{R}}^n)$, a positive integer $k$ and a point $\Omega\in
  \mathfrak{D}^{m+\mu+k}({\mathbf{R}}^n)$, if for a regular local chart $(\mathcal{U}_1,\varphi_1)$ of $\mathfrak{D}^{m+\mu}({\mathbf{R}}^n)$ at $\Omega$
  with a smooth base hypersurface $S_1$ and a regular local chart $(\mathcal{U}_2,\varphi_2)$ of $\mathfrak{D}^{m+\mu}({\mathbf{R}}^n)$ at $F(\Omega)$
  with a smooth base hypersurface $S_2$ such that $F(\mathcal{U}_1)\subseteq\mathcal{U}_2$, the mapping $\varphi_2\circ F\circ\varphi_1^{-1}:
  \varphi_1(\mathcal{U}_1)\subseteq C^{m+\mu}(S_1)\to C^{m+\mu}(S_2)$ is $k$-th order differentiable at $\varphi_1(\Omega)$, then we say $F$ is $k$-th
  order differentiable at $\Omega$. We denote $F\in\mathfrak{C}^k(\mathfrak{D}^{m+\mu+k}({\mathbf{R}}^n),\mathfrak{D}^{m+\mu}({\mathbf{R}}^n))$ if $F$ is
  $k$-th order differentiable at every point $\Omega\in\mathfrak{D}^{m+\mu+k}({\mathbf{R}}^n)$, and
$$
\begin{array}{rl}
  [\rho\mapsto (\varphi_2\circ F\circ\varphi_1^{-1})^{(j)}(\rho)]\in & C^{k-j}(\varphi_1(\mathcal{U}_1)\cap C^{m+\mu+k}(S_1),
  L^j(C^{m+\mu}(S_1),C^{m+\mu}(S_2))),\\
  & \qquad\qquad\qquad j=1,2,\cdots,k
\end{array}
$$
  for any regular local chars $(\mathcal{U}_1,\varphi_1)$, $(\mathcal{U}_2,\varphi_2)$ of $\mathfrak{D}^{m+\mu}({\mathbf{R}}^n)$ such that
  $\varphi_1(\mathcal{U}_1)\cap C^{m+\mu+k}(S_1)\neq\varnothing$ and $F(\mathcal{U}_1)\subseteq\mathcal{U}_2$.

\subsection{Lie group actions to $\mathfrak{D}^{m+\mu}({\mathbf{R}}^n)$}

\hskip 2em
  It is clear that $\mathfrak{D}^{m+\mu}({\mathbf{R}}^n)$ is invariant under dilations, translations and rotations in ${\mathbf{R}}^n$. In what follows
  we study smoothness of these Lie group actions to $\mathfrak{D}^{m+\mu}({\mathbf{R}}^n)$.

  We first consider the action of the translation group. Let $G_{tl}={\mathbf{R}}^n$ be the additive group of $n$-vectors. Given $z\in{\mathbf{R}}^n$ and
  $\Omega\in\mathfrak{D}^{m+\mu}({\mathbf{R}}^n)$, let
$$
  p(z,\Omega)=\Omega+z=\{x+z:\,x\in\Omega\}.
$$
  It is clear that $p(z,\Omega)\in\mathfrak{D}^{m+\mu}({\mathbf{R}}^n)$, $\forall\Omega\in\mathfrak{D}^{m+\mu}({\mathbf{R}}^n)$, $\forall z\in{\mathbf{R}}^n$. Moreover, it is also
  clear that
$$
  p(0,\Omega)=\Omega, \quad p(z_1+z_2,\Omega)=p(z_1,p(z_2,\Omega)), \quad \forall\Omega\in\mathfrak{D}^{m+\mu}({\mathbf{R}}^n), \;\; \forall z_1,z_2\in
  {\mathbf{R}}^n.
$$
  Hence $(G_{tl},p)$ is a Lie group action on $\mathfrak{D}^{m+\mu}({\mathbf{R}}^n)$.
\medskip

  {\bf Lemma 2.6}\ \ {\em The Lie group action $(G_{tl},p)$ on $\mathfrak{D}^{m+\mu}({\mathbf{R}}^n)$ satisfies the following properties: For any
  nonnegative integers $k$ and $l$,
\begin{equation}
  [z\mapsto p(z,\cdot)]\in C^k(G_{tl},\mathfrak{C}^{\,l}(\mathfrak{D}^{m+k+l+\mu}({\mathbf{R}}^n),\mathfrak{D}^{m+\mu}({\mathbf{R}}^n))).
\end{equation}
  In particular, $p\in C(G_{tl}\times\mathfrak{D}^{m+\mu}({\mathbf{R}}^n),\mathfrak{D}^{m+\mu}({\mathbf{R}}^n))$.}
\medskip

  {\em Proof}.\ \ We first compute the representation of the mapping $p:G_{tl}\times\mathfrak{D}^{m+\mu}({\mathbf{R}}^n)\to
  \mathfrak{D}^{m+\mu}({\mathbf{R}}^n)$ in local charts of $G_{tl}\times\mathfrak{D}^{m+\mu}({\mathbf{R}}^n)$ and $\mathfrak{D}^{m+\mu}({\mathbf{R}}^n)$.
  Let $z_0\in G_{tl}$ and $\Omega\in\mathfrak{D}^{m+\mu}({\mathbf{R}}^n)$ be given. Choose a closed $C^{\infty}$ hypersurface $S\subseteq{\mathbf{R}}^n$
  sufficiently closed to $\partial\Omega$ such that it satisfies the conditions $(a)$--$(c)$ in the previous subsection. Let $(\mathcal{U},\varphi)$
  be the local chart of $\mathfrak{D}^{m+\mu}({\mathbf{R}}^n)$ at the point $\Omega$ as defined by (2.1) and (2.2). Let
$$
  \hat{\Omega}=\Omega+z_0, \quad  \hat{S}=S+z_0=\{x+z_0:x\in S\},
$$
  and $(\hat{\mathcal{U}},\hat{\varphi})$ be the local chart of $\mathfrak{D}^{m+\mu}({\mathbf{R}}^n)$ at the point $\hat{\Omega}$ as defined by (2.1)
  and (2.2), with $S$ there replaced by $\hat{S}$. Take $\varepsilon>0$ sufficiently small such that, by slightly shrinking the neighborhood
  $\mathcal{U}$ of $\Omega$ when necessary, we have $Q+z\in\hat{\mathcal{U}}$ for all $Q\in\mathcal{U}$ and $z\in B(z_0,\varepsilon)$. Here
  $B(z_0,\varepsilon)$ denotes the open sphere in $G_{tl}={\mathbf{R}}^n$ with center $z_0$ and radius $\varepsilon$. Now let $Q\in\mathcal{U}$ and
  $z\in B(z_0,\varepsilon)$. There exists a unique $\rho\in\mathcal{O}$, where $\mathcal{O}$ is a neighborhood of the origin in $C^{m+\mu}(S)$, such
  that
\begin{equation}
  \partial Q=\{x+\rho(x)\bfn(x):x\in S\},
\end{equation}
  where $\bfn$ is the outward unit normal field of $S$. By definition we have $\varphi(Q)=\rho$. Similarly, there exists a unique $\hat{\rho}\in
  \hat{\mathcal{O}}$, where $\hat{\mathcal{O}}$ is a neighborhood of the origin in $C^{m+\mu}(\hat{S})$, such that
\begin{equation}
  \partial(Q+z)=\{x+\hat{\rho}(x)\hat{\bfn}(x):x\in\hat{S}\},
\end{equation}
  where $\hat{\bfn}$ is the outward unit normal field of $\hat{S}$. By definition we have $\hat{\varphi}(Q+z)=\hat{\rho}$. Hence, letting
$$
  \tilde{p}(z,\rho)=(\hat{\varphi}\circ p)(z,\varphi^{-1}(\rho)), \quad \forall\rho\in\mathcal{O},\;\; \forall z\in B(z_0,\varepsilon),
$$
  we have
$$
  \tilde{p}(z,\rho)=\hat{\rho}, \quad \forall\rho\in\mathcal{O},\;\; \forall z\in B(z_0,\varepsilon).
$$
  Next, from (2.17) we have
$$
  \partial(Q+z)=\{x+z+\rho(x)\bfn(x):x\in S\}=\{y-z_0+z+\rho(y-z_0)\bfn(y-z_0):y\in\hat{S}\}.
$$
  Comparing this expression with (2.18) we get
\begin{equation}
  x+\hat{\rho}(x)\hat{\bfn}(x)=y+z-z_0+\rho(y-z_0)\bfn(y-z_0), \quad x,y\in\hat{S},
\end{equation}
  which, similarly as before, means that for any $x\in\hat{S}$ there exists a unique $y\in\hat{S}$ such that the above equality holds, and vice versa.
  It follows that there exists a function $y=\mu(\rho,z,x)$ uniquely determined by $(\rho,z)$, mapping $x\in\hat{S}$ to $y\in\hat{S}$,
  such that
\begin{equation}
  \hat{\rho}(x)=\langle\mu(\rho,z,x)-x+z-z_0,\hat{\bfn}(x)\rangle+\rho(\mu(\rho,z,x)-z_0)\langle\bfn(\mu(\rho,z,x)-z_0),\hat{\bfn}(x)\rangle,
  \quad x\in\hat{S}.
\end{equation}
  Similarly as in the proof of Theorem 2.2 (cf. (2.5)), we have
$$
  [(\rho,z)\mapsto\mu(\rho,z,\cdot)]\in C^k(\mathcal{O}\times B(z_0,\varepsilon),C^{m-k+\mu}(\hat{S},\hat{S})), \quad k=0,1,\cdots,m.
$$
  Using this fact, the expression (2.20) and a similar argument as in the proof of Lemma 2.3 we get
$$
\begin{array}{rl}
  \{z\mapsto[\rho\mapsto\partial_{\rho}^j\tilde{p}(z,\rho)]\}\in &C^k(B(z_0,\varepsilon),C^{\,l-j}(\mathcal{O}\cap C^{m+k+l+\mu}(S),
  L^j(C^{m+\mu}(S),C^{m+\mu}(\hat{S})))),\\
  & \qquad\qquad\qquad  k,l=0,1,\cdots,\;\; j=0,1,\cdots,l.
\end{array}
$$
  This proves (2.16). $\quad\Box$
\medskip

  Next we consider the action of the dilation group. Let $G_{\!dl}={\mathbf{R}}_+=(0,\infty)$ be the multiplicative group of all positive numbers. Given
  $\lambda\in{\mathbf{R}}_+$ and $\Omega\in\mathfrak{D}^{m+\mu}({\mathbf{R}}^n)$, let
$$
  q(\lambda,\Omega)=\lambda\Omega=\{\lambda x:\,x\in\Omega\}.
$$
  Clearly $q(\lambda,\Omega)\in\mathfrak{D}^{m+\mu}({\mathbf{R}}^n)$, $\forall\Omega\in\mathfrak{D}^{m+\mu}({\mathbf{R}}^n)$, $\forall\lambda\in
  {\mathbf{R}}_+$. Moreover, it is also clear that
$$
  q(1,\Omega)=\Omega, \quad q(\lambda_1\lambda_2,\Omega)=q(\lambda_1,q(\lambda_2,\Omega)), \quad \forall\Omega\in\mathfrak{D}^{m+\mu}({\mathbf{R}}^n),
   \;\; \forall\lambda_1,\lambda_2\in{\mathbf{R}}_+.
$$
  Hence $(G_{\!dl},p)$ is a Lie group action on $\mathfrak{D}^{m+\mu}({\mathbf{R}}^n)$. Similar to Lemma 2.6 we have
\medskip

  {\bf Lemma 2.7}\ \ {\em The Lie group action $(G_{\!dl},q)$ in $\mathfrak{D}^{m+\mu}({\mathbf{R}}^n)$ satisfies the following properties: For any
  nonnegative integers $k$ and $l$,
$$
  [\lambda\mapsto q(\lambda,\cdot)]\in C^k(G_{\!dl},\mathfrak{C}^{\,l}(\mathfrak{D}^{m+k+l+\mu}({\mathbf{R}}^n),\mathfrak{D}^{m+\mu}({\mathbf{R}}^n))).
$$
  In particular, $q\in C(G_{\!dl}\times\mathfrak{D}^{m+\mu}({\mathbf{R}}^n),\mathfrak{D}^{m+\mu}({\mathbf{R}}^n))$.}
\medskip

  The group actions $(G_{tl},p)$ and $(G_{\!dl},q)$ to $\mathfrak{D}^{m+\mu}({\mathbf{R}}^n)$ are not commutative. However, it is clear that the
  following relation holds:
\begin{equation}
  q(\lambda,p(z,\Omega))=p(\lambda z, q(\lambda,\Omega)), \quad \forall\Omega\in\mathfrak{D}^{m+\mu}({\mathbf{R}}^n),\;\;
  \forall z\in G_{tl},\;\; \forall\lambda\in G_{\!dl}.
\end{equation}
  Hence the actions $p$ of $G_{tl}$ and $q$ of $G_{\!dl}$ to $\mathfrak{D}^{m+\mu}({\mathbf{R}}^n)$ are quasi-commutative. Besides, denoting
$$
  g(z,\lambda,\Omega)=p(z, q(\lambda,\Omega)), \quad \forall\Omega\in\mathfrak{D}^{m+\mu}({\mathbf{R}}^n),\;\;
  \forall z\in G_{tl},\;\; \forall\lambda\in G_{\!dl},
$$
  we easily see that the following relation holds:
\begin{equation}
  {\rm rank}\,\partial_{(z,\lambda)}g(z,\lambda,\Omega)=n+1, \quad \forall\Omega\in\mathfrak{D}^{m+1+\mu}({\mathbf{R}}^n),\;\;
  \forall z\in G_{tl},\;\; \forall\lambda\in G_{\!dl}.
\end{equation}

  Finally, let $O(n)$ be the Lie group of all $n\times n$ orthogonal matrices. For any $A\in O(n)$ and $\Omega\in\mathfrak{D}^{m+\mu}({\mathbf{R}}^n)$,
  let
$$
  r(A,\Omega)=\mathbf{A}(\Omega)=\{\mathbf{A}(x):\,x\in\Omega\},
$$
  where $\mathbf{A}$ denotes the orthogonal transformation in ${\mathbf{R}}^n$ induced by $A$, i.e, $\mathbf{A}(x)=(Ax^T)^T$ by regarding vectors in
  ${\mathbf{R}}^n$ as $1\times n$ matrix. Clearly $r(A,\Omega)\in\mathfrak{D}^{m+\mu}({\mathbf{R}}^n)$, $\forall\Omega\in\mathfrak{D}^{m+\mu}({\mathbf{R}}^n)$,
  $\forall A\in O(n)$. Moreover, it is also
  clear that
$$
  r(I,\Omega)=\Omega, \quad r(A_1A_2,\Omega)=r(A_1,r(A_2,\Omega)), \quad \forall\Omega\in\mathfrak{D}^{m+\mu}({\mathbf{R}}^n), \;\;
  \forall A_1,A_2\in O(n).
$$
  Hence $(O(n),r)$ is a Lie group action on $\mathfrak{D}^{m+\mu}({\mathbf{R}}^n)$. Similar to Lemmas 2.6 and 2.7 we have
\medskip

  {\bf Lemma 2.8}\ \ {\em Let $m,n\geqslant 2$ be integers and $0\leqslant\mu\leqslant 1$. The Lie group action $(O(n),q)$ in
  $\mathfrak{D}^{m+\mu}({\mathbf{R}}^n)$ satisfies the following properties: For any nonnegative integers $k$ and $l$,
$$
  [A\mapsto r(A,\cdot)]\in C^k(O(n),\mathfrak{C}^{\,l}(\mathfrak{D}^{m+k+l+\mu}({\mathbf{R}}^n),\mathfrak{D}^{m+\mu}({\mathbf{R}}^n))).
$$
  In particular, $r\in C(O(n)\times\mathfrak{D}^{m+\mu}({\mathbf{R}}^n),\mathfrak{D}^{m+\mu}({\mathbf{R}}^n))$.}
\medskip

  The proof is similar to that of Lemma 2.6; we omit it. $\quad\Box$

  As we shall not use the group action $(O(n),r)$ later on, we do not make further discussion to it here.

  Note that all discussions made in this section apply to $\dot{\mathfrak{D}}^{m+\mu}({\mathbf{R}}^n)$ when all $C^{m+\mu}$-spaces are replaced with
  corresponding $\dot{C}^{m+\mu}$-spaces. Here we do not repeat the details. Later on we shall use such results without further discussion.

\section{Parabolic differential equations in Banach manifolds}
\setcounter{equation}{0}

\hskip 2em
  In this section we abstractly study parabolic differential equations in Banach manifolds. In particular, we shall prove the linearized stability
  theorem for parabolic differential equations in Banach manifolds which are invariant or quasi-invariant properties under a finite number of Lie
  group actions. Here we particularly mention that, since the Banach manifold $\mathfrak{D}^{m+\mu}({\mathbf{R}}^n)$ is not differentiable, for the
  purpose of applications of the abstract result to free boundary problems, in this section differential equations will be considered in
  non-differentiable Banach manifolds.

\subsection{Basic concepts}

\hskip 2em
  Let $X$ and $X_0$ be two Banach spaces. We say $X_0$ is an embedded Banach subspace of $X$ if $X_0$ is a linear subspace of $X$ and the
  norm of $X$ is majorized by that of $X_0$ when restricted to $X_0$, i.e., there exists a constant $C>0$ such that
$$
  \|x\|_{X}\leqslant C\|x\|_{X_0}, \qquad \forall x\in X_0.
$$
  In what follows, for open subset $U$ of $X_0$, any subset $E$ of $X$ and positive integer $k$, we use the notation $\mathfrak{C}^k(U,E)$ to denote
  the set of mappings $F:U\to E$ which are $k$-th order differentiable at every point $x\in U$ in the topology of $X$ and, furthermore,
$$
  [x\mapsto F^{(j)}(x)]\in C^{k-j}(U,L^j(X,X)), \quad j=1,2,\cdots,k.
$$
  Note that the condition $F\in\mathfrak{C}^k(U,E)$ is stronger than the condition $F\in{C}^k(U,E)$; the latter means that $[x\mapsto F^{(j)}(x)]\in
  C^{k-j}(U,L^j(X_0,X))$, $j=1,2,\cdots,k$.

\medskip
  {\bf Definition 3.1}\ \ {\em Let $\mathfrak{M}$ and $\mathfrak{M}_0$ be two topological $($i.e., they need not have differentiable structure$)$
  Banach manifolds built on Banach spaces $X$ and $X_0$, respectively. Let $k$ be a positive integer. We say $\mathfrak{M}_0$ is a $C^k$-embedded
  Banach submanifold of $\mathfrak{M}$ if the following four conditions are satisfied:
\begin{enumerate}
\item[]$(D1)$\ \ $X_0$ is an embedded Banach subspace of $X$.\vspace*{-0.2cm}
\item[]$(D2)$\ \ $\mathfrak{M}_0$ is a topological subspace of $\mathfrak{M}$.\vspace*{-0.2cm}
\item[]$(D3)$\ \ For any point $\eta\in\mathfrak{M}_0$ there exists a local chart $(\mathcal{U},\varphi)$ of $\mathfrak{M}$ at $\eta$ and a
  neighborhood $\mathcal{U}_0$ of $\eta$ in $\mathfrak{M}_0$ such that $\mathcal{U}_0\subseteq\mathcal{U}$ and
  $(\mathcal{U}_0,\varphi|_{\mathcal{U}_0})$ is a local chart of $\mathfrak{M}_0$ at $\eta$. We call such a local chart $(\mathcal{U},\varphi)$ of
  $\mathfrak{M}$ at $\eta$ a {\em $\mathfrak{M}_0$-regular local chart} of $\mathfrak{M}$ at $\eta$.  \vspace*{-0.2cm}
\item[]$(D4)$\ \ For any point $\eta\in\mathfrak{M}_0$, if $(\mathcal{U},\varphi)$ and $(\mathcal{V},\psi)$ are two $\mathfrak{M}_0$-regular local
  charts of $\mathfrak{M}$ at $\eta$ with corresponding neighborhoods $\mathcal{U}_0$, $\mathcal{V}_0$ of $\eta$ in $\mathfrak{M}_0$ such that
  $\mathcal{U}_0\subseteq\mathcal{U}$, $\mathcal{V}_0\subseteq\mathcal{V}$ and $(\mathcal{U}_0,\varphi|_{\mathcal{U}_0})$,
  $(\mathcal{V}_0,\varphi|_{\mathcal{V}_0})$ are local charts of $\mathfrak{M}_0$ at $\eta$, respectively, then the following relations hold:
$$
  \psi\circ\varphi^{-1}\in\mathfrak{C}^k(\varphi(\mathcal{U}_0\cap\mathcal{V}_0),\psi(\mathcal{U}\cap\mathcal{V})) \quad \mbox{and} \quad
  \varphi\circ\psi^{-1}\in\mathfrak{C}^k(\psi(\mathcal{U}_0\cap\mathcal{V}_0),\varphi(\mathcal{U}\cap\mathcal{V})).
$$
\end{enumerate}
}

  Let $\mathfrak{M}$ be a Banach manifold and $\mathfrak{M}_0$ a $C^k$-embedded Banach submanifold of $\mathfrak{M}$, $k\geqslant 1$.
  It follows that for any $\eta\in\mathfrak{M}_0$ and a curve $f:(-\varepsilon,\varepsilon)\to\mathfrak{M}$ ($\varepsilon>0$) passing $\eta$:
  $f(0)=\eta$, the concept of $f(t)$ being differentiable at $t=0$ makes sense and is defined to be that $\varphi\circ f:(-\varepsilon,\varepsilon)
  \to X$ is differentiable at $t=0$ for any $\mathfrak{M}_0$-regular local chart $(\mathcal{U},\varphi)$ of $\mathfrak{M}$ at $\eta$. The tangent
  vector $f'(0)$ of this curve at $\eta$, or derivative of $f(t)$ at $t=0$, is defined similarly as in Definition 2.4. The tangent space
  $\mathcal{T}_{\eta}(\mathfrak{M})$ of $\mathfrak{M}$ at $\eta$ is defined as follows:
$$
  \mathcal{T}_{\eta}(\mathfrak{M})=\{f'(0):f:(-\varepsilon,\varepsilon)\to\mathfrak{M},\; f(0)=\eta,\; \mbox{$f(t)$ is differentiable at $t=0$}\}.
$$
  If $(\mathcal{U},\varphi)$ is a $\mathfrak{M}_0$-regular local chart of $\mathfrak{M}$ at $\eta\in\mathfrak{M}_0$, then the derivative of
  $\varphi: \mathcal{U}\to X$ at $\eta$ also makes sense and is defined to be a linear mapping $\varphi'(\eta):\mathcal{T}_{\eta}(\mathfrak{M})\to X$
  such that
$$
  \varphi'(\eta)\upsilon=\varphi'(\eta)f'(0)=(\varphi\circ f)'(0) \quad \mbox{for}\;\, \upsilon=f'(0)\in \mathcal{T}_{\eta}(\mathfrak{M}).
$$
  We endow $\mathcal{T}_{\eta}(\mathfrak{M})$ with the topology naturally induced by the strong topology of $X$, i.e., induced by the norm
$$
  \|\upsilon\|_{\varphi}=\|\varphi'(\eta)\upsilon\|_X, \quad \forall\upsilon\in\mathcal{T}_{\eta}(\mathfrak{M}).
$$
  It follows that for any $\eta\in\mathfrak{M}_0$, $\mathcal{T}_{\eta}(\mathfrak{M})$ is a Banach space isomorphic to $X$. Let
$$
  \mathcal{T}_{\mathfrak{M}_0}(\mathfrak{M})=\bigcup_{\eta\in\mathfrak{M}_0}\mathcal{T}_{\eta}(\mathfrak{M})\;(\mbox{disjoint union})\;
  :=\{(\eta,\upsilon):\eta\in\mathfrak{M}_0,\upsilon\in\mathcal{T}_{\eta}(\mathfrak{M})\}.
$$
  We can use the standard method to make $\mathcal{T}_{\mathfrak{M}_0}(\mathfrak{M})$ into a (topological) Banach space built on the Banach space
  $X_0\times X$. More precisely, for any $(\eta_0,\upsilon_0)\in\mathcal{T}_{\mathfrak{M}_0}(\mathfrak{M})$, let $(\mathcal{U},\varphi)$ be a
  $\mathfrak{M}_0$-regular local chart of $\mathfrak{M}$ at $\eta_0$ and $\mathcal{U}_0$ a neighborhood of $\eta_0$ in $\mathfrak{M}_0$ such that
  $\mathcal{U}_0\subseteq\mathcal{U}$ and $(\mathcal{U}_0,\varphi|_{\mathcal{U}_0})$ is a local chart of $\mathfrak{M}_0$ at $\eta$.
  Let
$$
  \mathcal{O}=\{(\eta,\upsilon):\eta\in\mathcal{U}_0,\upsilon\in\mathcal{T}_{\eta}(\mathfrak{M})\},
$$
  and define $\Phi:\mathcal{O}\to X_0\times X$ as follows:
$$
  \Phi(\eta,\upsilon)=(\varphi(\eta),\varphi'(\eta)\upsilon), \quad \forall (\eta,\upsilon)\in\mathcal{O}.
$$
  We use $(\mathcal{O},\Phi)$ as a local chart of $\mathcal{T}_{\mathfrak{M}_0}(\mathfrak{M})$ at the point $(\eta_0,\upsilon_0)$. It is not hard to
  check that this indeed makes $\mathcal{T}_{\mathfrak{M}_0}(\mathfrak{M})$ into a (topological) Banach manifold. We call
  $\mathcal{T}_{\mathfrak{M}_0}(\mathfrak{M})$ endowed with this topological structure the tangent bundle of $\mathfrak{M}_0$ in $\mathfrak{M}$.
  Later on we often regard $\mathcal{T}_{\mathfrak{M}_0}(\mathfrak{M})$ as the set $\{\upsilon:\upsilon\in\mathcal{T}_{\eta}(\mathfrak{M}),\eta\in
  \mathcal{U}_0\}$, i.e., we identify the point $(\eta,\upsilon)$ with $\upsilon$, because for different $\eta,\xi\in\mathfrak{M}_0$, the sets
  $\mathcal{T}_{\eta}(\mathfrak{M})$, $\mathcal{T}_{\xi}(\mathfrak{M})$ do not have common point.

  Let $\mathfrak{M}$ and $\mathfrak{M}_0$ be as above. Given an open interval $I\subseteq {\mathbf{R}}$, we use the notation $C^1\!(I,\mathfrak{M}_0)$
  to denote the set of mappings $f:I\to\mathfrak{M}_0$ such that for any $t_0\in I$, there exists a $\mathfrak{M}_0$-regular local chart $(\mathcal{U},
  \varphi)$ of $\mathfrak{M}$ at $f(t_0)$ and $\varepsilon>0$ sufficiently small such that $\varphi\circ f\in C((t_0-\varepsilon,t_0+\varepsilon),X_0)
  \cap C^1((t_0-\varepsilon,t_0+\varepsilon),X)$. Note that this in particular implies $f\in C(I,\mathfrak{M}_0)\subseteq C(I,\mathfrak{M})$ and
  $f(t)$ is differentiable at every point $t\in I$.
  Let $\mathscr{F}$ be a mapping from $\mathfrak{M}_0$ to $\mathcal{T}_{\mathfrak{M}_0}(\mathfrak{M})$ such that for any $\eta\in\mathfrak{M}_0$,
  $\mathscr{F}(\eta)\in\mathcal{T}_{\eta}(\mathfrak{M})$. We call $\mathscr{F}$ a vector field in $\mathfrak{M}$ with domain $\mathfrak{M}_0$. The
  purpose of this section is to study the following differential equation in the Banach manifold $\mathfrak{M}$:
\begin{equation}
   \eta'=\mathscr{F}(\eta),
\end{equation}
  or more precisely, the initial value problem of the above equation:
\begin{equation}
\left\{
\begin{array}{rlll}
   \eta'(t)&=&\mathscr{F}(\eta(t)), &\quad t>0,\\
   \eta(0)&=&\eta_0, &
\end{array}
\right.
\end{equation}
  where $\eta_0$ is a given point in $\mathfrak{M}$. By a solution of the equation (3.1) in an open interval $I\subseteq{\mathbf{R}}$ we mean a
  function $\eta\in C^1\!(I,\mathfrak{M}_0)$ such that $\eta'(t)=\mathscr{F}(\eta(t))$ for all $t\in I$, and by a solution of the problem (3.2)
  in a interval $[0,T)$ ($0<T\leqslant\infty$) we mean a function $\eta\in C([0,T),\mathfrak{M})\cap C^1((0,T),\mathfrak{M}_0)$ such that $\eta(0)=
  \eta_0$ and $\eta'(t)=\mathscr{F}(\eta(t))$ for all $t\in (0,T)$. Note that these conditions imply that $\eta_0\in\overline{\mathfrak{M}}_0$, where
  $\overline{\mathfrak{M}}_0$ denotes the closure of $\mathfrak{M}_0$ in $\mathfrak{M}$.

  Recall that a differential equation $x'=F(x)$ in  a Banach space $X$, where $F\in C^1(O,X)$ with $O$ being an open subset of an embedded Banach
  subspace $X_0$ of $X$, is said to be of parabolic type if for any $x\in O$, $F'(x)$ is a sectorial operator in $X$ with domain $X_0$, and
  the graph norm of $X_0={\rm Dom}\,F'(x)$ is equivalent to the norm of $X_0$, i.e., there exist positive constants $C_1,C_2$ such that
$$
  C_1\|y\|_{X_0}\leqslant\|y\|_{X}+\|F'(x)y\|_{X}\leqslant C_2\|y\|_{X_0}, \quad \forall y\in X_0.
$$
  We introduce the following concept:

\medskip
  {\bf Definition 3.2}\ \ {\em Let $\mathfrak{M}$ be a Banach manifold and $\mathfrak{M}_0$ a $C^k$-embedded Banach submanifold of $\mathfrak{M}$,
  $k\geqslant 1$. Let $\mathscr{F}$ be a vector field in $\mathfrak{M}$ with domain $\mathfrak{M}_0$. We say the differential equation $(3.1)$
  is of parabolic type if for any $\eta_0\in\mathfrak{M}_0$ there exists a $\mathfrak{M}_0$-regular local chart of $\mathfrak{M}$ at $\eta_0$ such
  that its representation in that local chart is of parabolic type.}

\medskip
  More precisely, for $\eta_0\in\mathfrak{M}_0$, let $(\mathcal{U},\varphi)$ be a $\mathfrak{M}_0$-regular local chart of $\mathfrak{M}$ at $\eta_0$
  and $\mathcal{U}_0$ a neighborhood of $\eta_0$ in $\mathfrak{M}_0$ such that $\mathcal{U}_0\subseteq\mathcal{U}$ and
  $(\mathcal{U}_0,\varphi|_{\mathcal{U}_0})$ is a local chart of $\mathfrak{M}_0$ at $\eta_0$. The representation of the vector field
  $\mathscr{F}$ in the local chart $(\mathcal{U},\varphi)$ is the mapping $F:\varphi(\mathcal{U}_0)\subseteq X_0\to X$ defined as follows:
$$
  F(x)=\varphi'(\varphi^{-1}(x))\mathscr{F}(\varphi^{-1}(x)), \quad x\in\varphi(\mathcal{U}_0).
$$
  The representation of the differential equation (3.1) in the local chart $(\mathcal{U},\varphi)$ is the following differential equation in $X$:
\begin{equation}
  x'=F(x).
\end{equation}
  Definition 3.1 says that the differential equation $(3.1)$ in the Banach manifold $\mathfrak{M}$ is of parabolic type if for any $\eta_0\in
  \mathfrak{M}_0$ there exists a $\mathfrak{M}_0$-regular local chart of $\mathfrak{M}$ at $\eta_0$ such that the differential equation (3.3) in the
  Banach space $X$ is of parabolic type. Note that this does not exclude the possibility that a representation of the equation (3.1) in some other
  $\mathfrak{M}_0$-regular local chart is not of parabolic type or even the representation $F$ of the vector field $\mathscr{F}$ is not differentiable.
  However, if $k\geqslant 2$ then if the representation of (3.1) in one $\mathfrak{M}_0$-regular local chart of a point $\eta_0\in\mathfrak{M}_0$ is
  of parabolic type then its representation in any other $\mathfrak{M}_0$-regular local chart of $\eta_0$ is also of parabolic type. To see this
  let $(\mathcal{V},\psi)$ be another $\mathfrak{M}_0$-regular local chart of $\mathfrak{M}$ at $\eta_0$, with $\mathcal{V}_0$
  being a neighborhood of $\eta_0$ in $\mathfrak{M}_0$ such that $\mathcal{V}_0\subseteq\mathcal{V}$ and $(\mathcal{V}_0,\psi|_{\mathcal{V}_0})$ is
  a local chart of $\mathfrak{M}_0$ at $\eta_0$. Let
$$
  F_U(x)=\varphi'(\varphi^{-1}(x))\mathscr{F}(\varphi^{-1}(x)), \quad x\in\varphi(\mathcal{U}_0),
$$
$$
  F_V(y)=\psi'(\psi^{-1}(y))\mathscr{F}(\psi^{-1}(y)), \quad y\in\psi(\mathcal{V}_0).
$$
  It is clear that
$$
  F_V(y)=\Psi(y)F_U(\Phi(y)), \quad y\in\psi(\mathcal{U}_0\cap\mathcal{V}_0),
$$
  where $\Psi(y)=(\psi\circ\varphi^{-1})'(\Phi(y))\in L(X)$ and $\Phi(y)=\varphi(\psi^{-1}(y))$ (so that $\Phi:\psi(\mathcal{U}_0\cap\mathcal{V}_0)
  \to\varphi(\mathcal{U}_0\cap\mathcal{V}_0)$). Since $k\geqslant 2$, we have $\Psi\in C^1(\psi(\mathcal{U}_0\cap\mathcal{V}_0),L(X))$ and, clearly,
\begin{equation}
   F_V'(y)z=\Psi(y)F_U'(\Phi(y))\Psi(y)^{-1}z+\{\Psi'(y)z\}F_U(\Phi(y)),  \quad y\in\psi(\mathcal{U}\cap\mathcal{V})\cap X_0, \;\; z\in X_0.
\end{equation}
  It is easy to see that if $F_U'(x)$ is a sectorial operator in $X$ (for any $x\in\varphi(\mathcal{U}_0)$) then $\Psi(y)F_U'(\Phi(y))\Psi(y)^{-1}$
  is also a sectorial operator in $X$ (for any $y\in\psi(\mathcal{U}_0\cap\mathcal{V}_0)$). Moreover, it is clear that for any $y\in\psi(\mathcal{U}_0
  \cap\mathcal{V}_0)$, the operator $z\mapsto \{\Psi'(y)z\}F_U(\Phi(y))$ is a bounded linear operator in $X$. Hence, by a standard perturbation
  theorem, it follows that $F_V'(\eta)$ is also a sectorial operator in $X$.

  Recall that for two Banach spaces $X,Y$ and an open subset $U$ of $X$, the notation $C^{2-0}(U,Y)$ denotes the set of all $C^1$-mappings $F:U\to Y$
  such that the mapping $x\mapsto F'(x)$ from $U$ to $L(X,Y)$ is locally Lipschitz continuous: For any $x_0\in U$ there exist corresponding constant
  $C=C(x_0)>0$ and neighborhood $B(x_0,\varepsilon)\subseteq U$ ($\varepsilon>0$) of $x_0$ such that
$$
  \|F'(x)-F'(y)\|_{L(X,Y)}\leqslant C\|x-y\|_X, \quad \forall x,y\in B(x_0,\varepsilon).
$$
  Local theory of parabolic differential equations in Banach manifolds can be easily established by using local chart to transform the problem into
  corresponding problem in Banach spaces. For instance, we have the following basic result:

\medskip
  {\em Assume that the differential equation $(3.1)$ is of parabolic type and for any $\eta_0\in\mathfrak{M}_0$ there exists a $\mathfrak{M}_0$-regular
  local chart $(\mathcal{U},\varphi)$ of $\mathfrak{M}$ at $\eta_0$ such that the representation $F$ of $\mathscr{F}$ in this local chart belongs to
  $C^{2-0}(\varphi(\mathcal{U}_0),X)$, where $\mathcal{U}_0$ is a neighborhood of $\eta_0$ in $\mathfrak{M}_0$ such that $\mathcal{U}_0\subseteq
  \mathcal{U}$ and $(\mathcal{U}_0,\varphi|_{\mathcal{U}_0})$ is a local chart of $\mathfrak{M}_0$ at $\eta_0$. Assume further that $X_0$ is dense in
  $X$. Then for any $\eta_0\in\mathfrak{M}_0$ there exists a corresponding $\delta>0$ such that the problem $(3.2)$ has a unique solution $\eta\in
  C([0,\delta),\mathfrak{M}_0)\cap C^1((0,\delta),\mathfrak{M}_0)$.}

\medskip
  As usual a point $\eta^*\in\mathfrak{M}_0$ such that $\mathscr{F}(\eta^*)=0$ is called a singular point of the vector field $\mathscr{F}$
  and a stationary point or stationary solution of the equation (3.1). Asymptotic stability of stationary solutions of differential
  equations in Banach manifolds can be defined by using local charts, and the linearized stability criterion for asymptotic stability of isolated
  stationary solutions can also be stated and proved by using local charts. We omit these discussions here.

\subsection{Invariant and quasi-invariant differential equations in Banach manifolds}

\hskip 2em
  We now turn to consider asymptotic stability of stationary solutions of invariant and quasi-invariant parabolic differential equations in Banach
  manifolds.

  Let $\mathfrak{M}$ be a (topological or $C^0$) Banach manifold and $G$ a Lie group of dimension $n$. An {\em action} of $G$ to $\mathfrak{M}$ is a
  mapping $p:G\times\mathfrak{M}\to\mathfrak{M}$ satisfying the following three conditions:
\begin{enumerate}
\item[]$(L1)$\ \ $p\in C(G\times\mathfrak{M},\mathfrak{M})$.\vspace*{-0.2cm}
\item[]$(L2)$\ \ $p(e,\eta)=\eta$, $\forall\eta\in\mathfrak{M}$, where $e$ denotes the unit element of $G$, and
$$
  p(a,p(b,\eta))=p(ab,\eta), \quad \forall a,b\in G, \;\; \forall\eta\in\mathfrak{M}.
$$
\item[]$(L3)$\ \ If $a,b\in G$ are such that $p(a,\eta)=p(b,\eta)$ for some $\eta\in\mathfrak{M}$, then $a=b$.\vspace*{-0.2cm}
\end{enumerate}
  Assume further that $\mathfrak{M}$ has a $C^k$-embedded Banach submanifold $\mathfrak{M}_0$, $k\geqslant 1$. We say that the Lie group action
  $(G,p)$ to $\mathfrak{M}$ is {\em $\mathfrak{M}_0$-regular} if the following three additional conditions are also satisfied:
\begin{enumerate}
\item[]$(L4)$\ \ $p(a,\mathfrak{M}_0)\subseteq\mathfrak{M}_0$, $\forall a\in G$, and $p\in C(G\times\mathfrak{M}_0,\mathfrak{M}_0)$.\vspace*{-0.2cm}
\item[]$(L5)$\ \ For any $a\in G$, the mapping $\eta\mapsto p(a,\eta)$ is differentiable at every point $\eta\in\mathfrak{M}_0$ ($\Rightarrow
  \partial_{\eta}p(a,\eta)\in L(T_{\eta}(\mathfrak{M}),T_{p(a,\eta)}(\mathfrak{M}))$, $\forall\eta\in\mathfrak{M}_0$), and $[(a,(\eta,\upsilon))
  \mapsto (p(a,\eta),\partial_{\eta}p(a,\eta)\upsilon)]\in C(G\times \mathcal{T}_{\mathfrak{M}_0}(\mathfrak{M}),\mathcal{T}_{\mathfrak{M}_0}(\mathfrak{M}))$.\vspace*{-0.2cm}
\item[]$(L6)$\ \ For any $\eta\in\mathfrak{M}_0$, the mapping $a\mapsto p(a,\eta)$ is differentiable at every point $a\in G$ ($\Rightarrow
  \partial_ap(a,\eta)\in L(T_{a}(G),T_{p(a,\eta)}(\mathfrak{M}))$, $\forall a\in G$), $[((a,z),\eta)\mapsto(p(a,\eta),\partial_ap(a,\eta)z)]\in
  C(T(G)\times\mathfrak{M}_0, \mathcal{T}_{\mathfrak{M}_0}(\mathfrak{M}))$, and
$$
  {\rm rank}\,\partial_{a}p(a,\eta)=n, \quad \forall a\in G, \;\; \forall\eta\in\mathfrak{M}_0.
$$
\end{enumerate}

  {\bf Definition 3.3}\ \ {\em Let $\mathfrak{M}$ be a Banach manifold and $\mathfrak{M}_0$ a $C^k$-embedded Banach submanifold of $\mathfrak{M}$,
  $k\geqslant 1$. Let $\mathscr{F}$ be a vector field in $\mathfrak{M}$ with domain $\mathfrak{M}_0$. Let $(G,p)$ be a $\mathfrak{M}_0$-regular Lie
  group action to $\mathfrak{M}$. We say $\mathscr{F}$ is quasi-invariant under the Lie group action $(G,p)$ if there exists a positive-valued
  function $\theta$ defined in $G$ such that the following condition is satisfied:
$$
  \mathscr{F}(p(a,\eta))=\theta(a)\partial_{\eta}p(a,\eta)\mathscr{F}(\eta), \quad \forall a\in G, \;\; \forall\eta\in\mathfrak{M}_0.
$$
  In this case we also say $\mathscr{F}$ is $\theta$-quasi-invariant and call $\theta$ quasi-invariance factor, and also say the differential
  equation $(3.1)$ is $\theta$-quasi-invariant under the Lie group action $(G,p)$. If in particular $\theta(a)=1$, $\forall a\in G$, then we simply
  say the vector field $\mathscr{F}$ and the differential equation $(3.1)$ are invariant under the Lie group action $(G,p)$.}
\medskip

  Clearly, $\theta$ is a group homomorphism from $G$ to the multiplicative group ${\mathbf{R}}_+$: For any $a,b\in G$, $\theta(ab)=\theta(a)\theta(b)$.
  The following result is obvious:
\medskip

  {\em Let $\mathscr{F}$ be a vector field in $\mathfrak{M}$ with domain $\mathfrak{M}_0$. Assume that $\mathscr{F}$ is $\theta$-quasi-invariant
  under the Lie group action $(G,p)$. If $t\mapsto\eta(t)$ is a solution of the differential equation $(3.1)$, then for any $a\in G$,
  $t\mapsto p(a,\eta(t\theta(a)))$ is also a solution of this equation, and if $\eta^*\in\mathfrak{M}_0$ is a stationary point of $(3.1)$, then
  for any $a\in G$, $p(a,\eta^*)$ is also a stationary point of $(3.1)$.} $\quad\Box$

\medskip
  As a consequence, no stationary point of a quasi-invariant differential equation is isolated, and if $\eta^*\in\mathfrak{M}_0$ is a stationary
  point of (3.1), then all points in the $n$-dimensional manifold $\{p(a,\eta^*):a\in G\}$ are stationary points of (3.1).

  Inspired by potential applications to free boundary problems such as (1.1), we consider a general situation where the vector field $\mathscr{F}$ is
  quasi-invariant under a number of Lie group actions $(G_i,p_i)$, $i=1,2,\cdots,N$, with possibly different quasi-invariance factors $\theta_i$,
  $i=1,2,\cdots,N$, respectively. We assume that the combined action of these Lie group actions satisfies the following conditions:
\begin{enumerate}
\item[]$(L7)$\ \ For any $1\leqslant i,j\leqslant N$ with $i\neq j$ there exists corresponding smooth function $f_{ij}:G_i\times G_j\to G_i$ such
  that
$$
   p_j(b,p_i(a,\eta))=p_i(f_{ij}(a,b),p_j(b,\eta)), \quad \forall\eta\in\mathfrak{M},\;\; \forall a\in G_i,\;\; \forall b\in G_j.
$$
\item[]$(L8)$\ \  Let $g:G\times\mathfrak{M}:=(G_1\times G_2\times\cdots\times G_N)\times\mathfrak{M}\to\mathfrak{M}$ be the function
$$
   g(a,\eta)=p_1(a_1,p_2(a_2,\cdots,p_N(a_N,\eta)))
$$
  for $\eta\in\mathfrak{M}$ and $a=(a_1,a_2,\cdots,a_N)\in G:=G_1\times G_2\times\cdots\times G_N$. Then
$$
   {\rm rank}\,\partial_ag(a,\eta)=n_1+n_2+\cdots+n_N,
$$
  $\forall\eta\in\mathfrak{M}$, $\forall a\in G$, where $n_i=\dim G_i$, $i=1,2,\cdots,N$.
\end{enumerate}

  The main result of this section is as follows:
\medskip

  {\bf Theorem 3.4}\ \ {\em Let $\mathfrak{M}$ and $\mathfrak{M}_0$ be two Banach manifolds built on the Banach spaces $X$ and $X_0$, respectively,
  such that $\mathfrak{M}_0$ is a $C^k$-embedded Banach submanifold of $\mathfrak{M}$, $k\geqslant 1$, where $X_0$ is a densely embedded Banach
  subspace of $X$. Let $(G_i,p_i)$, $i=1,2,\cdots,N$, be a finite number of $\mathfrak{M}_0$-regular Lie group actions to $\mathfrak{M}$ satisfying
  conditions $(L7)$ and $(L8)$. Let $\mathscr{F}$ be a vector field in $\mathfrak{M}$ with domain $\mathfrak{M}_0$, and $\eta^*\in
  \mathfrak{M}_0$ a singular point of $\mathscr{F}$. Assume that $\mathscr{F}$ is quasi-invariant under all Lie group actions $(G_i,p_i)$,
  $i=1,2,\cdots,N$. Assume further that there exists a $\mathfrak{M}_0$-regular local chart $(\mathcal{U},\varphi)$ of $\mathfrak{M}$ at $\eta^*$
  such that the representation $F$ of $\mathscr{F}$ in this local chart satisfies the following four conditions:
\begin{enumerate}
\item[]$(G_1)$\ \ The differential equation $x'=F(x)$ in $X$ is of parabolic type.\vspace*{-0.2cm}
\item[]$(G_2)$\ \ $F\in C^{2-0}(O,X)$, where $O=\varphi(\mathcal{U}_0)$ and $\mathcal{U}_0$ is a neighborhood of $\eta^*$ in
  $\mathfrak{M}_0$ as in Definition 3.1 $(ii)$.\vspace*{-0.2cm}
\item[]$(G_3)$\ \ Let $x^*=\varphi(\eta^*)$. Then $\dim{\rm Ker}F'(x^*)=n:=n_1+n_2+\cdots+n_N$, ${\rm Range}F'(x^*)$ is closed in $X$, and
\begin{equation}
  X={\rm Ker}F'(x^*)\oplus{\rm Range}F'(x^*).
\end{equation}
\item[]$(G_4)$\ \ $\sup\{{\rm Re}\lambda:\lambda\in\sigma(F'(x^*))\backslash\{0\}\}<0$.
\end{enumerate}
  Then we have the following assertions:

  $(1)$\ \ The set $\mathcal{M}_c=\{g(a,\eta^*):a\in G:=G_1\times G_2\times\cdots\times G_N\}$ is a $n$-dimensional submanifold of $\mathfrak{M}_0$.

  $(2)$\ \ There is a neighborhood $\mathcal{O}$ of $\mathcal{M}_c$ in $\mathfrak{M}_0$ such that for any $\eta_0\in\mathcal{O}$, the initial value
  problem $(3.2)$ has a unique solution $\eta\in C([0,\infty),\mathfrak{M}_0)\cap C^1((0,\infty),\mathfrak{M}_0)$.

  $(3)$\ \ There exists a submanifold $\mathcal{M}_s$ of $\mathfrak{M}_0$ of codimension $n$ passing $\eta^*$ such that for any $\eta_0\in
  \mathcal{M}_s$, the solution of the problem $(3.2)$ satisfies $\displaystyle\lim_{t\to\infty}\eta(t)=\eta^*$ and, conversely, if the solution of
  $(3.2)$ satisfies this property then $\eta_0\in\mathcal{M}_s$.

  $(4)$\ \ For any $\eta_0\in\mathcal{O}$ there exist unique $a\in G$ and $\xi_0\in\mathcal{M}_s$ such that $\eta_0=g(a,\xi_0)$ and, for the solution
  $\eta=\eta(t)$ of $(3.2)$,}
\begin{equation}
   \lim_{t\to\infty}\eta(t)=g(a,\eta^*).
\end{equation}

  To prove this theorem we need a preliminary lemma. Recall that for a given Banach space $X$ and given numbers $0<\alpha<1$ and $T>0$, the notation
  $C^{\alpha}_{\alpha}((0,T],X)$ denotes the Banach space of bounded vector functions $u:(0,T]\to X$ such that the vector function
  $t\mapsto t^{\alpha}u(t)$ is uniformly $\alpha$-H\"{o}lder continuous in $(0,T]$, with norm
$$
  \|u\|_{C^{\alpha}_{\alpha}((0,T],X)}=\sup_{0<t\leqslant T}\|u(t)\|_X+\sup_{0<s<t\leqslant T}
  \frac{\|t^{\alpha}u(t)-s^{\alpha}u(s)\|_X}{(t-s)^{\alpha}}
$$
  (cf. the introduction of Chapter 4 of \cite{Lun2}), and for given $\omega>0$, the notation $C^{\alpha}([T,\infty),X,-\omega)$ denotes the Banach
  space of vector functions $u:[T,\infty)\to X$ such that the vector function $t\mapsto e^{\omega t}u(t)$ is uniformly $\alpha$-H\"{o}lder continuous
  in $([T,\infty)$, with norm
$$
  \|u\|_{C^{\alpha}([T,\infty),X,-\omega)}=\sup_{t\geqslant T}\|e^{\omega t}u(t)\|_X+\sup_{t>s\geqslant T}
  \frac{\|e^{\omega t}u(t)-e^{\omega s}u(s)\|_X}{(t-s)^{\alpha}}
$$
  (cf. Section 4.4 of \cite{Lun2}).
\medskip

  {\bf Lemma 3.5}\ \ {\em Let $X$ be a Banach spaces and $X_0$ an embedded Banach subspace of $X$. Let $A$ be a sectorial operator in $X$ with
  domain $X_0$. Assume $\omega_-=-\sup\{{\rm Re}\lambda:\lambda\in\sigma(A)\}>0$ and $f\in C^{\alpha}_{\alpha}((0,1],X)\cap
  C^{\alpha}([1,\infty),X,-\omega)$, where  $0<\alpha<1$ and $\omega\in (0,\omega_-)$. Given $u_0\in X_0$, let $u(t)=e^{tA}u_0+\int_0^t
  e^{(t-s)A}f(s)ds$. Then $u\in C^{\alpha}_{\alpha}((0,1],X_0)\cap C^{\alpha}([1,\infty),X_0,-\omega)$, and there exists a constant $C=
  C(\alpha,\omega)>0$ such that}
\begin{equation}
   \|u\|_{C^{\alpha}_{\alpha}((0,T],X_0)}+\|u\|_{C^{\alpha}([T,\infty),X_0,-\omega)}\leqslant
   C(\|u_0\|_{X_0}+\|f\|_{C^{\alpha}_{\alpha}((0,T],X)}+\|f\|_{C^{\alpha}([T,\infty),X,-\omega)}).
\end{equation}

  {\em Proof}.\ \ See Lemma 2.2 of \cite{Cui2}. $\quad\Box$
\medskip

  {\em Proof of Theorem 3.4}.\ \ The assertion (1) is an immediate consequence of the properties $(L4)$ and $(L6)$ of the Lie group actions
  $(G_i,p_i)$, $i=1,2,\cdots,N$, and the condition $(L8)$. In what follows we prove the assertions (2) $\sim$ (4). For simplicity of notations we
  only consider the case $N=2$; for general $N$ the proof is similar. Besides, we assume $x^*=0$, so that $F(0)=0$.

  Firstly, by Theorem 8.1.1 of \cite{Lun2}, there exists a neighborhood $U_0$ of the origin of $X_0$ such that for any $x_0\in U_0$, the initial
  value problem
\begin{equation}
\left\{
\begin{array}{ll}
  x'(t)=F(x(t)), &\quad t>0,\\
  x(0)=x_0
\end{array}
\right.
\end{equation}
  has a unique local solution $x\in C([0,T],X_0)\cap C^1([0,T],X)\cap C^{\alpha}_{\alpha}((0,T],X_0)$, where $T=T(x_0)$ depends on $x_0$ and $\alpha$
  is an arbitrary number in $(0,1)$. Furthermore, denoting by $T^*(x_0)$ the supreme of all such $T$, then by Proposition 9.1.1 of \cite{Lun2} there
  exists $\varepsilon>0$ independent of $x_0$ such that if $\|x(t)\|_{X_0}<\varepsilon$ for all $t\in [0,T^*(x_0))$, then $T^*(x_0)=\infty$.

  In what follows we denote $M_c=\{\varphi(\eta):\eta\in \mathcal{M}_c\cap\mathcal{U}_0\cap\varphi^{-1}(U_0)\}$. Note that $0\in M_c$.

  Let $G_i^0$ be a small neighborhood of the unit element $e$ in $G_i$, $i=1,2$, and $\mathcal{U}'$, $\mathcal{U}'_0$ small neighborhoods of $\eta^*$
  in $\mathfrak{M}$ and $\mathfrak{M}_0$, respectively, such that $p_i(a_i,\eta)\in\mathcal{U},\mathcal{U}_0$, respectively, for all $a_i\in G_i^0$,
  $i=1,2$, and $\eta\in\mathcal{U}',\mathcal{U}'_0$, respectively. The Lie group actions $(G_1,p_1)$ and $(G_2,p_2)$ in the Banach manifold
  $\mathfrak{M}$ induce local Lie group actions $(G_1^0,\tilde{p}_1)$ and $(G_2^0,\tilde{p}_2)$, respectively, in the open subset
  $U'=\varphi(\mathcal{U}')$ of the Banach spaces $X$ as follows: For any $a_i\in G_i^0$ and $x\in U'$, define
$$
  \tilde{p}_i(a_i,x)=\varphi(p_i(a_i,\varphi^{-1}(x))), \quad i=1,2.
$$
  When restricted to the $U_0'=\varphi(\mathcal{U}_0')$, they are local Lie group actions in this open subset of $X_0$.
  A simple computation easily shows that for each $i=1,2$, $F$ is $\theta_i$-quasi-invariant with respect to the local Lie group
  action $(G_i^0,\tilde{p}_i)$, i.e.\footnotemark[2],
\footnotetext[2]{Recall that a mapping $F:O\subseteq X\to X$, where $X$ is a Banach space and $O$ is a subset of $X$ such that a local Lie group
  action $(G,p)$ is acted to $O$, is said to be invariant under the local Lie group action $(G,p)$ if it satisfies the relation $F(p(a,x))=F(x)$
  for all $a\in G$ and $x\in O$, and quasi-invariant (in the sense of \cite{Cui3}) under the local Lie group action $(G,p)$ if it satisfies the
  relation $F(p(a,x))=\partial_xp(a,x)F(x)$ for all $a\in G$ and $x\in O$. The relation (3.9) shows that if a vector field $\mathscr{F}$ in a Banach
  manifold $\mathfrak{M}$ is invariant under a Lie group action to $\mathfrak{M}$ then its representation in local chart is only quasi-invariant but
  not invariant under the representation of the Lie group action which is a local Lie group action to the base Banach space.}
\begin{equation}
  F(\tilde{p}_i(a_i,x))=\theta_i(a_i)\partial_{x}\tilde{p}_i(a_i,x)F(x), \quad
  \forall a_i\in G_i^0, \;\; \forall x\in U_0'.
\end{equation}

  Let $A=F'(0)$ and $N(x)=F(x)-Ax$. The $N\in C^2(U_0,X)$ and $N(0)=0$, $N'(0)=0$, so that by shrinking $U_0$ when necessary, there exists a
  constant $C>0$ such that
\begin{equation}
  \|N(x)\|_{X}\leqslant C\|x\|_{X_0}^2, \quad \|N(x)-N(y)\|_{X}\leqslant C(\|x\|_{X_0}+\|y\|_{X_0})\|x-y\|_{X_0},
  \quad \forall x,y\in U_0.
\end{equation}
  The equation $(3.8)_1$ can be rewritten as follows:
\begin{equation}
  x'=Ax+N(x).
\end{equation}
  Let $\sigma_-(A)=\sigma(A)\backslash\{0\}$. Choose a closed smooth curve in the complex plane such that it encloses the origin and separates it
  from $\sigma_-(A)$. We denote by $\Gamma$ this curve with anticlockwise orientation. Let $P\in L(X)$ be the following operator:
$$
  P=\frac{1}{2\pi i}\int_{\Gamma}R(\lambda,A){\rm d}\lambda.
$$
  We know that $P$ is a bounded projection: $P^2=P$. Since $0$ is the unique element in $\sigma(A)$ enclosed by $\Gamma$ and $X={\rm Ker}A\oplus
  {\rm Range}A$ with ${\rm Range}A$ closed, we have (cf. Proposition A.2.2 of \cite{Lun2})
$$
  PX=PX_0={\rm Ker}A, \quad (I-P)X={\rm Range}A.
$$
  Note that these relations imply that $APX=0$ and $A(I-P)X_0=(I-P)X$. Let $A_-=A|_{(I-P)X_0}$. Then $A_-:(I-P)X_0\to (I-P)X$ is an isomorphism and
  $\sigma(A_-)=\sigma_-(A)$, so that
\begin{equation}
  \sup\{{\rm Re}\lambda:\lambda\in\sigma(A_-)\}=-\omega_-<0.
\end{equation}
  Given $\delta,\delta'>0$, we denote by $B_1(0,\delta)$ and $B_2(0,\delta')$ spheres in $PX$ and $(I-P)X_0$, respectively, both centered at the
  origin but with radius $\delta$ and $\delta'$, respectively, i.e.,
$$
  B_1(0,\delta)=\{u\in PX=PX_0:\|u\|_{X_0}<\delta\}, \quad B_2(0,\delta')=\{v\in (I-P)X_0:\|v\|_{X_0}<\delta'\},
$$
  and let $O=\{u+v:u\in B_1(0,\delta),v\in B_2(0,\delta')\}\subseteq X_0$, $D=B_1(0,\delta)\times B_2(0,\delta')\subseteq PX\times(I-P)X_0$. We prove
  that if $\delta,\delta'>0$ are sufficiently small then there exists a mapping $h\in C^{2-0}(B_1(0,\delta),B_2(0,\delta'))$ such that
\begin{equation}
  M_c\cap O=\{u+h(u):u\in B_1(0,\delta)\}.
\end{equation}
  We first choose $\delta,\delta'>0$ small enough such that $O\subseteq U'$, and define $F_1:D\to PX$ and $F_2:D\to (I-P)X$ respectively by
$$
  F_1(u,v)=PF(u+v), \quad F_2(u,v)=(I-P)F(u+v), \quad \forall(u,v)\in D.
$$
  Clearly $F_2\in C^{2-0}(D,(I-P)X)$ (surely also $F_1\in C^{2-0}(D,PX)$) and $F_2(0,0)=0$. Since $\partial_vF_2(0,0)=A_-:(I-P)X_0\to (I-P)X$ is an
  isomorphism, by the implicit function theorem we see that if $\delta,\delta'$ are sufficiently small then there exists a unique mapping $h\in
  C^{2-0}(B_1(0,\delta),B_2(0,\delta'))$ such that $h(0)=0$, $F_2(u,h(u))=0$, $\forall u\in B_1(0,\delta)$, and, furthermore, for any $u\in
  B_1(0,\delta)$, $v=h(u)$ is the unique solution of the equation $F_2(u,v)=0$ in $B_2(0,\delta')$. For $x\in M_c\cap O$ let $u=Px$ and $v=(I-P)x$.
  Since $F(x)=0$, which implies $F_2(u,v)=(I-P)F(x)=0$, we infer that $v=h(u)$. Hence (3.13) is true. Note that this further implies that
$$
  F_1(u,h(u))=PF(u+h(u))=0, \quad \forall u\in B_1(0,\delta).
$$
  Note also that since $\partial_uF_2(0,0)PX=(I-P)APX_0=0$, we have $h'(0)=-[\partial_vF_2(0,0)]^{-1}\partial_uF_2(0,0)=0$.

  Let
$$
  N_1(u,v)=PN(u+v), \quad N_2(u,v)=(I-P)N(u+v), \quad \forall(u,v)\in D.
$$
  Since $AX_0=(I-P)X$, we have $PA=0$. Using this fact we see the differential equation (3.11) reduces into the following system of differential
  equations:
\begin{equation}
\left\{
\begin{array}{l}
  u'=N_1(u,v),\\
  v'=A_-v+N_2(u,v).
\end{array}
\right.
\end{equation}
  Given $x_0\in O$, set $u_0=Px_0$, $v_0=(I-P)x_0$, and let $(u,v)=(u(t),v(t))$ be the solution of (3.14) under the initial condition $(u(0),v(0))=
  (u_0,v_0)$ defined in a maximal interval $[0,T^*)$ such that $(u(t),v(t))\in D$ for all $t\in[0,T^*)$. Since $(u,v)\equiv(0,0)$ is a solution of
  (3.14) defined for all $t\geqslant 0$, by continuous dependence of the solution on initial data we infer that by replacing $\delta,\delta'$ with
  smaller numbers when necessary, we may assume that $T^*=T^*(x_0)>1$ for all $x_0\in D$. Let $\phi(t)=v(t)-h(u(t))$, $t\in [0,T^*)$. Since
  $A_-h(u)+N_2(u,h(u))=F_2(u,h(u))=0$ and $N_1(u,h(u))=F_1(u,h(u))=0$ for all $u\in B_1(0,\delta)$, we see that $\phi'(t)=A_-\phi(t)+\sigma(t)$,
  $\forall t\in [0,T^*)$, where
$$
  \sigma(t)=[N_2(u(t),v(t))-N_2(u(t),h(u(t)))]-h'(u(t))[N_1(u(t),v(t))-N_1(u(t),h(u(t)))].
$$
  By D'hamal formula, it follows that
$$
  \phi(t)={\rm e}^{tA_-}\phi(0)+\int_0^t{\rm e}^{(t-s)A_-}\sigma(s){\rm d}s, \quad t\in [0,T^*).
$$
  Using this expression and applying Lemma 3.5, we see that for any $0<\alpha<1$ and $\omega\in (0,\omega_-)$,
\begin{equation}
   \|\phi\|_{C^{\alpha}_{\alpha}((0,1],X_0)}+\|\phi\|_{C^{\alpha}([1,T^*),X_0,-\omega)}\leqslant
   C(\|\phi(0)\|_{X_0}+\|\sigma\|_{C^{\alpha}_{\alpha}((0,1],X)}+\|\sigma\|_{C^{\alpha}([1,T^*),X,-\omega)}).
\end{equation}
  From (3.10) it is not hard to deduce that (cf. the proof of Theorem 9.1.2 of \cite{Lun2})
$$
   \|\sigma\|_{C^{\alpha}_{\alpha}((0,1],X)}+\|\sigma\|_{C^{\alpha}([1,T^*),X,-\omega)}\leqslant
   C\max\{\delta,\delta'\}(\|\phi\|_{C^{\alpha}_{\alpha}((0,1],X_0)}+\|\phi\|_{C^{\alpha}([1,T^*),X_0,-\omega)}).
$$
  Substituting this estimate into (3.15) we see that if $\delta,\delta'$ are chosen sufficiently small then
$$
  \|\phi\|_{C^{\alpha}_{\alpha}((0,1],X_0)}+\|\phi\|_{C^{\alpha}([1,T^*),X_0,-\omega)}\leqslant C\|\phi(0)\|_{X_0},
$$
  which implies, in particular, that
\begin{equation}
  \|\phi(t)\|_{X_0}\leqslant C{\rm e}^{-\omega t}\|\phi(0)\|_{X_0}, \quad \forall t\in [0,T^*).
\end{equation}
  Next, since $F_1(u,h(u))=0$, $\forall u\in B_1(0,\delta)$, we have
$$
  u'(t)=N_1(u(t),v(t))-N_1(u(t),h(u(t)))=:\sigma_1(t), \quad \forall t\in [0,T^*).
$$
  Again from (3.10) we have
$$
  \|\sigma_1(t)\|_{X}\leqslant C(\|u(t)\|_{X_0}+\|v(t)\|_{X_0})\|\phi(t)\|_{X_0}\leqslant
   C\|\phi(t)\|_{X_0}, \quad \forall t\in [0,T^*).
$$
   Hence
\begin{eqnarray}
  \|u(t)\|_{X_0}&\leqslant& C\|u(t)\|_{X}\leqslant C(\|u_0\|_{X}+\int_0^t\|\sigma_1(s)\|_{X}{\rm d}s)
  \leqslant C(\|u_0\|_{X}+\int_0^t\|\phi(s)\|_{X_0}{\rm d}s) \nonumber \\
  &\leqslant& C(\|u_0\|_{X}+\|\phi(0)\|_{X_0}), \quad \forall t\in [0,T^*).
\end{eqnarray}
  In getting the first inequality we used the fact that $PX_0=PX={\rm Ker}A$ is a finite-dimensional space so that all norms in it are mutually
  equivalent. From (3.16) and (3.17) we see that if $\delta,\delta'$ are chosen sufficiently small then $\|x(t)\|_{X_0}<\varepsilon$ for all
  $t\in [0,T^*)$, so that $T^*=\infty$.

  Next, similarly as in the proof of (3.17) we see that for any $s>t\geqslant 0$,
\begin{equation}
  \|u(s)-u(t)\|_{X_0}\leqslant C\int_t^s\|\sigma_1(\tau)\|_{X}{\rm d}\tau\leqslant C\int_t^s\|\phi(\tau)\|_{X_0}{\rm d}\tau
  \leqslant C({\rm e}^{-\omega t}-{\rm e}^{-\omega s})\|\phi(0)\|_{X_0}.
\end{equation}
  Hence $\displaystyle\lim_{t\to\infty}u(t)$ exists in $X_0$, which we denote as $\bar{u}$. Let $\bar{x}=\bar{u}+h(\bar{u})$. Then
  $\bar{x}\in M_c$ and
\begin{equation}
  \lim_{t\to\infty}x(t)=\lim_{t\to\infty}u(t)+\lim_{t\to\infty}h(u(t))+\lim_{t\to\infty}\phi(t)=\bar{u}+h(\bar{u})=\bar{x}.
\end{equation}
  This shows that the solution of the equation (3.11) starting from an initial point located in the neighborhood $O$ of the origin converges to a
  point lying in $M_c$ as $t\to\infty$.

  Note that by letting $s\to\infty$ in (3.18), we have
\begin{equation}
  \|u(t)-\bar{u}\|_{X_0}\leqslant C{\rm e}^{-\omega t}\|\phi(0)\|_{X_0}, \quad \forall t\geqslant 0.
\end{equation}
  From (3.16) and (3.20) we easily obtain
\begin{equation}
  \|x(t)-\bar{x}\|_{X_0}\leqslant C{\rm e}^{-\omega t}\|\phi(0)\|_{X_0}, \quad \forall t\geqslant 0.
\end{equation}
  Hence, the solution of the equation (3.11) converges to its limit as $t\to\infty$ in exponential speed.

  We now prove that there exists a $C^{2-0}$-Banach manifold $M_s\subseteq O$ of codimension $n_1+n_2$, such that the solution $x=x(t)$ of the
  equation (3.19) converges to the origin as $t\to\infty$ if and only if $x(0)\in M_s$. To this end, for any $(u_0,v_0)\in D$ we denote by $u=
  u(t;u_0,v_0)$, $v=v(t;u_0,v_0)$ the unique solution of (3.22) with initial data $u(0)=u_0$, $v(0)=v_0$, and set $x(t;u_0,v_0)=u(t;u_0,v_0)+
  v(t;u_0,v_0)$, $\forall t\geqslant 0$. Since $x=x(t;u_0,v_0)$ is the unique solution of (3.11) with initial data $x(0;u_0,v_0)=x_0:=u_0+v_0$ and
$$
  {\rm e}^{tA}P=P, \quad {\rm e}^{tA}(I-P)={\rm e}^{tA_-}(I-P),
$$
  we have
\begin{eqnarray*}
  x(t;u_0,v_0)&=&\, {\rm e}^{tA}x_0+\int_0^t{\rm e}^{(t-\tau)A}N(x(\tau;u_0,v_0)){\rm d}\tau\\
  &=&\, u_0+{\rm e}^{tA_-}v_0+\int_0^tPN(x(\tau;u_0,v_0)){\rm d}\tau
  +\int_0^t{\rm e}^{(t-\tau)A_-}(I-P)N(x(\tau;u_0,v_0)){\rm d}\tau, \quad \forall t\geqslant 0.
\end{eqnarray*}
  Using (3.16), (3.18), (3.20) and (3.21) we see that the improper integral $\displaystyle\int_0^{\infty}PN(x(\tau;u_0,v_0)){\rm d}\tau$ is convergent
  and $\displaystyle\lim_{t\to\infty}x(t;u_0,v_0)=0$ if and only if the following relation holds:
$$
  u_0+\int_0^{\infty}PN(x(\tau;u_0,v_0)){\rm d}\tau=0.
$$
  By applying the implicit function theorem, we can easily show that if $\delta,\delta'$ are sufficiently small then this equation defines an implicit
  function $u_0=q(v_0)$, where $q\in C^{2-0}(B_2(0,\delta'),B_1(0,\delta))$ and $q(0)=0$, cf. the proof of Theorem 2.1 of \cite{Cui3} for details.
  Hence, the above equation defines a $C^{2-0}$-Banach manifold in $D$. We denote by $M_s$ the corresponding $C^{2-0}$-Banach manifold in $O$. Note
  that its definition ensures that $\displaystyle\lim_{t\to\infty}x(t;u_0,v_0)=0$ if and only if $x_0=u_0+v_0\in M_s$.

  We now prove that for any $x_0\in O$ there exist unique $a\in G_1$, $b\in G_2$ and $y_0\in M_s$ such that $x_0=\tilde{p}_1(a,\tilde{p}_2(b,y_0))$
  and
$$
  \lim_{t\to\infty}x(t)=\tilde{p}_1(a,\tilde{p}_2(b,0)),
$$
  where $x=x(t)$ is the solution of (3.11) with initial data $x(0)=x_0$. Indeed, let $\bar{x}$ be as in (3.19). Since $\bar{x}\in M_c$, there exist
  unique $a\in G_1^0$, $b\in G_2^0$ such that $\tilde{p}_1(a,\tilde{p}_2(b,0))=\bar{x}$. Let $y_0=\tilde{p}_2(b^{-1},\tilde{p}_1(a^{-1},x_0))$ and
$$
  y(t)=\tilde{p}_2(b^{-1},\tilde{p}_1(a^{-1},x(t\theta_1(a^{-1})\theta_2(b^{-1})))) \quad  \mbox{for}\;\, t\geqslant 0.
$$
  $y(t)$ is the solution of (3.19) with initial data $y(0)=y_0$. From the condition $\displaystyle\lim_{t\to\infty}x(t)=\bar{x}$ we have
$$
  \lim_{t\to\infty}y(t)=\tilde{p}_2(b^{-1},\tilde{p}_1(a^{-1},\bar{x}))=0.
$$
  Hence $y_0\in M_s$. This proves existence. Uniqueness is obvious.

  Now let $\mathcal{O}_{\eta^*}=\varphi^{-1}(O)$ and $\mathcal{M}_s=\varphi^{-1}(M_s)$. Then $\mathcal{O}_{\eta^*}$ is a neighborhood of $\eta^*$
  and $\mathcal{M}_s$ is a submanifold of $\mathfrak{M}_0$ satisfying the condition (3). Moreover, for any $\eta_0\in\mathcal{O}_{\eta^*}$
  the initial value problem (3.2) has a unique solution $\eta\in C([0,\infty),\mathfrak{M})\cap C^1((0,\infty),\mathfrak{M}_0)$ and there exist
  unique $a\in G_1$, $b\in G_2$ and $\xi_0\in\mathcal{M}_s$ such that $\eta_0=g(a,b,\xi_0)$ and the relation (3.6) holds. To finish the proof of
  Theorem 3.4 we now only need to repeat the above argument for every point $\eta\in\mathcal{M}_c$ and then glue all the open sets $\mathcal{O}_{\eta}$
  together to form the neighborhood $\mathcal{O}$ of $\mathcal{M}_c$. This completes the proof.  $\quad\Box$
\medskip

  {\em Remark}.\ \ $\mathcal{M}_c$ is called the {\em center manifold} of the equation (3.1), and $\mathcal{M}_s$ is called the {\em stable
  manifold} of the equation (3.1) corresponding to the stationary point $\eta^*$.

\section{Simple applications of Theorem 3.4}
\setcounter{equation}{0}

\hskip 2em
  The simplest application of Theorems 3.4 is to the following planer system:
\begin{equation}
  x'=0, \quad y'=-y.
\end{equation}
  This is a differential equation in the Euclidean space ${\mathbf{R}}^2$ (so that $\mathfrak{M}=\mathfrak{M}_0=X=X_0={\mathbf{R}}^2$). Let $F:{\mathbf{R}}^2
  \to{\mathbf{R}}^2$ be the mapping $F(x,y)=(0,-y)$, $\forall (x,y)\in{\mathbf{R}}^2$. Then the above system can be rewritten as the following differential
  equation in ${\mathbf{R}}^2$:
\begin{equation}
  u'=F(u).
\end{equation}
  Let $G={\mathbf{R}}$ be the usual one-dimensional additive Lie group. We introduce an action $p$ of $G$ to ${\mathbf{R}}^2$ as follows:
$$
  p(z,(x,y))=(x+z,y), \quad \forall (x,y)\in{\mathbf{R}}^2, \;\; \forall z\in G.
$$
  Since $\partial_up(z,u)=id$ and $F(p(z,u))=F(u)$, $\forall u\in{\mathbf{R}}^2$, $\forall z\in G$, we see the equation (4.2) is invariant in the group
  action $(G,p)$. Hence Theorem 3.4 (with $N=1$) applies to it (it is easy to verify that all conditions in Theorem 3.4 are satisfied). Clearly, the
  center manifold is the $x$-axis: $\mathcal{M}_c=\{(x,0):x\in{\mathbf{R}}\}$, and the stable manifold corresponding to every stationary point $(c,0)
  \in\mathcal{M}_c$ is the straight line passing it parallel to the $y$-axis; in particular, the stable manifold $\mathcal{M}_s$ corresponding to
  the stationary point $u^*=(0,0)$ is $y$-axis: $\mathcal{M}_s=\{(0,y):y\in{\mathbf{R}}\}$. For any $u_0=(x_0,y_0)\in{\mathbf{R}}^2$, let $z=x_0$ and
  $v_0=(0,y_0)$. Then  $v_0\in\mathcal{M}_s$, $u_0=p(z,v_0)$, and the solution $u(t)=(x_0,y_0{\rm e}^{-t})$ of (4.2) with initial data $u(0)=u_0$ can
  be expressed via group action as $u(t)=p(z,v(t))$, where $v(t)=(0,y_0{\rm e}^{-t})$ has initial data $v(0)=v_0$ and converges to the stationary point
  $u^*$ as $t\to\infty$, and $u(t)$ does not converge to $u^*$, but instead converges to the stationary point $(x_0,0)=p(z,u^*)$.

  The simplest application of Theorem 3.4 to partial differential equations is as follows. Let $\Omega$ be a given bounded domain in ${\mathbf{R}}^n$
  with a $C^2$-boundary. Consider the initial-boundary value problem
\begin{equation}
\left\{
\begin{array}{ll}
  \partial_tu(x,t)=\Delta u(x,t), &\quad x\in\Omega, \;t>0,\\
  \partial_{\bfn}u(x,t)=0, &\quad x\in\partial\Omega, \;t>0,\\
  u(x,0)=u_0(x), &\quad x\in\Omega,
\end{array}
\right.
\end{equation}
  where $\bfn$ denotes the outward-pointing unit normal field of $\partial\Omega$ and $u_0\in L^2(\Omega)$.
  It is well-known that stationary solutions of this problem is not isolated but make up an one-dimensional manifold $\mathcal{M}_c=
  {\mathbf{R}}{\mathbf{1}}$, where ${\mathbf{1}}$ denotes the function in $\overline{\Omega}$ with values identically $1$, and the solution of the
  above problem has the following asymptotic behavior as $t\to\infty$:
\begin{equation}
   \lim_{t\to\infty}u(x,t)=\frac{1}{|\Omega|}\int_{\Omega}u_0(x){\rm d}x \quad \mbox{in $L^2(\Omega)$ norm}.
\end{equation}
  To apply Theorem 3.4 to this problem we let $\mathfrak{M}=X=L^2(\Omega)$ and $\mathfrak{M}_0=X_0=\{u\in H^2(\Omega):
  \partial_{\bfn}u|_{\partial\Omega}=0\}$. Let $G={\mathbf{R}}$ be as before. We introduce an action $p$ of $G$ to $X=L^2(\Omega)$ as follows:
$$
  p(a,u)=a{\mathbf{1}}+u, \quad \forall u\in L^2(\Omega), \quad \forall a\in{\mathbf{R}}.
$$
  Note that the restriction of this group action to $X_0$ is also a group action to $X_0$. It is easy to see that the differential equation in the
  Banach space $X=L^2(\Omega)$ corresponding to the problem (4.3) is invariant under the Lie group action $(G,p)$ defined here, and Theorem 3.4
  applies to it. The stable manifold of this equation corresponding to the stationary point $u^*=0$ is
$$
   \mathcal{M}_s=\Big\{u\in X_0:\int_{\Omega}u(x){\rm d}x=0\Big\}.
$$
  For any $u_0\in X_0$, let $v_0=u_0-\displaystyle\Big(\frac{1}{|\Omega|}\int_{\Omega}u_0(x){\rm d}x\Big){\mathbf{1}}$ and $a=
  \displaystyle\frac{1}{|\Omega|}\int_{\Omega}u_0(x){\rm d}x$. Clearly $v_0\in\mathcal{M}_c$, $u_0=p(a,v_0)$, and $u(t)=p(a,v(t))$, $\forall t\geqslant
  0$, where $u(t)$ and $v(t)$ are solutions of (4.3) with respect to initial data $u_0$ and $v_0$, respectively. The last relation implies that
$$
   \lim_{t\to\infty}u(t)=p(a,0)=\Big(\frac{1}{|\Omega|}\int_{\Omega}u_0(x){\rm d}x\Big){\mathbf{1}} \quad \mbox{in $L^2(\Omega)$ norm},
$$
  which recovers the formula (4.4).

  To show an application of Theorem 3.4 in the case $N>1$ we consider the one phase Hele-Shaw problem in the whole space ${\mathbf{R}}^n$:
\begin{equation}
\left\{
\begin{array}{ll}
   \Delta u(x,t)=0, &\qquad x\in\Omega(t), \;\; t>0,\\
    u(x,t)=-\kappa(x,t),   &\qquad x\in\partial\Omega(t), \;\; t>0,\\
    V_n(x,t)=\partial_{\bfn}u(x,t), &\qquad x\in\partial\Omega(t), \;\; t>0,\\
   \Omega(0)=\Omega_0, &
\end{array}
\right.
\end{equation}
  where for each $t>0$, $\Omega(t)$ is an unknown domain in ${\mathbf{R}}^n$, $\Delta$ is the Laplacian in $n$ variables, $u=u(x,t)$ is an
  unknown function defined for $x\in\overline{\Omega(t)}$ and $t\geqslant 0$, $\kappa(\cdot,t)$ is the mean curvature of the boundary
  $\partial\Omega(t)$ of $\Omega(t)$, $V_n$ is the normal velocity of the free boundary $\partial\Omega(t)$, $\bfn$ denotes the outward-pointing
  normal field of $\partial\Omega(t)$, and $\Omega_0$ is a given initial domain. As before we take the convention that for a convex domain the mean
  curvature of its boundary takes nonnegative values. For simplicity we only consider the case that $\Omega_0$ is a simple domain sufficiently close
  to a sphere.

  The above problem has been intensively studied during the past fifty years. It has been proved that this problem is locally well-posed in H\"{o}lder
  and Sobolev spaces, cf. \cite{Chen},  \cite{ChenHY}, \cite{EscSim2}, \cite{EscSim3} and references therein. Moreover, it has also been proved that
  if $\Omega_0$ is a small perturbation of a sphere, then the solution of the above problem exists for all $t\geqslant 0$ and $\Omega(t)$ converges to
  a sphere with same volume as $\Omega_0$ as $t\to\infty$, cf. \cite{Chen}, \cite{EscSim4}, \cite{EscSim5} and references therein. In what follows we
  use Theorem 3.4 to give a more prcise description of the local phase diagram of the above problem near radial stationary solutions.

  Let $m$ be a positive integer $\geqslant 2$ and $0<\mu<1$. Let $\mathfrak{M}:=\dot{\mathfrak{D}}^{m+\mu}({\mathbf{R}}^n)$ and $\mathfrak{M}_0:=
  \dot{\mathfrak{D}}^{m+3+\mu}({\mathbf{R}}^n)$. We know that $\mathfrak{M}$ and $\mathfrak{M}_0$ are Banach manifolds built on the Banach spaces
  $\dot{C}^{m+\mu}({\mathbf{S}}^{n-1})$ and $\dot{C}^{m+3+\mu}({\mathbf{S}}^{n-1})$, respectively, and $\dot{C}^{m+3+\mu}({\mathbf{S}}^{n-1})$ is dense in
  $\dot{C}^{m+\mu}({\mathbf{S}}^{n-1})$. The problem (4.5) can be reduced into a differential equation in the Banach manifold $\mathfrak{M}$. Indeed,
  for any $Q\in\mathfrak{M}_0\subseteq\mathfrak{M}$ we use the standard local chart of $\mathfrak{M}$ at $Q$ to identify the tangent space
  $\mathcal{T}_{Q}(\mathfrak{M})$ with the Banach space $\dot{C}^{m+\mu}(\partial Q)$. It follows that if $I\subseteq{\mathbf{R}}$ is an open interval
  and $\Omega:I\to\mathfrak{M}_0$ is a $C^1$ curve, then
$$
   \Omega'(t)=V_n(\cdot,t) \quad \mbox{for}\;\; t\in I.
$$
  Given $\Omega\in\mathfrak{M}_0$ we denote by $u_{\Omega}$ the unique solution of the following Dirichlet problem:
$$
\left\{
\begin{array}{ll}
   \Delta u_{\Omega}=0 &\quad\;\; \mbox{in} \;\; \Omega,\\
   u_{\Omega}=-\kappa_{\partial\Omega}   &\quad \mbox{on} \;\; \partial\Omega,
\end{array}
\right.
$$
  where $\kappa_{\partial\Omega}$ denotes the mean curvature of $\partial\Omega$. It is known that $u_{\Omega}\in\dot{C}^{m+\mu+1}(\overline{\Omega})$.
  Now introduce a vector field $\mathscr{F}$ in $\mathfrak{M}$ with domain $\mathfrak{M}_0$ as follows: For any $\Omega\in\mathfrak{M}_0$ we define
$$
   \mathscr{F}(\Omega)=\partial_{\bfn}u_{\Omega}|_{\partial\Omega}\in\dot{C}^{m+\mu}(\partial\Omega)=\mathcal{T}_{\Omega}(\mathfrak{M}).
$$
  It follows that the problem (4.5) reduces into the following initial value problem of a differential equation in the Banach manifold $\mathfrak{M}$:
\begin{equation}
\left\{
\begin{array}{ll}
   \Omega'(t)=\mathscr{F}(\Omega(t)), &\quad  t>0,\\
    \Omega(0)=\Omega_0. &
\end{array}
\right.
\end{equation}

  Given $Q\in\mathfrak{M}_0$, let $(\mathcal{U},\varphi)$ be a regular local chart of $\mathfrak{M}$ at $Q$ with smooth base hypersurface
  $S$ (so that $O:=\varphi(\mathcal{U})$ is an open subset of $X:=\dot{C}^{m+\mu}(S)$), and $\mathcal{U}_0\subseteq\mathcal{U}$ a neighborhood of
  $Q$ in $\mathfrak{M}_0$ such that $(\mathcal{U}_0,\varphi|_{\mathcal{U}_0})$ is a local chart of $\mathfrak{M}_0$ at $Q$ (so that $O_0:=
  \varphi(\mathcal{U}_0)$ is an open subset of $X_0:=\dot{C}^{m+3+\mu}(S)$). Let $F:O_0\to X$ be representation of $\mathscr{F}$ in this local
  chart, i.e.,
$$
  F(u)=\varphi'(\varphi^{-1}(u))\mathscr{F}(\varphi^{-1}(u)), \quad \forall u\in O_0.
$$
  Then representation of the problem (4.6) in the neighborhood $\mathcal{U}$ of $Q$ in $\mathfrak{M}$ is the following initial value problem in the
  Banach space $X=\dot{C}^{m+\mu}(S)$:
$$
\left\{
\begin{array}{ll}
   u'(t)=F(u(t)), &\quad  t>0,\\
   u(0)=u_0, &
\end{array}
\right.
$$
  where $u(t)=\varphi(\Omega(t))$ and $u_0=\varphi(\Omega_0)$.

  From references \cite{EscSim2}--\cite{EscSim5} we know that $F\in C^{\infty}(O_0,X)$, and for any $u\in O_0$, $F'(u)$ is a sectorial
  operator in $X$ with domain $X_0$ and the graph norm of ${\rm Dom}F'(u)=X_0$ is equivalent to the norm of $X_0$. Hence the equation $\Omega'=
  \mathscr{F}(\Omega)$ is a parabolic differential equation in the Banach manifold $\mathfrak{M}$. Moreover, it is easy to check that the vector field
  $\mathscr{F}$ is invariant under the translation group action $(G_{tl},p)$, and quasi-invariant under the dilation group action $(G_{dl},q)$ with
  quasi-invariant factor $\theta(\lambda)=\lambda^{-3}$, $\lambda>0$. The discussion in Section 2 shows that these Lie group actions satisfy the
  conditions $(L1)$--$(L8)$. Concerning the other conditions in Theorem 3.4, we have the following preliminary result:
\medskip

  {\bf Lemma 4.1}\ \ {\em Let $Q=B(0,1)$ and $S={\mathbf{S}}^{n-1}$ $($so that $X=\dot{C}^{m+\mu}({\mathbf{S}}^{n-1})$ and $X_0=
  \dot{C}^{m+3+\mu}({\mathbf{S}}^{n-1}))$. We have the following assertions:

  $(1)$\ \ $F'(0)\eta=\displaystyle\frac{1}{n\!-\!1}\Delta_{\omega}D\eta+D\eta$, where $D$ is the Dirichlet-Neumann operator
  on the sphere ${\mathbf{S}}^{n-1}$.

  $(2)$\ \ $\sigma(F'(0))=\{0\}\cup\{\mu_k:k=2,3,\cdots\}$, where $\mu_k=\displaystyle-\frac{k(k\!-\!1)(k\!+\!n\!-\!1)}{n\!-\!1}$, $k=2,3,\cdots$.

  $(3)$\ \ ${\rm Ker}\,F'(0)={\rm span}\{{\mathbf{1}},Y_{11}(\omega),Y_{12}(\omega),\cdots,Y_{1n}(\omega)\}$, where $Y_{11}(\omega)$,
  $Y_{12}(\omega)$, $\cdots$, $Y_{1n}(\omega)$ are the basis of the linear space of first-order sphere harmonics. Hence $\dim{\rm Ker}\,F'(0)=n+1$.

  $(4)$\ \ ${\rm Range}\,F'(0)$ is closed, and $X={\rm Ker}\,F'(0)\oplus{\rm Range}\,F'(0)$.}
\medskip

  {\em Proof.}\ \ The proof is similar to those of Lemma 6.3, Corollary 6.4 and Lemma 6.7 in the next section (cf. also \cite{EscSim4},
  \cite{EscSim5} for similar discussion). To save spaces we omit it here. $\quad\Box$
\medskip

  The above lemma shows that Theorem 3.4 applies to the equation (4.6). Hence, by applying Theorem 3.4 we get the following result:
\medskip

  {\bf Theorem 4.2}\ \ {\em Let $\mathcal{M}_c$ be the $(n\!+\!1)$-dimensional submanifold of $\mathfrak{M}_0$ consisting of all spheres in
  ${\mathbf{R}}^n$. We have the following assertions:

  $(1)$\ \ There is a neighborhood $\mathcal{O}$ of $\mathcal{M}_c$ in $\mathfrak{M}_0$ such that for any $\Omega_0\in\mathcal{O}$, the initial value
  problem $(4.6)$ has a unique solution $\Omega\in C([0,\infty),\mathfrak{M}_0)\cap C^1((0,\infty),\mathfrak{M}_0)$.

  $(2)$\ \ There exists a submanifold $\mathcal{M}_s$ of $\mathfrak{M}_0$ of codimension $n\!+\!1$ passing $\Omega_s=B(0,1)$ such that for any
  $\Omega_0\in\mathcal{M}_s$, the solution of the problem $(4.6)$ satisfies $\displaystyle\lim_{t\to\infty}\Omega(t)=\Omega_s$ and, conversely, if
  the solution of $(4.6)$ satisfies this property then $\Omega_0\in\mathcal{M}_s$.

  $(3)$\ \ For any $\Omega_0\in\mathcal{O}$ there exist unique $x_0\in {\mathbf{R}}^n$, $R>0$ and $Q_0\in\mathcal{M}_s$ such that
  $\Omega_0=x_0+RQ_0$ and, for the solution $\Omega=\Omega(t)$ of $(4.6)$, we have}
$$
   \lim_{t\to\infty}\Omega(t)=B(x_0,R).
$$
$\Box$
\medskip

  {\em Remark}.\ \ From the conservation $|\Omega(t)|=|\Omega_0|$ ($\forall t\geqslant 0$) we see that $R=\sqrt[n]{|\Omega_0|/\omega_n}$, where
  $\omega_n$ denotes the volume of the $n$-dimensional unit sphere. However, how to compute $x_0$ is unknown.

\section{An obstacle problem}
\setcounter{equation}{0}

\hskip 2em
  In this section we study the obstacle problem (1.6). The purpose is through a such study to prove that for the solution $(\sigma,\pi_0)$ of the
  boundary value problem
\begin{equation}
\left\{
\begin{array}{rll}
   \Delta\sigma=&f(\sigma) &\quad\;\; \mbox{in} \;\; \Omega,\\
   -\Delta\pi_0=&g(\sigma)   &\quad\;\; \mbox{in} \;\; \Omega,\\
   \sigma=&1   &\quad\;\; \mbox{on} \;\; \partial\Omega,\\
   \pi_0=&0  &\quad\;\; \mbox{on} \;\; \partial\Omega,
\end{array}
\right.
\end{equation}
  where $f$, $g$ are the discontinuous functions given in (1.5) (with $\lambda=1$), the mapping $\Omega\mapsto\partial_{\bfn}\pi_0|_{\partial\Omega}$
  from a neighborhood of a sphere in $\mathfrak{M}_0:=\dot{\mathfrak{D}}^{m+3+\mu}({\mathbf{R}}^3)\subseteq\mathfrak{M}:=
  \dot{\mathfrak{D}}^{m+\mu}({\mathbf{R}}^3)$ to $\mathcal{T}_{\mathfrak{M}_0}(\mathfrak{M})$ is smooth (i.e., representation of this mapping in some
  regular local chart of $\mathfrak{M}$ at every sphere is smooth), where $\partial_{\bfn}$ denotes the derivative in the outward normal direction
  $\bfn$ of $\partial\Omega$. This result will be used in the next section to study the free boundary problem (1.3). We note that since the functions
  $f$, $g$ are discontinuous, such a result apparently looks unbelievable.

  We point out that although here we only consider the three dimension case, a similar discussion also works for general dimension $n\geqslant 2$
  case; in order to do so the discussion in \cite{Cui2} must be first extended, which is not hard.

  Let $m$ be a positive integer, $m\geqslant 2$, and $0<\mu<1$. Let $\Omega\in\mathfrak{D}^{m+\mu}({\mathbf{R}}^3)$ be given. It is easy to see that
  the obstacle problem (1.6) is equivalent to the following boundary value problem with discontinuous function $f(\sigma)=\sigma
  H(\sigma-\hat{\sigma})$:
\begin{equation}
\left\{
\begin{array}{rll}
   \Delta\sigma=&f(\sigma) &\quad\;\; \mbox{in} \;\; \Omega,\\
   \sigma=&1   &\quad\;\; \mbox{on} \;\; \partial\Omega.
\end{array}
\right.
\end{equation}
  Since $0<\hat{\sigma}<1$, it is clear that $\underline{\sigma}=\hat{\sigma}$ and $\overline{\sigma}=1$ are lower and upper solutions of this problem,
  respectively, so that by the upper and lower solution method\footnotemark[3], the above problem has a solution $\sigma$ satisfying the following
  conditions:
\footnotetext[3]{Note that discontinuity of the function $f$ does not produce an essential obstacle for validity of the upper and lower solution method,
  because this difficulty can be easily overcome through a mollification process.}
$$
  \hat{\sigma}\leqslant\sigma(x)\leqslant 1 \;\;  \mbox{for}\;\; x\in\overline{\Omega}, \quad \mbox{and} \quad
  \sigma\in W^{2,p}(\Omega) \;\;  \mbox{for any}\;\; 1\leqslant p<\infty.
$$
  Moreover, applying the regularity theory for elliptic boundary value problems we easily see that $\sigma$ has also the following properties: Letting
  $\Omega_{\rm liv}$ and $\Omega_{\rm nec}$ be as defined in Section 1, i.e., $\Omega_{\rm liv}=\{x\in\Omega: \sigma(x)>\hat{\sigma}\}$ and
  $\Omega_{\rm nec}={\rm int}\{x\in\Omega: \sigma(x)=\hat{\sigma}\}$, then
\begin{equation}
  \sigma\in C^{\infty}(\Omega_{\rm liv}\cup\Omega_{\rm nec}) \quad \mbox{and} \quad
  \sigma\in C^{m+\mu}(\Omega_{\rm liv}\cup\partial\Omega).
\end{equation}
  The last property is implied by the condition that $\partial\Omega$ is of $C^{m+\mu}$-class, the fact that $f(\sigma)$ is smooth for
  $\sigma\in (\hat{\sigma},\infty)$ and the condition $\sigma=1$ on $\partial\Omega$. Note that since $f$ is a nondecreasing function, the solution
  of the problem (5.2) is unique. In what follows we prove that if $\Omega$ is a small perturbation of a sphere, the interface $\Gamma=\partial
  \Omega_{\rm liv}\cap\Omega_{\rm nec}$ is smooth and the mapping $\partial\Omega\to\Gamma$ is smooth.

  We begin by rewriting the problem (5.2) into an equivalent problem which is easier to treat. We first recall (cf. \cite{Cui2}) that if
  $\Omega=B(0,R)$ for some $R>0$, then the unique solution of the problem (5.2) is a radial function, given by $\sigma(x)=U(|x|,R)$, where
  $U(r,R)=R\sinh r/r\sinh R$ for $R<R^*$ and
$$
  U(r,R)=\left\{
\begin{array}{ll}
    \displaystyle\hat{\sigma}[\sinh(r- K)+ K\cosh(r- K)]/r &\quad\;\;
    \mbox{for} \;\;  K\leqslant r\leqslant R\\ [0.2cm]
   \;\;\hat{\sigma}   &\quad\;\; \mbox{for} \;\; r<K
\end{array}
\right.
$$
  for $R>R^*$, where $R^*$ is the unique positive number solving the equation $\sinh\! R^*/R^*=1/\hat{\sigma}$, and for $R>R^*$, $K=K(R)$ is the
  unique solution of the following equation in the interval $(0,R)$:
$$
  \sinh(R-K)+K\cosh(R-K)=\frac{R}{\hat{\sigma}}.
$$
  Given $R>R^*$ and $\rho,\eta\in C^{2+\mu}({\mathbf{S}}^2)$ ($0<\mu<1$) with $\|\rho\|_{C^{2+\mu}({\mathbf{S}}^2)}$ and
  $\|\eta\|_{C^{2+\mu}({\mathbf{S}}^2)}$ sufficiently small, we denote
$$
  \Omega_{\rho}=\{x\in{\mathbf{R}}^3: r<R[1+\rho(\omega)]\}, \quad D_{\rho,\eta}=\{x\in{\mathbf{R}}^3:  K[1+\eta(\omega)]<r<R[1+\rho(\omega)]\},
$$
$$
   S_{\rho}=\{x\in{\mathbf{R}}^3: r=R[1+\rho(\omega)]\}, \quad \mbox{and} \quad \Gamma_{\eta}=\{x\in{\mathbf{R}}^3: r= K[1+\eta(\omega)]\}.
$$
  Moreover, let $\sigma_{\rho}$ be the solution of the problem (5.2) for $\Omega=\Omega_{\rho}$. The equivalent problem mentioned above is as
  follows:
\begin{equation}
\left\{
\begin{array}{rll}
   \Delta\sigma=&\sigma &  \quad  \mbox{in}\;\; D_{\rho,\eta},\\
   \sigma=&1 &  \quad  \mbox{on}\;\; S_{\rho},\\
   \sigma=&\hat{\sigma} &  \quad  \mbox{on}\;\; \Gamma_{\eta},\\
    \partial_{r}\sigma=&0 & \quad  \mbox{on}\;\; \Gamma_{\eta},
\end{array}
\right.
\end{equation}
  where $\partial_{r}$ denotes the derivative in radial direction. Later on we shall also use the following abbreviations:
$$
  D=D_{0,0}=B(0,R)\backslash\overline{B(0,K)}, \quad S_0=\partial B(0,R), \quad \mbox{and} \quad  \Gamma_0=\partial B(0,K).
$$

  {\bf Lemma 5.1}\ \ {\em For $\rho,\eta\in C^{2+\mu}({\mathbf{S}}^2)$ $(0<\mu<1)$ with $\|\rho\|_{C^{2+\mu}({\mathbf{S}}^2)}$ and
  $\|\eta\|_{C^{2+\mu}({\mathbf{S}}^2)}$ sufficiently small, if $\Omega=\Omega_{\rho}$ and $\Gamma=\Gamma_{\eta}$ then the problem $(5.2)$ is
  equivalent to the problem $(5.4)$.}
\medskip

  {\em Proof}.\ \ Let $\nu$ be the outward unit normal field of the inner part boundary $\Gamma_{\eta}$ of $D_{\rho,\eta}$. It is easy to see that
  under the assumption $\rho,\eta\in C^{2+\mu}({\mathbf{S}}^2)$, $\Omega=\Omega_{\rho}$ and $\Gamma=\Gamma_{\eta}$, the problem (5.2) is
  equivalent to the following problem:
\begin{equation}
\left\{
\begin{array}{rll}
   \Delta\sigma=&\sigma &  \quad  \mbox{in}\;\; D_{\rho,\eta},\\
   \sigma=&1 &  \quad  \mbox{on}\;\; S_{\rho},\\
   \sigma=&\hat{\sigma} &  \quad  \mbox{on}\;\; \Gamma_{\eta},\\
    \partial_{\nu}\sigma=&0 & \quad  \mbox{on}\;\; \Gamma_{\eta}.
\end{array}
\right.
\end{equation}
  That is, if $\sigma$ is a solution of (5.2) for $\Omega=\Omega_{\rho}$ and its free boundary $\Gamma$ has the form $\Gamma=\Gamma_{\eta}$ with
  $\eta\in C^{2+\mu}({\mathbf{S}}^{n-1})$, then the restriction of $\sigma$ to $\overline{D}_{\rho,\eta}$ is a solution of (5.5); conversely, if
  $\sigma$ is a solution of (5.5) then by extending it into the whole domain $\Omega=\Omega_{\rho}$ such that it identically takes the value
  $\hat{\sigma}$ in $\Omega\backslash\overline{D}_{\rho,\eta}$, then after such extension $\sigma$ is a solution of (5.2).
  The condition $\sigma=\hat{\sigma}$ on $\Gamma_{\eta}$ implies $\nabla_{\omega}\sigma(\eta(\omega),\omega)=-\partial_r\sigma(\eta(\omega),\omega)
  \nabla_{\omega}\eta(\omega)$ for $\omega\in{\mathbf{S}}^{n-1}$. Since $\nu=[\omega-(1/r)\nabla_{\omega}\eta(\omega)]/\sqrt{1+(1/r^2)
  |\nabla_{\omega}\eta(\omega)|^2}|_{r=\eta(\omega)}$ and $\nabla\sigma=(\partial_{r}\sigma)\omega+(1/r)\nabla_{\omega}\sigma$, so that
$$
  \partial_{\nu}\sigma|_{r=\eta(\omega)}=\frac{\partial_{r}\sigma-(1/r^2)\nabla_{\omega}\eta(\omega)\cdot\nabla_{\omega}\sigma}{\sqrt{1+(1/r^2)
  |\nabla_{\omega}\eta(\omega)|^2}}\Big|_{r=\eta(\omega)},
$$
  we see that the condition $\sigma=\hat{\sigma}$ on $\Gamma_{\eta}$ implies $\partial_{\nu}\sigma|_{r=\eta(\omega)}=\sqrt{1+(1/r^2)
  |\nabla_{\omega}\eta(\omega)|^2}\partial_{r}\sigma|_{r=\eta(\omega)}$. Hence the problems (5.4) and (5.5) are equivalent. This proves the desired
  assertion. $\quad\Box$
\medskip

  The main result of this section is the following theorem:
\medskip

  {\bf Theorem 5.2}\ \ {\em For any $R>R^*$, integer $m\geqslant 2$ and $0<\mu<1$ there exists corresponding $\delta>0$ such that for any $\rho\in
  C^{m+\mu}({\mathbf{S}}^2)$ with $\|\rho\|_{C^{m+\mu}({\mathbf{S}}^2)}<\delta$, the problem $(5.4)$ has a unique solution $(\sigma,\eta)$ with
  $\eta\in C^{\infty}({\mathbf{S}}^2)$ and $\sigma\in C^{m+\mu}(\overline{D}_{\rho,\eta})\cap C^{\infty}(\overline{D}_{\rho,\eta}\backslash
  S_{\rho})$, and the mapping $\rho\mapsto\eta$ from the open set $\|\rho\|_{C^{m+\mu}({\mathbf{S}}^2)}<\delta$ in $C^{m+\mu}({\mathbf{S}}^2)$
  to $C^{\infty}({\mathbf{S}}^2)$ is smooth.}
\medskip

  {\em Proof}.\ \ We know that $C^{\infty}({\mathbf{S}}^2)$ with the family of seminorms $\{\|\cdot\|_{C^({\mathbf{S}}^2)}\}\cup
  \{\|\cdot\|_{C^{m+\mu}({\mathbf{S}}^2)}\}_{m=1}^{\infty}$ is a tame Frech\'{e}t space. We also regard the Banach space
  $C^{m+\mu}({\mathbf{S}}^2)$ as a tame Frech\'{e}t space. Let $R>R^*$ be given and set $K=K(R)$. For sufficiently small $\delta,\delta'>0$
  we denote
$$
  O_{\delta}=\{\rho\in C^{m+\mu}({\mathbf{S}}^2):\|\rho\|_{C^{m+\mu}({\mathbf{S}}^2)}<\delta\}, \quad
  O_{\delta'}'=\{\eta\in C^{\infty}({\mathbf{S}}^2):\|\eta\|_{C^{2+\mu}({\mathbf{S}}^2)}<\delta'\};
$$
  they are open subsets of $C^{m+\mu}({\mathbf{S}}^2)$ and $C^{\infty}({\mathbf{S}}^2)$, respectively. We define a map $A: O_{\delta}\times
  O_{\delta'}'\subseteq C^{m+\mu}({\mathbf{S}}^2)\times C^{\infty}({\mathbf{S}}^2)\to C^{\infty}({\mathbf{S}}^2)$ as follows: Given $\rho\in
  O_{\delta}$ and $\eta\in O_{\delta'}'$, let $\sigma=\sigma(r,\omega;\rho,\eta)$ be the unique solution of the equations $(5.4)_1$, $(5.4)_2$ and
  $(5.4)_4$, and define
$$
  A(\rho,\eta)=[\omega\mapsto\sigma(K[1+\eta(\omega)],\omega;\rho,\eta)-\hat{\sigma},\omega\in{\mathbf{S}}^2].
$$
  Clearly $A(0,0)=0$. Moreover, $A$ is a smooth tame map. To prove this assertion we choose a function $\phi\in C^{\infty}[ K,R]$ such that it
  satisfies the following conditions:
$$
  0\leqslant\phi\leqslant 1; \quad \phi(R)=\phi( K)=1; \quad
  \phi(t)=0 \;\; \mbox{for}\;\,\frac{3}{4} K+\frac{1}{4}R\leqslant t\leqslant\frac{1}{4} K+\frac{3}{4}R;
$$
$$
  \phi'(t)\leqslant 0 \;\; \mbox{for}\;\, K\leqslant t\leqslant\frac{3}{4} K+\frac{1}{4}R; \quad
  \phi'(t)\geqslant 0 \;\; \mbox{for}\;\,\frac{1}{4} K+\frac{3}{4}R\leqslant t\leqslant R.
$$
  Let $M_0=\displaystyle\max_{ K\leqslant t\leqslant R}|\phi'(t)|$ and assume $\delta,\delta'$ are small enough such that $\delta<
  (1+M_0R)^{-1}$, $\delta'<(1+M_0 K)^{-1}$ and $\max\{\delta,\delta'\}<\displaystyle\frac{1}{3}\frac{R- K}{R+ K}$. Consider
  the variable transformation $y=\Psi_{\rho,\eta}(x)$ from $\overline{D}_{\rho,\eta}$ to $\overline{D}$, where for $x\in\overline{D}_{\rho,\eta}$,
\begin{equation}
  \Psi_{\rho,\eta}(x)=\left\{
\begin{array}{ll}
   \displaystyle x-R\rho(\omega)\phi\Big(\frac{r}{1+\rho(\omega)}\Big)\omega &\quad\;\;\mbox{if} \;\; r\geqslant\frac{1}{2}( K+R),\\ [0.3cm]
   \displaystyle x- K\eta(\omega)\phi\Big(\frac{r}{1+\eta(\omega)}\Big)\omega &\quad\;\;\mbox{if} \;\; r<\frac{1}{2}( K+R).
\end{array}
\right.
\end{equation}
  It is easy to see that $\Psi_{\rho,\eta}$ is a $C^{m+\mu}$ diffeomorphism from $\overline{D}_{\rho,\eta}$ onto $\overline{D}$. Moreover, denoting
$$
  E_{\eta}=\{x\in{\mathbf{R}}^3:  K[1+\eta(\omega)]<r<\frac{1}{2}( K+R)\}, \quad
  E=\{x\in{\mathbf{R}}^3:  K<r<\frac{1}{2}( K+R)\},
$$
  we see that the restriction of $\Psi_{\rho,\eta}$ on $\overline{E}_{\eta}$ is a $C^{\infty}$-diffeomorphism from $\overline{E}_{\eta}$ onto
  $\overline{E}$, due to the facts that $\eta\in C^{\infty}({\mathbf{S}}^2)$ and that this restriction is independent of $\rho$. Because of
  the latter property, we re-denote the restriction of $\Psi_{\rho,\eta}$ on $\overline{E}_{\eta}$ as $\Psi_{\eta}$, and denote by $\psi_{\eta}$ the
  restriction of $\Psi_{\eta}$ to $\Gamma_{\eta}$, which is clearly a $C^{\infty}$-diffeomorphism from $\Gamma_{\eta}$ onto $\Gamma_0$. Now define
  operators $\mathscr{A}(\rho,\eta):C^{m+\mu}(\overline{D})\cap C^{\infty}(\overline{E})\to C^{m-2+\mu}(\overline{D})\cap C^{\infty}(\overline{E})$ and
  $\mathscr{N}(\eta):C^{\infty}(\overline{E})\to C^{\infty}(\Gamma_0)$ respectively as follows:
$$
  \mathscr{A}(\rho,\eta)u=[\Delta(u\circ\Psi_{\rho,\eta})]\circ\Psi_{\rho,\eta}^{-1} \quad\;\;
   \mbox{for}\;\,u\in C^{m+\mu}(\overline{D})\cap C^{\infty}(\overline{E}),
$$
$$
  \mathscr{N}(\eta)u=[\partial_{\nu}(u\circ\Psi_{\eta})|_{\Gamma_{\eta}}]\circ\psi_{\eta}^{-1} \quad\;\;
  \mbox{for}\;\,u\in C^{\infty}(\overline{E}).
$$
  Let $u=\sigma\circ\Psi_{\rho,\eta}^{-1}$. After the variable transformation $x\mapsto\Psi_{\rho,\eta}(x)$, the problem $(5.4)_1$, $(5.4)_2$ and
  $(5.4)_4$ transforms into the following problem:
\begin{equation}
\left\{
\begin{array}{rll}
   \mathscr{A}(\rho,\eta)u=&u &  \quad  \mbox{in}\;\; D,\\
   u=&1 &  \quad  \mbox{on}\;\; S_0,\\
   \mathscr{N}(\eta)u=&0 & \quad  \mbox{on}\;\; \Gamma_0.
\end{array}
\right.
\end{equation}
  Note that for any integer $k\geqslant 2$, if we denote $O_{\delta'}^k=\{\eta\in C^{k+\mu}({\mathbf{S}}^2):\|\eta\|_{C^{2+\mu}({\mathbf{S}}^2)}<
  \delta'\}$, then
\begin{equation}
\left\{
\begin{array}{l}
   [(\rho,\eta)\mapsto\mathscr{A}(\rho,\eta)]\in C^{\infty}(O_{\delta}\times O_{\delta'}^k,L(C^{m+\mu}(\overline{D})\cap C^{k+\mu}(\overline{E}),
   C^{m-2+\mu}(\overline{D})\cap C^{k-2+\mu}(\overline{E}))),\\ [0.1cm]
   [\eta\mapsto\mathscr{N}(\eta)]\in C^{\infty}(O_{\delta'}^k,L(C^{k+\mu}(\overline{E}),C^{k-1+\mu}(\Gamma))).
\end{array}
\right.
\end{equation}
  From these properties and the standard theory for elliptic boundary value problems we easily see that the solution map $(\rho,\eta)\mapsto u$
  of the problem (5.7) is a smooth tame map (cf. Theorem 3.3.1 in Part II of \cite{Ham1}): For any integer $k\geqslant 2$,
\begin{equation}
   [(\rho,\eta)\mapsto u]\in C^{\infty}(O_{\delta}\times O_{\delta'}^k,C^{m+\mu}(\overline{D})\cap C^{k+\mu}(\overline{E})),
\end{equation}
  and for $k\geqslant m$ (cf. Lemma 3.3.2 in Part II of \cite{Ham1}),
$$
  \|u\|_{C^{m+\mu}(\overline{D})}+\|u\|_{C^{k+\mu}(\overline{E})}\leqslant C_k(1+\|\eta\|_{C^{k+\mu}({\mathbf{S}}^2)})
$$
  for some constant $C_k>0$, and similarly for any order Frech\'{e}t derivatives of $u$ in $(\rho,\eta)$. Since
$$
  A(\rho,\eta)=\Big[\omega\mapsto u(K\omega)-\hat{\sigma},
  \omega\in{\mathbf{S}}^2\Big],
$$
  the desired assertion immediately follows.

  Now, a simple computation shows that for $\rho\in O_{\delta}$, $\eta\in O_{\delta'}'$ and $\zeta\in C^{\infty}({\mathbf{S}}^2)$,
  $\partial_{\eta}A(\rho,\eta)\zeta=[\omega\mapsto u(K[1+\eta(\omega)]\omega;\rho,\eta),\omega\in{\mathbf{S}}^2]$, where $u=u(\cdot;\rho,\eta)$ is the
  solution of the following problem:
$$
  \Delta u=u \;\;  \mbox{in} \;\; D_{\rho,\eta}, \quad  u=0 \;\;  \mbox{on}\;\; S_{\rho},  \quad
   \partial_ru=-K\partial_r^2\sigma|_{\Gamma_{\eta}}\zeta\;\;\;  \mbox{on} \;\; \Gamma_{\eta}.
$$
  If $\rho=0$ and $\eta=0$ then $\sigma=U(r,R)$, so that $\partial_r^2\sigma(\cdot;0,0)|_{\Gamma_{\eta}}=\partial_r^2U(K,R)=U(K,R)=\hat{\sigma}$.
  We now choose $\delta,\delta'>0$ sufficiently small such that for all $\rho\in O_{\delta}$ and $\eta\in O_{\delta'}'$ there hold $(1/2)\hat{\sigma}
  \leqslant\partial_r^2\sigma(\cdot;\rho,\eta)|_{\Gamma_{\eta}}\leqslant 2\hat{\sigma}$. Then for any $\rho\in O_{\delta}$ and $\eta\in O_{\delta'}'$
  the operator $\partial_{\eta}A(\rho,\eta)$ is invertible, with
$$
  [\partial_{\eta}A(\rho,\eta)]^{-1}\xi=-\frac{\partial_rv(\cdot;\rho,\eta)}{K\partial_r^2\sigma(\cdot;\rho,\eta)}\Big|_{\Gamma_{\eta}},
  \quad \forall\xi\in C^{\infty}({\mathbf{S}}^{n-1}),
$$
  where $v=v(\cdot;\rho,\eta)$ is the solution of the following problem:
$$
  \Delta v=v \;\;  \mbox{in} \;\; D_{\rho,\eta}, \quad  v=0 \;\;  \mbox{on}\;\; S_{\rho},  \quad
  v=\xi\;\;  \mbox{on} \;\; \Gamma_{\eta}.
$$
  Similarly as before we can prove the map $(\rho,\eta,\xi)\mapsto [\partial_{\eta}A(\rho,\eta)]^{-1}\xi$ is tame. Hence, by the Nash-Moser implicit
  function theorem (cf. Theorem 3.3.1 in Part III of \cite{Ham1}) we conclude that for any $\rho\in O_{\delta}$ there exists $\eta\in O_{\delta'}'$
  such that it is the unique solution of the equation $A(\rho,\eta)=0$ in $O_{\delta'}'$, and the map $\rho\mapsto\eta$ from $O_{\delta}\subseteq
  C^{m+\mu}({\mathbf{S}}^2)$ to $O_{\delta'}'\subseteq C^{\infty}({\mathbf{S}}^2)$ is a smooth tame map. This shows that for any $\rho\in
  C^{m+\mu}({\mathbf{S}}^2)$ with $\|\rho\|_{C^{m+\mu}({\mathbf{S}}^2)}<\delta$, the free boundary $\Gamma_{\eta}$ is smooth and the mapping
  $\rho\mapsto\eta$ is also smooth. Having proved smoothness of the free boundary $\Gamma_{\eta}$, the assertion $\sigma\in
  C^{m+\mu}(\overline{D}_{\rho,\eta})\cap C^{\infty}(\overline{D}_{\rho,\eta}\backslash S_{\rho})$ follows immediately. $\quad\Box$
\medskip

  {\em Remark}.\ \ The reader might argue why we don't use the usual implicit function theorem, i.e., the implicit function theorem in Banach spaces,
  to get a simpler proof. The reason is that all our such efforts failed, despite that they are only for the purpose to get a weaker result that
  $\eta$ and $\rho\mapsto\eta$ are finite-order smooth.
\medskip

  Next we consider the problem (5.1). Let $\delta$ and $O_{\delta}$ be as in the proof of Theorem 5.2. Given $\rho\in O_{\delta}$, we first solve
  the problem $(5.1)_1$, $(5.1)_3$ and next substitute the solution $\sigma$ into $(5.1)_2$ and take $(5.1)_4$ into account. Then we obtain the
  following elliptic boundary value problem:
\begin{equation}
\left\{
\begin{array}{rll}
   -\Delta\pi_0=&g(\sigma) &\quad\;\; \mbox{in} \;\; \Omega,\\
   \pi_0=&0   &\quad\;\; \mbox{on} \;\; \partial\Omega.
\end{array}
\right.
\end{equation}
  By applying the standard theory for elliptic boundary value problems and using the properties of $\sigma$ in (5.3), we see the above problem
  has a unique solution satisfying the following properties:
\begin{equation}
  \pi_0\in W^{2,p}(\Omega) \;(\forall p\in [1,\infty)), \quad
  \pi_0\in C^{\infty}(\Omega_{\rm liv}\cup\Omega_{\rm nec}) \quad \mbox{and} \quad
  \pi_0\in C^{m+\mu}(\Omega_{\rm liv}\cup\partial\Omega).
\end{equation}
  It follows that $\partial_{\bfn}\pi|_{\partial\Omega}\in C^{m-1+\mu}(\partial\Omega)$, where $\bfn$ denotes the unit outward normal field
  of $\partial\Omega$. In this way we obtain a map $F_0:O_{\delta}\subseteq C^{m+\mu}({\mathbf{S}}^2)\to C^{m-1+\mu}({\mathbf{S}}^2)$ defined as
  follows: For any $\rho\in O_{\delta}$,
$$
  F_0(\rho)=[\omega\mapsto\partial_{\bfn}\pi_0(R[1+\rho(\omega)],\omega),\omega\in{\mathbf{S}}^2].
$$
  Our next goal of this section is to prove the following result:
\medskip

  {\bf Lemma 5.3}\ \ {\em $F_0\in C^{\infty}(O_{\delta},C^{m-1+\mu}({\mathbf{S}}^2))$.}
\medskip

  {\em Proof}.\ \ From the proof of Theorem 5.2 we easily see that not only the map $\rho\mapsto\eta$ is tame, but also the map
  $\rho\mapsto\sigma|_{\overline{E}}$ from  $O_{\delta}\subseteq C^{m+\mu}({\mathbf{S}}^2)$ to $C^{\infty}(\overline{E})$ is tame.
  Let $D_{\rho,\eta}$, $S_{\rho}$ and $\Gamma_{\eta}$ be as before and set
$$
  B_{\eta}=\{x\in{\mathbf{R}}^3: r< K[1+\eta(\omega)]\}, \qquad B_0=B(0,K).
$$
  Then (5.10) can be rewritten as the following equivalent problem:
\begin{equation}
\left\{
\begin{array}{rll}
   -\Delta\pi_0=&a(\sigma-\tilde{\sigma})-b & \;  \mbox{in}\;\; D_{\rho,\eta},\\
   -\Delta\pi_0=&-b & \;  \mbox{in}\;\; B_{\eta},\\
   \pi_0=&0 & \;  \mbox{on}\;\; S_{\rho},\\
   \pi_0,\partial_{\nu}\pi_0\;& \;\mbox{are continuous}\;& \; \mbox{across}\;\, \Gamma_{\eta},
\end{array}
\right.
\end{equation}
  where $\nu$ is as before (note that it is also the inward unit normal field of the boundary $\Gamma_{\eta}$ of $B_{\eta}$). Let $\Psi_{\rho,\eta}$,
  $\Psi_{\eta}$, $\psi_{\eta}$, $\mathscr{A}(\rho,\eta)$, $u$ be as in the proof of Theorem 5.2 and define
$$
  \mathscr{B}(\eta)u=[\partial_{\nu}(u\circ\Psi_{\eta})|_{\Gamma_{\eta}}]\circ\psi_{\eta}^{-1} \quad\;\;
  \mbox{for}\;\,u\in C^{\infty}(\overline{E}).
$$
  Choose another smooth function $\phi_1\in C^{\infty}[0, K]$ such that it satisfies the following conditions:
$$
  0\leqslant\phi_1\leqslant 1; \quad \phi_1'\geqslant 0; \quad \phi_1(t)=0 \;\; \mbox{for}\;\,0\leqslant t\leqslant\frac{1}{2} K; \quad
  \phi_1( K)=1.
$$
  Let $M_1=\displaystyle\max_{0\leqslant t\leqslant K}|\phi_1'(t)|$ and assume $\delta'>0$ is small enough such that in addition to the
  conditions appearing in the proof of Lemma 5.1 we have also $\delta'<(1+M_1 K)^{-1}$. Let $\Psi_{\eta}^1:\overline{B}_{\eta}\to \overline{B}_0$ be
  as follows:
\begin{equation}
  \Psi_{\eta}^1(x)=x- K\eta(\omega)\phi_1\Big(\frac{r}{1+\eta(\omega)}\Big)\omega \quad\;\;\mbox{for} \;\; x\in\overline{B}_{\eta}.
\end{equation}
  Define $\mathscr{A}_1(\eta):C^{\infty}(\overline{B}_0)\to C^{\infty}(\overline{B}_0)$ and $\mathscr{B}_1(\eta):C^{\infty}(\overline{B}_0)\to
  C^{\infty}(\Gamma_0)$ respectively as follows:
$$
  \mathscr{A}_1(\eta)u=[\Delta(u\circ\Psi_{\eta}^1)]\circ(\Psi_{\eta}^1)^{-1} \quad\;\; \mbox{for}\;\,u\in C^{\infty}(\overline{B}_0),
$$
$$
  \mathscr{B}_1(\eta)u=[\partial_{\nu}(u\circ\Psi_{\eta}^1)|_{\Gamma_{\eta}}]\circ(\psi_{\eta}^1)^{-1} \quad\;\;
  \mbox{for}\;\,u\in C^{\infty}(\overline{B}_0),
$$
  where $\psi_{\eta}^1=\Psi_{\eta}^1|_{\Gamma_{\eta}}$. Let $v=\pi_0|_{\overline{D}_{\rho,\eta}}\circ\Psi_{\rho,\eta}^{-1}$ and $v_1=
  \pi_0|_{\overline{B}_{\eta}}\circ(\Psi_{\eta}^1)^{-1}$. After the variable transformation $x\mapsto\Psi_{\rho,\eta}(x)$ (for $x\in
  \overline{D}_{\rho,\eta}$) and $x\mapsto\Psi_{\eta}^1(x)$ (for $x\in B_{\eta}$), the problem (5.12) transforms into the following problem:
\begin{equation}
\left\{
\begin{array}{rll}
   -\mathscr{A}(\rho,\eta)v=&a(u-\tilde{\sigma})-b &  \quad  \mbox{in}\;\; D,\\
   \mathscr{A}_1(\eta)v_1=&b &  \quad  \mbox{in}\;\; B_0,\\
   v=&0 &  \quad  \mbox{on}\;\; S_0,\\
   v=&v_1 &  \quad  \mbox{on}\;\; \Gamma_0,\\
   \mathscr{B}(\eta)v=& \mathscr{B}_1(\eta)v_1 &   \quad  \mbox{on}\;\; \Gamma_0.
\end{array}
\right.
\end{equation}
  Lemma 5.3 is an immediate consequence of Theorem 5.2 and the following preliminary result:
\medskip

  {\bf Lemma 5.4}\ \ {\em Let $\delta,\delta',O_{\delta}$ be as in the proof of Theorem 5.2 and set $O_{\delta'}''=\{\eta\in
  C^{m+\mu}({\mathbf{S}}^2): \|\eta\|_{C^{m+\mu}({\mathbf{S}}^2)}<\delta'\}$. Given $(\rho,\eta,u)\in O_{\delta}'\times O_{\delta'}''
  \times C^{m+\mu}(\overline{D})$, the problem $(5.14)$ has a unique solution $(v,v_1)\in C^{m+\mu}(\overline{D})\times C^{m+\mu}(\overline{B}_0)$,
  and the solution map $(\rho,\eta,u)\mapsto (v,v_1)$ from $O_{\delta}'\times O_{\delta'}''\times C^{m+\mu}(\overline{D})\subseteq
  C^{m+\mu}({\mathbf{S}}^2)\times C^{m+\mu}({\mathbf{S}}^2)\times C^{m+\mu}(\overline{D})$ to $C^{m+\mu}(\overline{D})\times
  C^{m+\mu}(\overline{B}_0)$ is smooth.}
\medskip

  {\em Proof.}\ \ We first note that for $\Omega=B(0,R)$ with $R>R^*$, the unique solution of the problem (5.10) is given by $\pi_0=V(r,R)$, where
$$
  V(r,R)=\left\{
\begin{array}{ll}
    \displaystyle D\Big(\frac{1}{R}\!-\!\frac{1}{r}\Big)-\frac{1}{6}(a\tilde{\sigma}+b)(R^2\!-\!r^2)-a\int_r^R\!\!\!\int_{\xi}^RU(\eta,R)
    \Big(\frac{\eta}{\xi}\Big)^2d\eta d\xi &\quad\;\;
    \mbox{for} \;\;  K\leqslant r\leqslant R,\\ [0.3cm]
   \displaystyle C+\frac{b}{6}r^2  &\quad\;\; \mbox{for} \;\; r< K,
\end{array}
\right.
$$
  where $U(r,R)$ and $K=K(R)$ are as before, and $C,D$ are constants such that the relations $V( K^+,R)=V( K^-,R)$ and $\partial_rV( K^+,R)=
  \partial_rV( K^-,R)$ hold; in particular,
$$
  D=-\frac{1}{3}a\tilde{\sigma}K^3-a\int_{ K}^RU(\eta,R)\eta^2d\eta.
$$
  (cf. Lemmas 3.1 and 3.2 of \cite{Cui2}).

  Given $(\rho,\eta,u,\xi)\in O_{\delta}'\times O_{\delta'}''\times C^{m+\mu}(\overline{D})\times C^{m+\mu}(\Gamma_0)$, we
  consider the following two elliptic boundary value problems:
\begin{equation}
\left\{
\begin{array}{rll}
   -\mathscr{A}(\rho,\eta)v=\,&a(u-\tilde{\sigma})-b &  \quad  \mbox{in}\;\; D,\\
   v=\,&0 &  \quad  \mbox{on}\;\; S_0,\\
   v=\,&\xi &  \quad  \mbox{on}\;\; \Gamma_0,
\end{array}
\right.
  \qquad\qquad
\left\{
\begin{array}{rll}
   \mathscr{A}_1(\eta)v_1=\,&b &  \quad  \mbox{in}\;\; B_0,\\
   v_1=\,&\xi &  \quad  \mbox{on}\;\; \Gamma_0.
\end{array}
\right.
\end{equation}
  Clearly, these problems have unique solutions $v\in C^{m+\mu}(\overline{D})$ and $v_1\in C^{m+\mu}(\overline{B}_0)$, respectively. Define
  $\mathscr{D}(\rho,\eta):C^{m+\mu}(\Gamma_0)\to C^{m-1+\mu}(\Gamma_0)$ and $\mathscr{D}_1(\eta):C^{m+\mu}(\Gamma_0)\to C^{m-1+\mu}(\Gamma_0)$
  respectively as follows:
$$
  \mathscr{D}(\rho,\eta)\xi=\mathscr{B}(\eta)v, \quad \mathscr{D}_1(\eta)\xi=\mathscr{B}_1(\eta)v_1 \quad
  \mbox{for}\;\,\xi\in C^{m+\mu}(\Gamma_0).
$$
  The problem (5.14) is equivalent to the following problem: Find $\xi\in C^{m+\mu}(\Gamma_0)$ such that
\begin{equation}
  \mathscr{D}(\rho,\eta)\xi=\mathscr{D}_1(\eta)\xi.
\end{equation}
  We introduce a mapping $\mathscr{G}: O_{\delta}'\times O_{\delta'}''\times C^{m+\mu}(\overline{D})\times C^{m+\mu}(\Gamma_0)\to
  C^{m-1+\mu}(\Gamma_0)$ by defining
$$
  \mathscr{G}(\rho,\eta,u,\xi)=\mathscr{D}(\rho,\eta)\xi-\mathscr{D}_1(\eta)\xi \quad
  \mbox{for}\;\,\rho\in O_{\delta}',\;\,\rho\in O_{\delta'}'',\;\, \xi\in C^{m+\mu}(\Gamma_0).
$$
  It is clear that $\mathscr{G}\in C^{\infty}(O_{\delta}'\times O_{\delta'}'' C^{m+\mu}(\overline{D})\times C^{m+\mu}(\Gamma_0),C^{m-1+\mu}(\Gamma_0))$,
  and
$$
  \mathscr{G}(0,0,u_0,\xi_0)=0,
$$
  where $u_0=U(r,R)$ and $\xi_0$ represents the following constant function in $\Gamma_0$: $\xi_0(x)=V( K,R)$ for $x\in\Gamma_0$. A simple computation
  shows that for any $\zeta\in C^{m+\mu}(\Gamma_0)$, $L\zeta:=\partial_{\xi}\mathscr{G}(0,0,u_0,\xi_0)\zeta=\partial_{\nu}w|_{\Gamma_0}-
  \partial_{\nu}w_1|_{\Gamma_0}=\partial_rw_1|_{\Gamma_0}-\partial_rw|_{\Gamma_0}$, where $w,w_1$ are respectively the unique solutions of the
  following problems:
\begin{equation}
\left\{
\begin{array}{rll}
   \Delta w=&0 &  \quad  \mbox{for}\;\;  K<r<R,\\
   w=&0 &  \quad  \mbox{for}\;\; r=R,\\
   w=&\zeta &  \quad  \mbox{for}\;\; r=K,
\end{array}
\right.
  \qquad\qquad
\left\{
\begin{array}{rll}
   \Delta w_1=&0 &  \quad  \mbox{for}\;\; r< K,\\
   w_1=&\zeta &  \quad  \mbox{for}\;\; r= K.
\end{array}
\right.
\end{equation}
  From this fact it is not hard to see that if $L\zeta=0$ then $\zeta=0$. Indeed, if $L\zeta=0$ then by letting $z=w$ for $ K\leqslant r\leqslant R$
  and $z=w_1$ for $r< K$, we get a weak solution of the boundary value problem
$$
\left\{
\begin{array}{rl}
   \Delta z=&0   \quad  \mbox{in}\;\; B(0,R),\\
   z=&0   \quad  \mbox{on}\;\;\partial B(0,R),
\end{array}
\right.
$$
  which implies, by Green's second identity, that $z=0$ and, consequently, $\zeta=0$. This shows that ${\rm Ker}\,L=\{0\}$. Since $L$ is a sum of
  two Dirichlet-Neumann operators, it is a first-order pseudo-differential operator of elliptic type (cf. \cite{Esc}, \cite{Gri}), so that standard
  Schauder estimate applies to it: There exists a positive constant $C>0$ such that
$$
  \|\zeta\|_{C^{m+\mu}(\Gamma_0)}\leqslant C(\|L\zeta\|_{C^{m-1+\mu}(\Gamma_0)}+\|\zeta\|_{L^{\infty}(\Gamma_0)}),
   \quad \forall\zeta\in C^{m+\mu}(\Gamma_0).
$$
   Since ${\rm Ker}\,L=\{0\}$, this implies, by a standard argument, the following estimate:
\begin{equation}
  \|\zeta\|_{C^{m+\mu}(\Gamma_0)}\leqslant C\|L\zeta\|_{C^{m-1+\mu}(\Gamma_0)}, \quad \forall\zeta\in C^{m+\mu}(\Gamma_0).
\end{equation}
  For every $k\in{\mathbf{Z}}_+$ let $\{Y_{kl}(\omega)\}_{l=1}^{2k+1}$ be the normalized orthogonal basis (in $L^2({\mathbf{S}}^2)$ inner product) of
  the linear space of $k$-th order spherical harmonics. A simple computation shows that for any $\zeta\in C^{\infty}(\Gamma_0)$,
$$
  L\zeta(\omega)=\sum_{k=0}^{\infty}\sum_{l=1}^{2k+1}\frac{(2k\!+\!1)c_{kl}}{[1-(K/R)^{2k+1}]K}Y_{kl}(\omega) \quad
  \mbox{if}\;\; \zeta(K\omega)=\sum_{k=0}^{\infty}\sum_{l=1}^{2k+1}c_{kl}Y_{kl}(\omega).
$$
  From this expression of $L$ it is easy to see that for any $\eta\in C^{\infty}(\Gamma_0)$ the equation $L\zeta=\eta$ has a unique solution
  $\zeta\in C^{\infty}(\Gamma_0)$ (see the proofs of Lemma 6.3 and Corollary 6.4 in the next section for more details in this argument). It follows,
  by using the estimate (5.18) and a standard approximation argument, that also for any $\eta\in C^{m-1+\mu}(\Gamma_0)$ the equation
  $L\zeta=\eta$ has a unique solution $\zeta\in C^{m+\mu}(\Gamma_0)$\footnotemark[4].
\footnotetext[4]{Since $C^{\infty}(\Gamma_0)$ is not dense in $C^{m-1+\mu}(\Gamma_0)$, one might argue validity of the approximation argument. It is
  as follows: For any $\eta\in C^{m-1+\mu}(\Gamma_0)$ choose a sequence $\{\eta_j\}_{j=1}^{\infty}\subseteq C^{\infty}(\Gamma_0)$ such that it is
  bounded in $C^{m-1+\mu}(\Gamma_0)$ and converges to $\eta$ in $C^{m-1+\mu'}(\Gamma_0)$ for any $0<\mu'<\mu$ (e.g., choose a mollification
  sequence of $\eta$). For every $j$ let $\zeta_j\in C^{\infty}(\Gamma_0)$ be the unique solution of the equation $L\zeta_j=\eta_j$. The estimate
  (5.18) ensures the sequence $\{\zeta_j\}_{j=1}^{\infty}$ is bounded in $C^{m+\mu}(\Gamma_0)$ and converges to a function $\zeta$ in
  $C^{m+\mu'}(\Gamma_0)$ for any $0<\mu'<\mu$, which implies that $\zeta$ is a solution of the equation $L\zeta=\eta$. Boundedness of the
  sequence $\{\zeta_j\}_{j=1}^{\infty}$ in $C^{m+\mu}(\Gamma_0)$ implies $\zeta\in C^{m+\mu}(\Gamma_0)$.}
  This shows that $L=\partial_{\xi}\mathscr{G}(0,0,u_0,\xi_0): C^{m+\mu}(\Gamma_0)
  \to C^{m-1+\mu}(\Gamma_0)$ is a (linear and topological) isomorphism. Hence, by applying the implicit function theorem in Banach spaces, we see that
  by choosing $\delta,\delta'>0$ further small when necessary, there exists a smooth mapping $\chi:O_{\delta}'\times O_{\delta'}''\times
  C^{m+\mu}(\overline{D})\subseteq C^{m+\mu}({\mathbf{S}}^2)\times C^{m+\mu}({\mathbf{S}}^2)\times C^{m+\mu}(\overline{D})\to C^{m+\mu}(\Gamma)$ such
  that $\chi(0,0,u_0)=\xi_0$ and for any $\rho\in O_{\delta}'$, $\eta\in O_{\delta'}''$ and $u\in C^{m+\mu}(\overline{D})$, $\xi=\chi(\rho,\eta,u)$
  is the unique solution of the equation $\mathscr{G}(\rho,\eta,u,\xi)=0$ in a small neighborhood of $\xi_0$ in $C^{m+\mu}(\Gamma)$. This proves unique
  solvability of the equation (5.16) and smoothness of the solution map $\xi=\chi(\rho,\eta,u)$. As a result, the solution map $(\rho,\eta,u)\mapsto
  (v,v_1)$ (from $O_{\delta}'\times O_{\delta'}''\times C^{m+\mu}(\overline{D})\subseteq C^{m+\mu}({\mathbf{S}}^2)\times C^{m+\mu}({\mathbf{S}}^2)\times
  C^{m+\mu}(\overline{D})$ to $C^{m+\mu}(\overline{D})\times C^{m+\mu}(\overline{B}))$) of the problem (5.14) is smooth. This proves Lemma 5.4 and
  also completes the proof of Lemma 5.3. $\quad\Box$
\medskip

   {\em Remark}.\ \ We note that all H\"{o}lder spaces appearing in this section can be replaced with corresponding little H\"{o}lder spaces. In the
   next section we shall use this fact without further explanation.

\section{A free boundary problem modeling necrotic tumor growth}
\setcounter{equation}{0}

\hskip 2em
  In this section we use Theorem 3.4 to study asymptotic stability of the radial stationary solution of the free boundary problem (1.3) in the case
  that $f$, $g$ are discontinuous functions given by (1.5). We shall prove that there exists a positive constant $\gamma^*$ such that if $\gamma>
  \gamma^*$ then the radial stationary solution is asymptotically stable module translations in ${\mathbf{R}}^3$, whereas if $\gamma<\gamma^*$ then
  it is unstable.

  We first recall that radial stationary solution $(\sigma_s,\pi_s,\Omega_s)$ of the problem (1.3) is given by
\begin{equation}
  \sigma_s(r)=U(r,R_s), \qquad  \pi_s(r)=\frac{\gamma}{R_s}+V(r,R_s), \qquad \Omega_s=B(0,R_s),
\end{equation}
  where $U$, $V$ are as in Section 5, and $R_s$ is the root of the equation $\pi_s'(R_s)=0$ or $\partial_rV(R_s,R_s)=0$, i.e.,
$$
  \frac{1}{3}a\tilde{\sigma}K^3(R_s)+a\int_{K(R_s)}^{R_s}U(\eta,R_s)\eta^2d\eta=\frac{1}{3}(a\tilde{\sigma}+b)R_s^3.
$$
  Since $0<\tilde{\sigma}<1$, by Lemma 4.2 of \cite{Cui2} we know that this equation has a unique solution $R_s>R^*$, which means the problem (1.3)
  with $f$, $g$ given in (1.5) has a unique radial stationary solution.

  In order to apply Theorem 3.4, let us first reduce the problem (1.3) into a differential equation in a Banach manifold. Let $m$ be a positive
  integer $\geqslant 2$ and $0<\mu<1$. Let $\mathfrak{M}:=\dot{\mathfrak{D}}^{m+\mu}({\mathbf{R}}^3)$ and $\mathfrak{M}_0:=
  \dot{\mathfrak{D}}^{m+\mu+3}({\mathbf{R}}^3)$. Given $\Omega\in\mathfrak{M}_0$, we have seen that the problem (5.1) has a unique solution
  $(\sigma,\pi_0)$ satisfying the following properties:
$$
  \sigma,\pi_0\in W^{2,p}(\Omega) \;(\forall p\in [1,\infty)), \quad
  \sigma,\pi_0\in\dot{C}^{\infty}(\Omega_{\rm liv}\cup\Omega_{\rm nec}) \quad \mbox{and} \quad
  \sigma,\pi_0\in\dot{C}^{m+3+\mu}(\Omega_{\rm liv}\cup\partial\Omega).
$$
  It follows that $\partial_{\bfn}\pi_0|_{\partial\Omega}\in \dot{C}^{m+\mu+2}(\partial\Omega)$.
  We define $\mathscr{F}_0:\mathfrak{M}_0\to\mathcal{T}_{\mathfrak{M}_0}(\mathfrak{M})$ by setting
$$
  \mathscr{F}_0(\Omega)=-\partial_{\bfn}\pi_0|_{\partial\Omega}, \quad  \forall\Omega\in\mathfrak{M}_0.
$$
  Now let $\mathscr{F}:\mathfrak{M}_0\to\mathcal{T}_{\mathfrak{M}_0}(\mathfrak{M})$ be the vector field in $\mathfrak{M}$ introduced in Section 4 in the study of
  Hele-Shaw problem (for $n=3$), and define $\mathscr{G}:\mathfrak{M}_0\to\mathcal{T}_{\mathfrak{M}_0}(\mathfrak{M})$ as follows:
\begin{equation}
  \mathscr{G}(\Omega)=\gamma\mathscr{F}(\Omega)+\mathscr{F}_0(\Omega), \quad  \forall\Omega\in\mathfrak{M}_0.
\end{equation}
  Then the problem (1.3) reduces into the following differential equation in the Banach manifold $\mathfrak{M}$:
\begin{equation}
\left\{
\begin{array}{ll}
   \Omega'(t)=\mathscr{G}(\Omega(t)), &\quad  t>0,\\
    \Omega(0)=\Omega_0. &
\end{array}
\right.
\end{equation}
  The fact that $(\sigma_s,\pi_s,\Omega_s)$ is a stationary solution of the problem (1.3) implies that $\Omega_s$ is a stationary solution
  of the equation $\Omega'=\mathscr{G}(\Omega)$.

  It is easy to check the equation $\Omega'=\mathscr{G}(\Omega)$ is quasi-invariant under the translation group action $(G_{tl},p)$. Hence,
  for all $a\in G_{tl}$, $\Omega=p(\Omega_s,a)$ are stationary solutions of the equation $\Omega'=\mathscr{G}(\Omega)$ and they form a three
  dimensional manifold $\mathcal{M}_c$.

  Next we consider the representation of the problem (6.3) in a regular local chart $(\mathcal{U},\varphi)$ of $\mathfrak{M}$ at the point $\Omega_s$.
  Let $\delta>0$ be a sufficiently small number. We denote
$$
  X=\dot{C}^{m+\mu}({\mathbf{S}}^2),  \qquad  X_0=\dot{C}^{m+\mu+3}({\mathbf{S}}^2),  \qquad  X_1=\dot{C}^{m+\mu+2}({\mathbf{S}}^2),
$$
$$
  O_{\delta}=\{\rho\in X,\|\rho\|_{X}<\delta\},  \qquad O'_{\delta}=\{\rho\in X_0,\|\rho\|_{X_0}<\delta\},
$$
$$
  \mathcal{U}=\{\Omega_{\rho}:\rho\in O_{\delta}\},  \qquad \mathcal{U}'=\{\Omega_{\rho}:\rho\in O'_{\delta}\},
$$
  where
\begin{equation}
  \Omega_{\rho}=\{x\in{\mathbf{R}}^3: r<R_s[1+\rho(\omega)]\}.
\end{equation}
  It is clear that $\Omega_{\rho}\in\mathfrak{M}$. Define $\varphi:\mathcal{U}\to X$ by letting $\varphi(\Omega_{\rho})=\rho$, $\forall\rho\in
  O_{\delta}$. Then $(\mathcal{U},\varphi)$ is a regular local chart of $\mathfrak{M}$ at the point $\Omega_s$, with base space
  $X$. We denote by $F$, $F_0$ and $G$ the representations of the vector fields $\mathscr{F}$, $\mathscr{F}_0$ and $\mathscr{G}$, respectively, in
  this local chart, i.e., for any $\rho\in O'_{\delta}$,
$$
  F(\rho)=\varphi'(\Omega_{\rho})\mathscr{F}(\Omega_{\rho}), \quad F_0(\rho)=\varphi'(\Omega_{\rho})\mathscr{F}_0(\Omega_{\rho})
   \quad \mbox{and} \quad G(\rho)=\varphi'(\Omega_{\rho})\mathscr{G}(\Omega_{\rho}).
$$
  Then $G(\rho)=\gamma F(\rho)+F_0(\rho)$, and representation of the problem (6.3) in the local chart $(\mathcal{U},\varphi)$ is the following
  initial value problem in the Banach space $X$:
\begin{equation}
\left\{
\begin{array}{ll}
   \rho'(t)=G(\rho(t)), &\quad  t>0,\\
   \rho(0)=\rho_0, &
\end{array}
\right.
\end{equation}
  where $\rho_0\in O'_{\delta}$ is the function such that $\Omega_0=\Omega_{\rho_0}$. From Section 4 we see that $F\in C^{\infty}(O'_{\delta},X)$.
  Lemma 5.3 shows that $F_0\in C^{\infty}(O'_{\delta},X_1)\subseteq C^{\infty}(O'_{\delta},X)$. Hence we have
\medskip

  {\bf Lemma 6.1}\ \ {\em $G\in C^{\infty}(O'_{\delta},X)$.} $\quad\Box$
\medskip

  {\bf Corollary 6.2}\ \ {\em The differential equation $(6.3)$ is of parabolic type.}
\medskip

  {\em Proof.}\ \ From the discussion in Section 4 we see that for any $\rho\in O'_{\delta}$, $F'(\rho)$ is a sectorial operator in $X$ with domain
  $X_0$. Lemma 5.3 ensures that for any $\rho\in O'_{\delta}$, $F_0'(\rho)\in L(X_0,X_1)$. Since $G'(\rho)=\gamma F'(\rho)+F_0'(\rho)$ and $X_1$
  is an intermediate space between $X_0$ and $X$, by a well-known perturbation theorem for sectorial operators it follows that $G'(\rho)$ is also
  a sectorial operator in $X$ with domain $X_0$. Besides, it is easy to check that the graph norm of ${\rm Dom}\,G'(\rho)=X_0$ is equivalent to the
  norm of $X_0$. Hence the desired assertion follows. $\quad\Box$

  Following \cite{FriRei2} and \cite{Cui4}, we compute $G'(0)$ as follows: Let
$$
  \rho(\omega)=\varepsilon\xi(\omega), \quad \eta(\omega)=\varepsilon\zeta(\omega), \quad \sigma(r,\omega)=\sigma_s(r)+\varepsilon u(r,\omega).
$$
  Substituting these expressions into (5.5) and the equation $\Delta\sigma=0$ for $r< K_s$, and comparing coefficients of first-order terms of
  $\varepsilon$, we obtain the following equations:
\begin{equation}
\left\{
\begin{array}{rll}
   \Delta u=&u &  \quad  \mbox{for}\;\;  K_s<r<R_s,\\
   \Delta u=&0 &  \quad  \mbox{for}\;\; r< K_s,\\
   u=&-R_s\sigma_s(R_s)\xi(\omega) &  \quad  \mbox{for}\;\; r=R_s,\\
   u=&0 &  \quad  \mbox{for}\;\; r= K_s,\\
   \partial_r^+u=&-\hat{\sigma} K_s\zeta(\omega) &  \quad  \mbox{for}\;\; r= K_s.
\end{array}
\right.
\end{equation}
  Here $K_s=K(R_s)$ and for $r= K_s$, $\partial_r^+u=\frac{\partial u}{\partial r}( K_s^+,\omega)$. Similarly, by letting $\pi(r,\omega)=\pi_s(r)+
  \varepsilon v(r,\omega)$, we get the following equations for $v$:
\begin{equation}
\left\{
\begin{array}{rll}
   \Delta v=&-au &  \quad  \mbox{for}\;\;  K_s<r<R_s,\\
   \Delta v=&0 &  \quad  \mbox{for}\;\; r< K_s,\\
   v=&\displaystyle-\frac{\gamma}{R_s}\Big(\xi(\omega)+\frac{1}{2}\Delta_{\omega}\xi(\omega)\Big) &  \quad  \mbox{for}\;\; r=R_s,\\
   v^+=&v^- &  \quad  \mbox{for}\;\; r= K_s,\\
   \partial_r^+v-\partial_r^-v=&a(\hat{\sigma}-\tilde{\sigma}) K_s\zeta(\omega) &  \quad  \mbox{for}\;\; r= K_s.
\end{array}
\right.
\end{equation}
  Here for $r=K_s$, $v^{\pm}=v( K_s^{\pm},\omega)$ and $\partial_r^{\pm}v=\frac{\partial v}{\partial r}( K_s^{\pm},\omega)$. After solving these
  equations (with $u$, $v$ and $\zeta$ being unknown functions) for any given function $\xi=\xi(\omega)$, we then have
\begin{equation}
  G'(0)\xi(\omega)=-\frac{\partial v}{\partial r}(R_s,\omega)+g(1)R_s\xi(\omega), \quad \forall\xi\in X_0
\end{equation}
  (cf. (4.4) in \cite{Cui4}, but be aware that here the perturbation of the sphere $r=R_s$ is given by $r=R_s[1+\varepsilon\xi(\omega)]$, not as in
  \cite{Cui4} given by $r=R_s+\varepsilon\xi(\omega)$).

  We now use spherical harmonics expansions of functions in ${\mathbf{S}}^2$ to solve the problems (6.6) and (6.7). Hence let $\{Y_{kl}(\omega):
  k=0,1,2,\cdots,\; l=1,2,\cdots,2k+1\}$ be a normalized orthogonal basis of $L^2({\mathbf{S}}^2)$ consisting spherical harmonics on ${\mathbf{S}}^2$,
  where for every $k\in{\mathbf{Z}}_+$, $Y_{kl}(\omega)$, $l=1,2,\cdots,2k+1$, are spherical harmonics of degree $k$, so that
$$
  \Delta_{\omega}Y_{kl}(\omega)=-\lambda_kY_{kl}(\omega), \quad k=0,1,2,\cdots,\;\; l=1,2,\cdots,2k+1.
$$
  where $\lambda_k=k(k+1)$, $k=0,1,2,\cdots$. A simple computation shows that if a given function $\xi\in C^{\infty}({\mathbf{S}}^2)$ has a spherical
  harmonics expansion
\begin{equation}
  \xi(\omega)=\sum_{k=0}^{\infty}\sum_{l=1}^{2k+1}c_{kl}Y_{kl}(\omega),
\end{equation}
  then the solution $(u,\zeta)$ of the problem (6.6) is given by
\begin{equation}
  u(r,\omega)=\left\{
\begin{array}{cl}
  \displaystyle -R_s\sigma_s(R_s)\sum_{k=0}^{\infty}\sum_{l=1}^{2k+1}\Big(\frac{r}{R_s}\Big)^k\bar{u}_k(r)c_{kl}Y_{kl}(\omega) &  \quad
   \mbox{for}\;\;  K_s<r\leqslant R_s\\
   0 &  \quad  \mbox{for}\;\; r\leqslant K_s
\end{array}
\right.
\end{equation}
  and $\zeta(\omega)=-\partial_r^+u/\hat{\sigma} K_s$, where $\bar{u}_k$ is the unique solution of the following boundary value problem:
\begin{equation}
  \left\{
\begin{array}{l}
  \displaystyle \bar{u}_k''(r)+\frac{2(k\!+\!1)}{r}\bar{u}_k'(r)=\bar{u}_k(r) \quad  \mbox{for}\;\;  K_s<r<R_s,\\
   \bar{u}_k( K_s)=0,   \quad  \bar{u}_k(R_s)=1.
\end{array}
\right.
\end{equation}
  Substituting these expressions of $u$ and $\zeta$ into (6.7), we easily obtain the following solution of that problem:
\begin{equation}
  v(r,\omega)=\sum_{k=0}^{\infty}\sum_{l=1}^{2k+1}\Big[\frac{\gamma}{2R_s}(k-1)(k+2)
  +aR_s\sigma_s'(R_s)\bar{v}_k(r)\Big]\Big(\frac{r}{R_s}\Big)^kc_{kl}Y_{kl}(\omega),
\end{equation}
  where $\bar{v}_k$ is the unique solution of the following boundary value problem:
\begin{equation}
  \left\{
\begin{array}{l}
  \displaystyle \bar{v}_k''(r)+\frac{2(k\!+\!1)}{r}\bar{v}_k'(r)=\bar{u}_k(r) \quad  \mbox{for}\;\;  K_s<r<R_s,\\
  \displaystyle \bar{v}_k''(r)+\frac{2(k\!+\!1)}{r}\bar{v}_k'(r)=0  \quad \mbox{for}\;\; r< K_s,\\
   \bar{v}_k(R_s)=0,   \\
   \bar{v}_k( K_s^+)=\bar{v}_k( K_s^-),  \\
  \displaystyle \bar{v}_k'( K_s^+)-\bar{v}_k'( K_s^-)=\frac{\hat{\sigma}-\tilde{\sigma}}{\hat{\sigma}}\bar{u}_k'( K_s^+).
\end{array}
\right.
\end{equation}
  From (6.13) we get the following relation:
\begin{equation}
  \bar{v}_k'(R_s)=\frac{\hat{\sigma}-\tilde{\sigma}}{\hat{\sigma}}\bar{u}_k'(K^+)\Big(\frac{K_s}{R_s}\Big)^{2(k+1)}
  +\int_{K_s}^{R_s}\bar{u}_k(\tau)\Big(\frac{\tau}{R_s}\Big)^{2(k+1)}d\tau.
\end{equation}
  Now let
\begin{equation}
  a_k(\gamma)=-\frac{\gamma}{2R_s^2}k(k-1)(k+2)-aR_s\sigma_s'(R_s)\bar{v}_k'(R_s)+g(1)R_s.
\end{equation}
  Then from (6.8), (6.12), $(6.13)_3$ and (6.14) we obtain the following preliminary result:
\medskip

  {\bf Lemma 6.3}\ \ {\em $G'(0)$ is a Fourier multiplier in the sense that if $\xi\in C^{\infty}({\mathbf{S}}^2)$ has expansion $(6.9)$, then
  $G'(0)\xi=\displaystyle\sum_{k=0}^{\infty}\sum_{l=1}^{2k+1}a_k(\gamma)c_{kl}Y_{kl}(\omega)\in C^{\infty}({\mathbf{S}}^2)$.}
  $\quad\Box$
\medskip

  {\bf Corollary 6.4}\ \ {\em $\sigma(G'(0))=\{a_k(\gamma):k=0,1,2,\cdots\}$.}
\medskip

  {\em Proof.}\ \ Since $G'(0)\in L(X_0,X)$ and it is a sectorial operator in $X$ with domain $X_0$, the inverse mapping theorem implies that for
  any $\lambda\in\rho(G'(0))$ we have $[\lambda I-G'(0)]^{-1}\in L(X,X_0)$. As a consequence, $R(\lambda,G'(0))=[\lambda I-G'(0)]^{-1}$ is a compact
  linear operator in $X$, so that its spectrum contains only eigenvalues. It follows that $\sigma(G'(0))$ also contains only eigenvalues. Next, for
  every $s\geqslant 0$ let $H^s({\mathbf{S}}^2)$ be the standard Sobolev space on ${\mathbf{S}}^2$ with index $s$. We know that $H^s({\mathbf{S}}^2)$ has
  an equivalent norm $\|\xi\|_{H^s({\mathbf{S}}^2)}=\displaystyle\Big[\sum_{k=0}^{\infty}\sum_{l=1}^{2k+1}(1\!+\!\lambda_k)^s|c_{kl}|^2\Big]^{1/2}$,
  if $\xi\in H^s({\mathbf{S}}^2)$ has the expression (6.9) as an element of $L^2({\mathbf{S}}^2)$. Since $C^{\infty}({\mathbf{S}}^2)=\displaystyle
  \bigcap_{s\geqslant 0}H^s({\mathbf{S}}^2)$ and it is dense in every $H^s({\mathbf{S}}^2)$, $s\geqslant 0$, Lemma 6.3 shows that the operator $G'(0):
  \dot{C}^{m+3+\mu}({\mathbf{S}}^2)\to\dot{C}^{m+\mu}({\mathbf{S}}^2)$ can be uniquely extended into a bounded linear operator from $H^3({\mathbf{S}}^2)$
  to $L^2({\mathbf{S}}^2)$, and after extension we have the relation $\sigma(G'(0))=\{a_k(\gamma):k=0,1,2,\cdots\}$. Moreover, the eigenspace
  corresponding to the eigenvalue $a_k(\gamma)$ is ${\rm span}\{Y_{kl}(\omega):l=1,2,\cdots,2k\!+\!1\}\subseteq C^{\infty}({\mathbf{S}}^2)$, i.e, all
  eigenfunctions are smooth. Hence, since if $\xi\in\dot{C}^{m+3+\mu}({\mathbf{S}}^2)$ is an eigenvector of the operator $G'(0):
  \dot{C}^{m+3+\mu}({\mathbf{S}}^2)\to\dot{C}^{m+\mu}({\mathbf{S}}^2)$ then it is also an eigenvector of the operator $G'(0):H^3({\mathbf{S}}^2)\to
  L^2({\mathbf{S}}^2)$, we obtain the desired assertion. $\quad\Box$
\medskip

  It is clear that $a_0(\gamma)$ and $a_1(\gamma)$ are independent of $\gamma$. Hence we re-denote them as $a_0$ and $a_1$, respectively, i.e.,
\begin{equation}
  a_0=g(1)R_s-aR_s\sigma_s'(R_s)\bar{v}_0'(R_s), \quad  a_1=g(1)R_s-aR_s\sigma_s'(R_s)\bar{v}_1'(R_s).
\end{equation}
   For $k\geqslant 2$ we denote
\begin{equation}
  \gamma_k=\frac{2R_s^3}{k(k-1)(k+2)}[g(1)-a\sigma_s'(R_s)\bar{v}_k'(R_s)].
\end{equation}
  Then from (6.15) we have
\begin{equation}
  a_k(\gamma)=-\frac{1}{2R_s^2}k(k-1)(k+2)(\gamma-\gamma_k), \qquad k=2,3,\cdots.
\end{equation}
  We shall prove $a_0<0$, $a_1=0$ and $\gamma_k>0$ for $k\geqslant 2$. For this purpose we need the following lemma:
\medskip

  {\bf Lemma 6.5}\ \ {\em For the solution of the problem $(6.11)$ we have the following assertions:

  $(1)$\ \ $0<\bar{u}_k(r)<1$ and $\bar{u}_k'(r)>0$ for $ K_s<r<R_s$.

  $(2)$\ \ If $k>l$ then $\bar{u}_k(r)>\bar{u}_l(r)$ for $ K_s<r<R_s$, and
  $\bar{u}_k'( K_s)\geqslant\bar{u}_l'( K_s)$, $\bar{u}_k'(R_s)\leqslant\bar{u}_l'(R_s)$.

  $(3)$\ \ If $k>l$ then $\bar{u}_k(r)(r/R_s)^{k+1}<\bar{u}_l(r)(r/R_s)^{k+1}$ for $ K_s<r<R_s$.}
\medskip

  {\em Proof.}\ \ The assertion $0<\bar{u}_k(r)<1$ for $ K_s<r<R_s$ is an immediate consequence of the maximum principle. Note that this
  assertion joint with the boundary value conditions $\bar{u}_k( K_s)=0$ and $\bar{u}_k(R_s)=1$ implies that $\bar{u}_k'( K_s)
  \geqslant0$ and $\bar{u}_k'(R_s)\geqslant 0$. Next we let $w_k(r)=\bar{u}_k'(r)$. A simple computation shows that $w_k$ satisfies the following
  equation:
$$
  w_k''(r)+\frac{2(k\!+\!1)}{r}w_k'(r)-\Big(\frac{2(k\!+\!1)}{r^2}+1\Big)w_k(r)=0 \quad  \mbox{for}\;\;  K_s<r<R_s.
$$
  Since $w_k( K_s)\geqslant0$ and $w_k(R_s)\geqslant 0$, again by the maximum principle we see that $w_k(r)>0$ for $ K_s<r<R_s$. This
  proves the assertion (1). From the property $\bar{u}_k'(r)>0$ for $ K_s<r<R_s$ it follows that if $k>l$ then
$$
  \bar{u}_k''(r)+\frac{2(l\!+\!1)}{r}\bar{u}_k'(r)-\bar{u}_k(r)<0 \quad  \mbox{for}\;\;  K_s<r<R_s.
$$
  Hence by the maximum principle we obtain $\bar{u}_k(r)>\bar{u}_l(r)$ for $ K_s<r<R_s$ and $k>l$, which easily implies that $\bar{u}_k'( K_s)
  \geqslant\bar{u}_l'( K_s)$ and $\bar{u}_k'(R_s)\leqslant\bar{u}_l'(R_s)$ for $k>l$. This proves the assertion (2). Finally we
  let $z_k(r)=\bar{u}_k(r)(r/R_s)^{k+1}$. It can be easily seen that $z_k(r)$ is a solution of the following problem:
$$
  \left\{
\begin{array}{l}
  \displaystyle z_k''(r)-\Big(\frac{k(k\!+\!1)}{r^2}+1\Big)z_k(r)=0 \quad  \mbox{for}\;\;  K_s<r<R_s,\\ [0.2cm]
   z_k( K_s)=0,   \quad  z_k(R_s)=1.
\end{array}
\right.
$$
  Since $z_k>0$ for $ K_s<r<R_s$, a similar argument as in the proof of the assertion (2) shows that if $k>l$ then $z_k(r)<z_l(r)$ for  for
  $ K_s<r<R_s$. This proves the assertion (3) and completes the proof of the lemma. $\quad\Box$
\medskip

  {\bf Lemma 6.6}\ \ {\em We have the following assertions:

  $(1)$\ \ $a_0<0$ and $a_1=0$.

  $(2)$\ \ $\gamma_k>0$ for all $k=2,3,\cdots$, and $\gamma_k\sim 2R_s^3g(1)k^{-3}$ as $k\to\infty$.}
\medskip

  {\em Proof.}\ \ From (6.11) we have
$$
  \bar{u}_k'(R_s)=\bar{u}_k'( K_s)\Big(\frac{ K_s}{R_s}\Big)^{2(k+1)}+
  \int_{ K_s}^{R_s}\bar{u}_k(\tau)\Big(\frac{\tau}{R_s}\Big)^{2(k+1)}d\tau.
$$
  Hence the relation (6.14) can be rewritten as follows:
$$
  \bar{v}_k'(R_s)=\frac{\hat{\sigma}-\tilde{\sigma}}{\hat{\sigma}}\bar{u}_k'(R_s)
  +\frac{\tilde{\sigma}}{\hat{\sigma}}\int_{ K_s}^{R_s}\bar{u}_k(\tau)\Big(\frac{\tau}{R_s}\Big)^{2(k+1)}d\tau.
$$
  Using this relation and the assertions (2) and (3) of Lemma 6.5, we conclude that
$$
  \bar{v}_k'(R_s)<\bar{v}_l'(R_s) \quad \mbox{if}\;\; k>l.
$$
  Hence, all the desired assertions will follow if we prove that $a_1=0$ and $\bar{v}_k'(R_s)=O(1/k)$ as $k\to\infty$. The proof that
  $\bar{v}_k'(R_s)=O(1/k)$ as $k\to\infty$ is easy and is omitted. In what follows we prove $a_1=0$.

  It is easy to check
\begin{equation}
  \bar{u}_1(r)=\frac{R_s\sigma_s'(r)}{r\sigma_s'(R_s)} \quad \mbox{for}\;\,  K_s\leqslant r\leqslant R_s.
\end{equation}
  In what follows we prove
\begin{equation}
  \bar{v}_1(r)=-\frac{R_s\pi_s'(r)}{ar\sigma_s'(R_s)} \quad \mbox{for}\;\,  K_s\leqslant r\leqslant R_s.
\end{equation}
  We use the notation $q(r)$ to denote the function on the right-hand side of the above relation. To prove the above relation, we only need to
  show that $q(r)$ is a solution of the following problem:
$$
  \left\{
\begin{array}{l}
  \displaystyle q''(r)+\frac{4}{r}q'(r)=\bar{u}_1(r) \quad  \mbox{for}\;\;  K_s<r<R_s,\\ [0.2cm]
  \displaystyle  q(R_s)=0,   \quad  q'(K_s^+)=\frac{\hat{\sigma}-\tilde{\sigma}}{\hat{\sigma}}\bar{u}_1'(K_s^+).
\end{array}
\right.
$$
  The equation in the first line is easy to check, and the boundary value condition $q(R_s)=0$ is clear. Since
$$
  q'(r)=-\frac{R_s}{a\sigma_s'(R_s)}\Big[\frac{\pi_s''(r)}{r}-\frac{\pi_s'(r)}{r^2}\Big]
  =\frac{R_s}{a\sigma_s'(R_s)}\Big[\frac{3\pi_s'(r)}{r^2}+\frac{g(\sigma_s(r))}{r}\Big],
$$
  we have
$$
  q'( K_s^+)=\frac{R_s}{a\sigma_s'(R_s)}\Big[\frac{3\pi_s'( K_s)}{ K_s^2}+\frac{g(\hat{\sigma}^+)}{ K_s}\Big]
  =\frac{3R_s\pi_s'(K_s)}{a\sigma_s'(R_s)K_s^2}.
$$
  From the equation $\Delta\pi_s=-b$ (for $r<K_s$) we see that $\pi_s'(K_s)=(1/3)bK_s=(a/3)(\hat{\sigma}-\tilde{\sigma})K_s$, and it is clear
  that $\bar{u}_1'(K_s^+)=R_s\sigma_s''(K_s^+)/K_s\sigma_s'(R_s)=\hat{\sigma}R_s/K_s\sigma_s'(R_s)$. Combining these relations, we see that the
  the boundary value condition $q'(K_s^+)=\frac{\hat{\sigma}-\tilde{\sigma}}{\hat{\sigma}}\bar{u}_1'( K_s^+)$ is also satisfied. Hence (6.19) is
  true. The assertion $a_1=0$ is an immediate consequence of the relation (6.19) and the fact that $\pi_s''(R_s)=-g(1)$. $\quad\Box$
\medskip

  {\bf Lemma 6.7}\ \ {\em Let $\gamma^*=\max\{\gamma_k:k=2,3,\cdots\}$. The following assertions hold:

  $(1)$\ \ If $\gamma>\gamma^*$ then $\sup\{{\rm Re}\lambda:\lambda\in\sigma(G'(0))\backslash\{0\}\}<0$ and ${\rm Ker}\,G'(0)=
  {\rm span}\{Y_{11}(\omega),Y_{12}(\omega),Y_{13}(\omega)\}$, so that $\dim{\rm Ker}\,G'(0)=3$. If instead $0<\gamma<\gamma^*$ then
  $\sup\{{\rm Re}\lambda:\lambda\in\sigma(G'(0))\}>0$.

  $(4)$\ \ Let $\gamma>\gamma^*$. Then ${\rm Range}\,G'(0)$ is closed, and $X={\rm Ker}\,G'(0)\oplus{\rm Range}\,G'(0)$.}
\medskip

  {\em Proof.}\ \ Assertions in (1) are immediate consequences of Corollary 6.4, the expression (6.18) and Lemma 6.6. To prove the assertion (2),
  we note that $G'(0)=\gamma F'(0)+F_0'(0)$. From \cite{EscSim2} we know that $F'(0)$ is a third-order elliptic pseudo-differential operator on the
  sphere ${\mathbf{S}}^2$ (cf. also Lemma 4.1 (1) and note that the Dirichlet-Neumann operator $D$ is a first-order elliptic pseudo-differential
  operator on the sphere, cf. \cite{Esc}), and Corollary 5.4 shows that $F_0'(0)$ is a lower-order perturbation (actually, a similar discussion as
  in \cite{Esc} and \cite{EscSim2} shows that it is a first-order pseudo-differential operator on the sphere). It follows that standard
  $C^{\mu}$-estimates work for $G'(0)$ and, consequently, the Fredholm alteration principle applies to it, by a similar argument as in the proof of
  Lemma 5.4. Hence the assertion (2) follows. $\quad\Box$
\medskip

  Lemmas 6.1, 6.6, 6.7 and Corollaries 6.2, 6.4 show that Theorem 3.4 (with $N=1$) applies to the equation $(6.3)_1$. Hence, by applying Theorem 3.4
  and the linearized instability criterion for parabolic equations in Banach spaces (i.e. Theorem 9.1.3 of \cite{Lun2}) we get the following result:
\medskip

  {\bf Theorem 6.8}\ \ {\em Assume that $\gamma>\gamma^*$. Let $\mathcal{M}_c$ be the $3$-dimensional submanifold of $\mathfrak{M}_0$ consisting
  of all spheres in ${\mathbf{R}}^3$ of radius $R_s$. We have the following assertions:

  $(1)$\ \ There is a neighborhood $\mathcal{O}$ of $\mathcal{M}_c$ in $\mathfrak{M}_0$ such that for any $\Omega_0\in\mathcal{O}$, the initial value
  problem $(6.3)$ has a unique solution $\Omega\in C([0,\infty),\mathfrak{M}_0)\cap C^1((0,\infty),\mathfrak{M}_0)$.

  $(2)$\ \ There exists a submanifold $\mathcal{M}_s$ of $\mathfrak{M}_0$ of codimension $3$ passing $\Omega_s=B(0,R_s)$ such that for any
  $\Omega_0\in\mathcal{M}_s$, the solution of the problem $(6.3)$ satisfies $\displaystyle\lim_{t\to\infty}\Omega(t)=\Omega_s$ and, conversely, if
  the solution of $(6.3)$ satisfies this property then $\Omega_0\in\mathcal{M}_s$.

  $(3)$\ \ For any $\Omega_0\in\mathcal{O}$ there exist unique $x_0\in{\mathbf{R}}^3$ and $Q_0\in\mathcal{M}_s$ such that $\Omega_0=x_0+Q_0$ and for
  the solution $\Omega=\Omega(t)$ of $(6.3)$ we have
$$
   \lim_{t\to\infty}\Omega(t)=B(x_0,R_s).
$$
   If on the contrary $0<\gamma<\gamma^*$ then the radial stationary solution of the problem $(1.3)$ is unstable.$\quad\Box$}
\medskip

   {\bf Acknowledgement}.\hskip 1em This work is supported by the National Natural  Science Foundation of China under grant numbers 11571381.


\medskip


\begin{thebibliography}{99}

\bibitem{BazFri1} B. Bazaliy and A. Friedman, A free boundary problem for an elliptic-parabolic system: application to a model of tumor growth,
  \textit{Comm. Part. Diff. Equa.}, \textbf{28}(2003), pp. 517--560.
\bibitem{BazFri2} B. Bazaliy and A. Friedman, Global existence and asymptotic stability for an elliptic-parabolic free boundary problem: an
  application to a model of tumor growth, \textit{Indiana Univ. Math. J.}, \textbf{52}(2003), pp. 1265--1304.
\bibitem{ByrC1} H. Byrne and M. Chaplain, Growth of nonnecrotic tumors in the presence and absence of inhibitors, \textit{Math. Biosci.},
  \textbf{130}(1995), pp. 151--181.
\bibitem{ByrC2} H. Byrne and M. Chaplain, Growth of necrotic tumors in the presence and absence of inhibitors, \textit{Math. Biosci.},
  \textbf{135}(1996), pp. 187--216.
\bibitem{Caff} L. A. Caffarelli, The regularity of free boundaries in higher dimension, \textit{Acta Math.}, \textbf{139}(1977), pp. 155--184.
\bibitem{Chen} X. Chen, The Hele-Shaw problem and area-preserving curve-shortening motions, \textit{Arch. Rat. Mech. Anal.}, \textbf{123}(1993),
  pp. 117--151.
\bibitem{ChenHY} X. Chen, J. Hong and F. Yi, Existence, uniqueness, and regularity of classical solutions of the Mullins-Sekerka problem,
  \textit{Comm. Part. Diff. Equa.}, \textbf{21}(1996), pp. 1705--1727.
\bibitem{Cui1} S. Cui, Analysis of a free boundary problem modeling tumor growth, \textit{Acta Math. Sinica Engl. Ser.}, \textbf{21}(2005),
  pp. 1071--1082.
\bibitem{Cui2} S. Cui, Formation of necrotic cores in the growth of tumors: analytic results, \textit{Acta Math. Sci. Ser. B Engl. Ed.},
  \textbf{26B}(2006), pp. 781--796.
\bibitem{Cui3} S. Cui, Lie group action and stability analysis of stationary solutions for a free boundary problem modeling tumor growth,
  \textit{J. Diff. Equa.},  \textbf{246}(2009), pp. 1845--1882.
\bibitem{Cui4} S. Cui and J. Escher, Asymptotic behavior of solutions of a multidimensional moving boundary problem modeling
  tumor growth, \textit{Comm. Part. Diff. Equa.}, \textbf{33}(2008), pp. 636--655.
\bibitem{Cui5} S. Cui and A. Friedman, Analysis of a mathematical model of the growth of necrotic tumors, {\em J. Math. Anal. Appl.},
  \textbf{255}(2001), pp. 636--677.
\bibitem{Pra} G. Da Prato and A. Lunardi, Stability, instability and center manifold theorem for fully nonlinear autonomous parabolic equations in
  Banach spaces, \textit{Arch. Rat. Mech. Anal.}, \textbf{101}(1988), pp. 115--141.
\bibitem{Dra} A. K. Drangeid, The principle of linearized stability for quasilinear parabolic evolution equations, \textit{Nonlinear Anal. T. M. A.},
  \textbf{13}(1989), pp. 1091--1113.
\bibitem{EngNag} K. J. Engel and R. Nagel, \textit{One-Parameter Semigroups for Linear Evolution Equations}, New York: Springer, 2000.
\bibitem{Esc} J. Escher, The Dirichlet-Neumann operator on continuous functions, \textit{Ann. Scuola Norm. Super. Pisa}, \text{21}(1994),
  pp. 235--266.
\bibitem{EscPro} J. Escher and G. Prokert, Analyticity of solutions to nonlinear parabolic equations on manifolds and an application to Stokes flow,
  \textit{J. Math. Fluid Mech.}, \textbf{8}(2006), 1--35.
\bibitem{EscSim2} J. Escher and G. Simonett, Classical solutions for Hele-Shaw models with surface tension, \textit{Adv. Diff. Equa.}, \textbf{2}(1997),
  pp. 619--642.
\bibitem{EscSim3} J. Escher and G. Simonett, Classical solutions for multidimensional Hele-Shaw models, \textit{SIAM J. Math. Anal.}, \textbf{28}(1997),
  pp. 1028--1047.
\bibitem{EscSim4} J. Escher and G. Simonett, The volume preserving mean curvature flow near spheres, \textit{Proc. Amer. Math. Soc.}, \textbf{126}(1998),
  pp. 2789--2796.
\bibitem{EscSim5} J. Escher and G. Simonett, A center manifold analysis for the Mullins-Sekerka model, \textit{J. Diff. Equa.}, \textbf{143}(1998),
  pp. 267--292.
\bibitem{Fried} A. Friedman, \textit{Variational Principles and Free Boundary Problems}, Malabar: Robert E. Krieger Publishing Company Inc., 1982.
\bibitem{FH1} A. Friedman and B. Hu, Asymptotic stability for a free boundary problem arising in a tumor model, \textit{J. Diff. Equa.},
  \textbf{227}(2006), pp. 598--639.
\bibitem{FH2} A. Friedman and B. Hu, Stability and instability of Liapunov-Schmidt and Hopf bifurcation for a free boundary problem arising in a
  tumor model, \textit{Trans. Amer. Math. Soc.}, \textbf{360}(2008), pp. 5291--5342.
\bibitem{FriRei2}  A. Friedman and F. Reitich, Symmetric-breaking bifurcation of analytic solutions to free boundary problems, \textit{Trans.
  Amer. Math. Soc.}, \textbf{353}(2000), pp. 1587--1634.
\bibitem{Gri} D. Grieser, The plasmonic eigenvalue problem, \textit{Reviews in Math. Phys.}, \textbf{26}(2014), pp. 1450005-1--1450005-26.
\bibitem{GunPro} M. G\"{u}nther and G. Prokert, Existence result for the quasistationary motion of a free capillary liquid drop, \textit{Z. Anal.
  Anwendungen}, \textbf{16}(1997), pp. 311--348.
\bibitem{Ham1} R. S. Hamilton, The inverse function theorem of Nash and Moser, \textit{Bull. Amer. Math. Soc.}, \textbf{7}(1982), pp. 65--222.
\bibitem{Hao1} W. Hao, J. D. Hauenstein, B. Hu and et al, Bifurcation for a free boundary problem modeling the growth of a tumor with a necrotic
  core, \textit{Nonliear Anal. Real World Appl.}, \textbf{13}(2012), pp. 694--709.
\bibitem{Lun1} A. Lunardi, Asymptotic exponential stability in quasilinear parabolic equations, \textit{Nonlinear Anal. T. M. A.}, \textbf{9}(1985),
  pp. 563--586.
\bibitem{Lun2} A. Lunardi, \textit{Analytic Semigroups and Optimal Regularity in Parabolic Problems}, Basel: Birkh\"{a}user, 1995.
\bibitem{Poi} M. Poitier-Ferry, The linearization principle for the stability of solutions of quasilinear parabolic equations I,
  \textit{Arch. Rat. Mech. Anal.}, \textbf{77}(1981), pp. 301--320.
\bibitem{PruSim} J. Pr\"{u}ss and G. Simonett, On the manifold of closed hypersurfaces in ${\mathbf{R}}^n$, \textit{Disc. Cont. Dyna. Syst. Ser. A},
  \textbf{33}(2013), pp. 5407--5428.
\bibitem{Sol1} V. A. Solonnikov, Lectures on evolution free boundary problems: classical solutions, \textit{Mathematical Aspects of Evolving
  Interfaces}, J. M. Morel, F. Takens and B. Teissier eds., Lecture Notes in Math., vol. 1812, Berlin: Springer, 2003, pp. 123--175.
\bibitem{Sol2} V. A. Solonnikov, On the stability of uniformly rotating viscous incompressible self-gravitating liquid, \textit{Lect. on Anal.
 \& Nonl. Part. Diff. Equa.}, vol.2 (2010), pp. 225--331.
\bibitem{WuCui} J. Wu and S. Cui, Asymptotic stability of stationary solutions of a free boundary problem modelling the growth of tumours with
  fluid tissues, {\em SIAM J. Math. Anal.}, \textbf{41}(2009), pp. 391--414.

\end{thebibliography}
\end{document}